\documentclass[11pt]{article}

\usepackage{amsmath,amssymb,amsthm}
\usepackage{dsfont}
\usepackage{graphicx}
\usepackage{geometry}
\usepackage[british]{babel}
\usepackage[utf8]{inputenc}
\usepackage{xcolor}
\usepackage{cite}
\usepackage{url}
\usepackage{diagbox}

\def\vs{\lambda_S}
\def\vi{\lambda_I}
\def\vr{\lambda_R}
\def\RR{\mathbb{R}}
\newcommand{\f}{\mathbf{f}}
\newcommand{\g}{\mathbf{g}}

\newcommand{\x}{\mathbf{x}}
\newcommand{\Q}{\mathbf{Q}}
\newcommand{\F}{\mathbf{F}}
\renewcommand{\S}{\mathbf{S}}

\newcommand{\SO}{S_u}
\newcommand{\IO}{I_u}
\newcommand{\EO}{E_u}
\newcommand{\RO}{R_u}
\newcommand{\DD}{D^u}

\newtheorem{remark}{Remark}

\newcommand{\U}{{\mathbf{U}}}

\renewcommand{\u}{{\bf{u}}}
\renewcommand{\v}{{\bf{v}}}

\newcommand{\E}{{\mathbf{E}}}
\newcommand{\D}{{\mathbf{D}}}
\newcommand{\J}{{\mathbf{J}}}

\begin{document}
%\title{Multiscale kinetic transport models and numerical methods for the spatial spread of infectious diseases}
\title{Modeling and simulating the spatial spread of an epidemic through multiscale kinetic transport equations}

\author{Walter Boscheri$^*$\and Giacomo Dimarco$^*$\and Lorenzo Pareschi\footnote{Department of Mathematics and Computer Science, University of Ferrara, Via Machiavelli 30 and Center for Modeling, Computing and Statistic CMCS, University of Ferrara, Via Muratori 9, 44121 - Ferrara, Italy (e-mail: {\tt walter.boscheri@unife.it}, {\tt giacomo.dimarco@unife.it}, {\tt lorenzo.pareschi@unife.it}).}}

\maketitle
\begin{abstract}{In this work we propose a novel space-dependent multiscale model for the spread of infectious diseases in a two-dimensional spatial context on realistic geographical scenarios. The model couples a system of kinetic transport equations describing a population of commuters moving on a large scale (extra-urban) with a system of diffusion equations characterizing the non commuting population acting over a small scale (urban).
The modeling approach permits to avoid unrealistic effects of traditional diffusion models in epidemiology, like infinite propagation speed on large scales and mass migration dynamics. A construction based on the transport formalism of kinetic theory allows to give a clear model interpretation to the interactions between infected and susceptible in compartmental space-dependent models. In addition, in a suitable scaling limit, our approach permits to couple the two populations through a consistent diffusion model acting at the urban scale. A discretization of the system based on finite volumes on unstructured grids, combined with an asymptotic preserving method in time, shows that the model is able to describe correctly the main features of the spatial expansion of an epidemic. An application to the initial spread of COVID-19 is finally presented.}
\end{abstract}


%\tableofcontents

%------------------------------------------------------------------------------------
%------------------------------------------------------------------------------------
\section{Introduction}
The study of epidemic models has certainly experienced enormous growth in recent times due to the impact of the COVID-19 outbreak \cite{APZ, bellomo2020multiscale, Bert, Franco2020, Gatto, FBK2020, PuSa, Tang, Veneziani2021}. Most models and results focused on the nature of epidemic interaction at a global level under the assumption of homogeneous mixing, thus ignoring other possible aspects of human behavior, particularly regarding mobility patterns and social interactions. Although meta-population approaches attempt to overcome such limitations they still suffer from the complexity of the model construction in developing interactions within the meta-populations \cite{BGRCV, CV, Gatto, FWF, REIM}. 

On the contrary, the geographical spread of epidemics is less understood and much less studied than the temporal development and the corresponding control of infectious diseases. The usefulness of realistic models for the spatio-temporal evolution of epidemics is evident if one considers the implementation of appropriate quarantine control strategies in the first phase of the epidemic or the progressive restart of productive activities in a second phase. Both phases require knowledge of the epidemic spread over the territory and therefore, it is of paramount importance to have predictive models that take into account spatial characteristics in order to efficiently deal with the consequences of the epidemic. On the other hand, the increasing availability of information on people's mobility through the use of GPS devices and the Internet, together with the increase in computational capabilities, makes it realistic today to think about the construction of mathematical models able to take into account the heterogeneity of the territory \cite{FMW, FWF, KGV, MWW, Veneziani2020, Veneziani2021}. 

Most models based on partial differential equations describing the spatial dynamics of the epidemic are based on reaction-diffusion equations with a single population per compartment \cite{ABLN, Cap, FMW, KGV, Liu, MWW, Sun, Veneziani2021, Veneziani2020, Wang2020}. These models have highlighted the ability to describe the formation of heterogeneous spatial patterns and the diffusion of the epidemic in geographical contexts where the entire population moves indistinctly. Recently, in order to avoid the paradox of the infinite speed of propagation typical of diffusion problems, alternative models based on hyperbolic equations have been proposed \cite{BCV13, Bert,Colombo}. However, preferential directions of displacement have not been considered while of paramount importance to describe correctly a population dynamic.    

Kinetic transport equations add a whole new level of description to our toolbox of mathematical models for spatial spread of populations. They are situated between individual based models, which act on the microscopic scale and reaction diffusion equations, which rank on the diffusive macroscopic scale. Transport equations are thus often associated with a mesoscopic description based on a statistical physic approach \cite{cs}. These equations use movement characteristics of individual agents (velocity, turning rate etc.), but they describe a population by a continuous density \cite{CMPS, HS, Per} which can be interpreted as the probability for an individual to be in a given position and to move in a certain direction at a given instant of time. They have the possibility to describe complex interaction dynamics in a similar way to particle collision dynamics in rarefied gas flows \cite{Albietal, bellomo2020multiscale, Bellomo2000a, Cer, Deli, DPTZ, PuSa, RS}.

In this work, taking inspiration from the Boltzmann theory of rarefied gas dynamic, we propose a novel multiscale kinetic model for the spread of infectious diseases in a two-dimensional setting. The model is characterized by a coupled system composed, from one side, of kinetic transport equations that describe a population of commuters moving on a large scale (extra-urban) and a set of diffusion equations that characterize the non commuting population on a small scale (urban). The formalism of kinetic theory permits to give a precise interpretation to the interactions between infected and susceptible in a compartmental space-dependent setting and, using a suitable scaling limit \cite{LK}, allows to highlight the relationship with existing models based on reaction-diffusion equations. In addition, our model avoids the unphysical feature of infinite propagation speed and thanks to the interaction between the two populations, it also avoids that the whole population in a compartment moves indiscriminately in the full space originating an unrealistic mass migration effect. This latter aspect is of paramount importance if one wants to track the effective motion of the infection between cities. 

Once the model is defined, its solution on a computational domain describing a realistic geographical scenario poses several difficulties, due to the large dimensionality of the system (in addition to space and time the system depends on the additional velocity variable), the irregular shape of the spatial region of interest, and the multiscale nature of the dynamics. A particular care has been then devoted to the design of an effective numerical solver which provides an accurate and computationally feasible model solution. To this aim  
we introduced a discretization of the system based on Gaussian quadrature points in velocity space \cite{GJL, JPT} and a finite volume approach on unstructured grids \cite{Dumbser2007693, ArepoTN}, combined with an asymptotic preserving method in time \cite{Bos1, Bos2}.

The rest of the manuscript is organized as follows. In Section 2, we introduce the basic features of our model in the case of a simple SIR compartmental interaction. In addition, at an appropriate scaling limit, we show that a two-population diffusion model can be recovered if needed. The model is subsequently extended to a more realistic SEIR type compartmental structure taking into account specific characteristics of the COVID-19 pandemic. Next, in Section 3 we present several numerical examples aimed at validating the model and its numerical solution. 
The details of the numerical scheme employed and a numerical study of its convergence properties are reported in a separate Appendix. 
%We study to that aim both accuracy of the proposed method and its capability to describe the spatial spread of an epidemic. 
Finally, an application to a realistic geographical scenario describing the spread of COVID-19 is analyzed and discussed. Some conclusions and future research directions are reported at the end of the paper.
% together with an Appendix describing the details of the numerical method designed for in the simulations and a numerical study of its convergence properties.  

%------------------------------------------------------------------------------------
%------------------------------------------------------------------------------------
\section{Multiscale kinetic transport models for epidemic spread}
%\subsection{Discrete velocity transport models}
Let $ \Omega\in \RR^2$ a two-dimensional domain of interest. Suppose that individuals can be separated into two separate populations, a commuter population typically moving over long distances (extra-urban) and a non commuter population moving only in small-scale urban areas. For simplicity, we first illustrate our model in the case of a classic SIR compartmental dynamic and subsequently we will extend our arguments to a more realistic SEIR model designed to take into account specific features of the COVID-19 pandemic. In general, our reasoning can naturally be extended to other, more structured, compartment models.

\subsection{A simple kinetic compartmental model}
We consider the population of commuters at position $x\in\Omega$ moving with velocity directions $v \in \mathbb{S}^1$ and denote by $f_S=f_S(x,v,t)$, $f_I=f_I(x,v,t)$ and $f_R=f_R(x,v,t)$, the respective kinetic densities of susceptible (individuals who may be infected by the disease), infected ( individuals who may transmit the disease) and recovered  (individuals healed or died due to the disease). 
%We refer to this model as multiscale kinetic SIR model (MK-SIR).  
%In more general cases, a fourth group of exposed individuals could be considered (resulting in a SEIR-type model). We point out, that the following arguments can be straightforwardly extended to other, more realistic, compartmental models.
The kinetic distribution of commuters is then given by
\[
f(x,v,t)=f_S(x,v,t)+f_I(x,v,t)+f_R(x,v,t),
\]
and we recover their total density by integration over the velocity space 
\[
\rho(x,t)=\int_{\mathbb{S}^1} f(x,v_*,t)\,dv_*.
\]
This then provides the total density of people regardless of their direction of motion. As a consequence 
\[
S(x,t)=\frac1{2\pi}\int_{\mathbb{S}^1}  f_S(x,v,t)\,dv,\,\, I(x,t)=\frac1{2\pi}\int_{\mathbb{S}^1}  f_I(x,v,t)\,dv,\,\, R(x,t)=\frac1{2\pi}\int_{\mathbb{S}^1}  f_R(x,v,t)\,dv,
\]
with $\rho(x,t)=S(x,t)+I(x,t)+R(x,t)$, 
denoting the density fractions of the commuter population at position $x$ and time $t>0$ that are susceptible, infected and recovered respectively. In this setting, the kinetic densities of the commuters satisfy the transport dynamic equations
\begin{eqnarray}
\nonumber
\frac{\partial f_S}{\partial t} + v_S \cdot \nabla_x f_S &=& -F(f_S, I_T) +\frac1{\tau_S}\left(S-f_S\right)\\
\label{eq:kineticc}
\frac{\partial f_I}{\partial t} + v_I \cdot\nabla_x f_I &=&  F(f_S, I_T)-\gamma f_I+\frac1{\tau_I}\left(I-f_I\right)\\
\nonumber
\frac{\partial f_R}{\partial t} + v_R \cdot\nabla_x f_R &=& \gamma f_I+\frac1{\tau_R}\left(R-f_R\right)
\end{eqnarray}
where we defined the total densities
\[
S_T(x,t)=S(x,t)+\SO(x,t),\,\,\, I_T(x,t)=I(x,t)+\IO(x,t),\,\,\, R_T(x,t)=R(x,t)+\RO(x,t),
\]
and $\SO(x,t)$, $\IO(x,t)$, $\RO(x,t)$ are the density fractions of the non commuter individuals moving only on an urban scale. These are  described accordingly to a diffusion dynamic acting on a local scale
\begin{eqnarray}
\nonumber
\frac{\partial \SO}{\partial t} &=& -F(\SO, I_T) + \nabla_x ({\DD_S}\nabla_x S) \\
\label{eq:diffuse}
\frac{\partial \IO}{\partial t}  &=&  F(\SO, I_T)-\gamma \IO+\nabla_x ({\DD_I}\nabla_x I)\\
\nonumber
\frac{\partial \RO}{\partial t}  &=& \gamma \IO+\nabla_x ({\DD_R}\nabla_x R).
\end{eqnarray}
In the above model, the velocities
$v_S=\vs v$, $v_I=\vi v$, $v_R=\vr v$, $\vs,\vi,\vr \geq 0$ in \eqref{eq:kineticc}, as well as the diffusion coefficients $\DD_S$, $\DD_I$, $\DD_R$ in \eqref{eq:diffuse}, are designed to take into account the heterogeneities of geographical areas, and are thus chosen dependent on the spatial location. Similarly, also the relaxation times $\tau_S$, $\tau_I$ and $\tau_R$ are space dependent. The quantity $\gamma=\gamma(x)$ is the recovery rate of infected, while the transmission of the infection is governed by an incidence function $F(\cdot,I_T)$ modeling the transmission of the disease \cite{HWH00}. We assume local interactions to characterize the general incidence function 
\begin{equation}
F(g,I_T)=\beta \frac{g I_T^p}{1+\kappa I_T},
\label{eq:incf}
\end{equation}
where the classic bilinear case corresponds to $p = 1$, $k=0$, even though it has been observed that an incidence rate that increases more than linearly with respect to the number of infected $I$ can occur under certain circumstances \cite{CS78,BCV13,KM05}. The parameter $\beta=\beta(x)$ characterizes the contact rate, whereas the parameter $\kappa=\kappa(x) > 0$ takes into account social distancing and other control effects which may occur during the progress of the disease \cite{Wang2020,Franco2020}. 
The resulting model \eqref{eq:kineticc}-\eqref{eq:diffuse} will be referred to as multiscale kinetic SIR (MK-SIR) model. Note that, because of the presence of two populations acting at different scale distances, the model allows a more realistic description of the typical commuting dynamic involving only a fraction of the population and distinguishes it from the epidemic process affecting the entire population.

%The function $\gamma=\gamma(x)$ is the recovery rate of infectious and the quantities 
%\[
%S=\frac1{2\pi}\int_{\mathbb{S}^1}  f_S(x,v_*,t)\,dv_*,\quad I=\frac1{2\pi}\int_{\mathbb{S}^1}  f_I(x,v_*,t)\,dv_*,\quad R=\frac1{2\pi}\int_{\mathbb{S}^1}  f_R(x,v_*,t)\,dv_*,
%\]
%are the macroscopic densities.
%
%
%
%and, additionally we consider a background of individuals not moving from a specific site, like a city, a village or another residential area, that we will denote using $\SO(x,t)$, $\IO(x,t)$ and $\RO(x,t)$ and which satisfy locally a conventional SIR model
%\begin{eqnarray}
%\nonumber
%\frac{\partial \SO}{\partial t}  &=& -\beta \SO (I+\IO)\\
%\label{eq:kinetic0}
%\frac{\partial \IO}{\partial t}   &=& \beta \SO (I+\IO)-\gamma \IO\\
%\nonumber
%\frac{\partial \RO}{\partial t}   &=& \gamma \IO.
%\end{eqnarray}

The standard threshold of epidemic models is the well-known reproduction number $R_0$, which defines the average number of secondary infections produced when one infected individual is introduced into a host population in which everyone is susceptible \cite{HWH00}. This number determines when an infection can invade and persist in a new host population. For many deterministic infectious disease models, an infection begins in a fully susceptible population if and only if $R_0 > 1$. 
%In the multiscale kinetic model \eqref{eq:kineticc}-\eqref{eq:diffuse} the reproduction number is characterized locally by the ratio ${\beta(x)}/{\gamma(x)}$. More precisely, 
Assuming no inflow/outflow boundary conditions in $\Omega$, integrating over velocity/space and summing up the second equation in \eqref{eq:kineticc} and \eqref{eq:diffuse} we have
\[
\frac{\partial}{\partial t} \int_{\Omega} I_T(x,t)\,dx =  \int_{\Omega} F(S_T,I_T)\,dx-\int_{\Omega} \gamma(x) I_T(x,t)\,dx \geq 0
\]
when
\begin{equation}
R_0(t)=\frac{\int_{\Omega} F(S_T,I_T)\,dx}{\int_{\Omega} \gamma(x) I_T(x,t)\,dx} \geq 1.
\label{eq:R0}
\end{equation}
The above quantity therefore, defines the basic reproduction number for system \eqref{eq:kineticc} describing the space averaged instantaneous variation of the number of infective individuals at time $t>0$. This definition naturally extends locally by integrating over any subset of the computational domain $\Omega$ if one ignores the boundary flows.

Let us also observe that, under the same no inflow/outflow boundary conditions, integrating in $\Omega$ equations \eqref{eq:kineticc}  and \eqref{eq:diffuse}   yields respectively the conservation of the total populations of commuters and non commuters
\[
\frac{\partial}{\partial t} \int_{\Omega} (S(x,t)+I(x,t)+R(x,t))\,dx =0,\quad 
\frac{\partial}{\partial t} \int_{\Omega} (\SO(x,t)+\IO(x,t)+\RO(x,t))\,dx =0. 
\] 

\subsubsection{Diffusion limit of the commuters dynamic}
In this part we discuss the multiscale nature of the model \eqref{eq:kineticc} in order to elucidate the different population behaviors in urban and nonurban areas and to emphasize the relationships with other classical space dependent epidemic models where typically the entire dynamic has a diffusive nature. 
%For simplicity, we consider the case in which $\lambda_S,\lambda_I,\lambda_R$ are space independent constant. T
To this aim, let us introduce the flux functions
\[
J_S=\frac{\lambda_S}{2\pi} \int_{\mathbb{S}^1}  v f_S(x,v,t)\,dv,\quad J_I=\frac{\lambda_I}{2\pi}\int_{\mathbb{S}^1}  v f_I(x,v,t)\,dv,\quad J_R=\frac{\lambda_R}{2\pi}\int_{\mathbb{S}^1}  v f_R(x,v,t)\,dv.
\]
Then, integrating the system \eqref{eq:kineticc} against $v$, it is straightforward to get the following set of equations for the macroscopic densities of commuters
\begin{eqnarray}
\nonumber
\frac{\partial S}{\partial t} + \nabla_x\cdot J_S &=& -F(S, I_T)\\
\label{eq:density}
\frac{\partial I}{\partial t} + \nabla_x\cdot J_I &=& F(S, I_T) -\gamma I\\
\nonumber
\frac{\partial R}{\partial t} + \nabla_x\cdot J_R &=& \gamma I
\end{eqnarray}
whereas the flux functions satisfy
\begin{eqnarray}
\nonumber
\frac{\partial J_S}{\partial t} +  \frac{\vs^2}{2\pi} \int_{\mathbb{S}^1}  (v\cdot \nabla_x f_S)v\,dv &=& -F(J_S, I_T)-\frac1{\tau_S} J_S\\
\label{eq:flux}
\frac{\partial J_I}{\partial t} +  \frac{\vi^2}{2\pi} \int_{\mathbb{S}^1}  (v\cdot \nabla_x f_I)v\,dv &=& -\frac{\lambda_I}{\lambda_S}F(J_S, I_T) - \gamma J_I-\frac1{\tau_I} J_I\\
\nonumber
\frac{\partial J_R}{\partial t} +  \frac{\vr^2}{2\pi} \int_{\mathbb{S}^1}  (v\cdot \nabla_x f_R)v\,dv &=& -\frac{\lambda_R}{\lambda_I} \gamma J_I-\frac1{\tau_R} J_R.
\end{eqnarray}
Clearly, the above system is not closed because the evolution of the fluxes in \eqref{eq:flux} involves higher order moments of the kinetic densities.
The diffusion limit can be formally recovered, by introducing the space dependent diffusion coefficients 
\begin{equation}
D_S=\frac12\lambda_S^2\tau_S,\quad D_I=\frac12\lambda_I^2\tau_I,\quad D_R=\frac12\lambda_R^2\tau_R,
\label{eq:diffcf}
\end{equation}
and keeping the above quantities fixed while letting the relaxation times $\tau_{S,I,R}$ to zero. We get from the r.h.s. in \eqref{eq:kineticc}
\[
\begin{split}
&f_S=S,\quad f_I=I ,\quad f_R=R,
\end{split}
\]
and consequently from \eqref{eq:flux} we recover Fick's law
\begin{equation}
J_S = -{D_S} \nabla_x S,\quad  J_I = -{D_I}\nabla_x I,\quad J_R = -{D_R} \nabla_x R,
\label{eq:flick}
\end{equation}
since
\[
\int_{\mathbb{S}^1}  (v\cdot \nabla_x S)v\,dv =
\int_{\mathbb{S}^1}  (v\otimes v)\,dv \nabla_x S  = \pi \nabla_x S
\]
and similarly for the other densities. Thus, substituting \eqref{eq:flick} into \eqref{eq:density} we get the diffusion system for the commuters' population \cite{MWW, Sun, Webb}
\begin{eqnarray}
\nonumber
\frac{\partial S}{\partial t} &=&  -F(S, I_T)+\nabla_x ({D_S}\nabla_x S)\\
\label{eq:diff}
\frac{\partial I}{\partial t} &=&  F(S, I_T)-\gamma I+\nabla_x ({D_I}\nabla_x I)\\
\nonumber
\frac{\partial R}{\partial t} &=&  \gamma I+\nabla_x ({D_R}\nabla_x R),
\end{eqnarray}
coupled with system \eqref{eq:diffuse} for the non commuters. Let us observe that the model's capability to account for different regimes, hyperbolic or parabolic, accordingly to the space dependent values $\tau_S$, $\tau_I$, $\tau_R$, makes it suitable for describing the dynamics of populations composed of human beings. Indeed, it is clear that the daily routine is a complex mixing of individuals moving at the scale of a city and individuals moving among different urban centers. 
In this situation, it seems reasonable to avoid, due to the lack of microscopic information and the high complexity, the description of the details of movements within an urban area and to describe this aspect through a diffusion operator. On the other hand, commuters when moving from one city to another follow well established connections for which a hyperbolic setting is certainly the most appropriate approach. Finally, we emphasize that commuters entering an urban area change their regime by adapting to a diffusive dynamic thanks to an appropriate choice of the scaling parameters $\tau_S,\tau_I,\tau_R$.

\subsection{Extension to more structured kinetic compartmental models}
It is clear that, one can consider more general compartmental subdivisions. For example, more realistic models for COVID-19 should take into account the exposed population as well as the asymptomatic fraction of infected (see for example \cite{Gatto, Tang}). In a space-dependent model, however, the increase in the number of compartments has some drawbacks, due to both the increasing computational complexity and the inherent difficulties in the process of identifying the parameters. Here, we describe the extension of the multiscale kinetic modeling presented in the previous section to a more general compartmental structure, even if still sufficiently simple, where the exposed population includes also the asymptomatic one (see \cite{Veneziani2021, Veneziani2020} for a similar approach). We denote the commuter individuals which belong to the newly introduced compartment of exposed by $f_E(x,v,t)$. The kinetic dynamic of the commuters then reads
\begin{eqnarray}
\nonumber
\frac{\partial f_S}{\partial t} + v_S \cdot \nabla_x f_S &=& -F_I(f_S,I_T)-\tilde F_E(f_S,E_T) +\frac1{\tau_S}\left(S-f_S\right)\\
\nonumber
\frac{\partial f_E}{\partial t} + v_E \cdot \nabla_x f_E &=& F_I(f_S,I_T)+\tilde F_E(f_S,E_T) - \tilde a f_E - \tilde \gamma_E  f_E +\frac1{\tau_E}\left(E-f_E\right)\\[-.25cm]
\label{eq:kineticc2}
\\[-.25cm]
\nonumber
\frac{\partial f_I}{\partial t} + v_I \cdot\nabla_x f_I &=&  \tilde a f_E - \gamma_I f_I +\frac1{\tau_I}\left(I-f_I\right)\\
\nonumber
\frac{\partial f_R}{\partial t} + v_R \cdot\nabla_x f_R &=& \tilde \gamma_E  f_E + \gamma_I f_I+\frac1{\tau_R}\left(R-f_R\right).
\end{eqnarray}
In such model, the kinetic compartment $f_E$ is characterized by a fraction $\zeta f_E$ of latently infected, i.e. individuals who are not yet contagious and a fraction $(1-\zeta) f_E$ of contagious asymptomatic individuals. The above model can be formally derived starting from a more general one that takes into account also the asymptomatic compartment as proposed in \cite{Gatto, Tang}. In this context, system \eqref{eq:kineticc2} is obtained by merging the asymptomatic and the exposed compartment in one single category while the parameter $\zeta$ still permits to distinguish the two classes of individuals. The parameter $\tilde a$ is then given, in this setting, by $\tilde a=\sigma \zeta a$ where $a^{-1}$ is the average incubation period and $\sigma$ is the probability that a latently infected individual develops symptoms. Finally, $\tilde\gamma_E=(1-\zeta)\gamma_E$ and $\gamma_I$ are the recovery rates of the fraction of asymptomatic and of the infected with symptoms. The so-called incident functions are given by 
\begin{equation}
\tilde F_E(g,E_T)={\beta_E}\frac{g(1-\zeta) E_T}{1+\kappa_E(1-\zeta)E_T},
\label{eq:incfE}
\end{equation} 
while $F_I$ has the form \eqref{eq:incf} for $p=1$ with contact rate parameters $\beta_I$ and $\kappa_I$. Assuming $\sigma=\zeta=1$, namely absence of asymptomatic individuals, we recover the corresponding kinetic version of the standard SEIR model \cite{HWH00}. 

The above system can then be coupled with an analogous dynamics of the non commuter population acting at a diffusion level on restricted regions identified by the urban centers. This reads
\begin{eqnarray}
\nonumber
\frac{\partial \SO}{\partial t} &=& -F_I(\SO, I_T)-\tilde F_E(\SO, E_T) + \nabla_x ({\DD_S}\nabla_x S) \\
\nonumber
\frac{\partial \EO}{\partial t} &=& F_I(\SO, I_T)+\tilde F_E(\SO, E_T) - \tilde a \EO - \tilde \gamma_E \EO + \nabla_x ({\DD_E}\nabla_x E) \\[-.25cm]
\label{eq:diffuse2}
\\[-.25cm]
\nonumber
\frac{\partial \IO}{\partial t}  &=&  \tilde a \EO-\gamma_I \IO+\nabla_x ({\DD_I}\nabla_x I)\\
\nonumber
\frac{\partial \RO}{\partial t}  &=& \tilde\gamma_E E +\gamma_I \IO+\nabla_x ({\DD_R}\nabla_x R).
\end{eqnarray}
In \eqref{eq:kineticc2}-\eqref{eq:diffuse2} we denoted by $E_T = E + \EO$, $E(x,t)=\frac1{2\pi}\int_{\mathbb{S}^1}  f_E(x,v,t)\,dv$, $v_E=\lambda_E v$, $\lambda_E > 0$, $\tau_E$ the relaxation time, and $\DD_E$ the diffusion coefficient associated to the non commuters population of exposed. All these quantities have to be intended spatial dependent. The resulting model \eqref{eq:kineticc2}-\eqref{eq:diffuse2} will be referred to as multiscale kinetic SEIR (MK-SEIR) model in the following.

Let us point out that when a more realistic spatial model is introduced to describe the spread of an epidemic, it is necessary to resort to the definition of generalized reproduction number based on the spectral radius of a suitable epidemiological matrix \cite{Diekmann,Gatto}. To this goal, summing up the second and the third equations in systems \eqref{eq:kineticc2} and 
\eqref{eq:diffuse2} and integrating against $v$ and $x$ we obtain 
\begin{eqnarray}
\nonumber
\frac{\partial}{\partial t} \int_{\Omega} E_T(x,t)\,dx &=&  \int_{\Omega} \left(F_I(S_T, I_T)+\tilde F_E(S_T, E_T)\right)\,dx - \int_{\Omega} \left(\tilde a(x)+\tilde \gamma_E(x)\right) E_T(x,t)\,dx\\[-.25cm]
\label{eq:SEIR}
\\[-.25cm]
\nonumber
\frac{\partial}{\partial t} \int_{\Omega} I_T(x,t)\,dx &=&   \int_{\Omega} \tilde a(x) E_T(x,t)\,dx-\int_{\Omega} \gamma_I(x) I_T(x,t)\,dx,
\end{eqnarray}
showing that both compartments now contribute to the disease transmission. Following the analysis in \cite{Diekmann,Gatto,Veneziani2021}, and omitting the details for brevity, we obtain a reproduction number as the sum of two contributions 
\begin{equation}
\begin{split}
R_0(t) =  &\frac{\int_{\Omega} \tilde F_E(S_T,E_T)\,dx}{\int_{\Omega} \left(\tilde a(x)+\tilde \gamma_E(x)\right)E_T(x,t)\,dx}\\
+&\frac{\int_{\Omega} F_I(S_T,I_T)\,dx}{\int_{\Omega}\left(\tilde a(x)+\tilde \gamma_E(x)\right) E_T(x,t)\,dx}\cdot
\frac{\int_{\Omega}  \tilde a(x)E_T(x,t)\,dx}{\int_{\Omega} \gamma_I(x) I_T(x,t)\,dx}.
\end{split}
\label{eq:R02}
\end{equation}
The above definition can be considered as an indicator of viral reproduction for the velocity and space averaged kinetic model. It is easy to see that, whenever both the exposed as well as the infected population are non decreasing in \eqref{eq:SEIR} we have $R_0(t) \geq 1$, on the other hand the converse is not true in general. In particular, for $\sigma=\zeta=1$ we have $F_E\equiv 0$, $\tilde \gamma_E=0$ and expression \eqref{eq:R02} reduces to \eqref{eq:R0}. Models with additional compartments can be analyzed similarly, we refer to \cite{Gatto, Tang} for more details. 

Finally, in a similar way to the case of the MK-SIR model, taking the limit where the relaxation parameters goes to zero, $\tau_{S,E,I,R}\to 0$, under the assumptions \eqref{eq:diffcf} and the additional one that $D_E=\frac12 \lambda_E^2 \tau_E$ remains fixed, we formally obtain the corresponding diffusion system for the commuting population \cite{SaSiLu,Veneziani2020,Veneziani2021} 
\begin{eqnarray}
\nonumber
\frac{\partial S}{\partial t} &=&  -F_I(S, I_T)-\tilde F_E(S, E_T)+\nabla_x ({D_S}\nabla_x S)\\
\nonumber
\frac{\partial E}{\partial t} &=&  F_I(S, I_T)+\tilde F_E(S, E_T)-\tilde a  E-\tilde \gamma_E  E +\nabla_x ({D_E}\nabla_x E)\\[-.25cm]
\label{eq:diff2}
\\[-.25cm]
\nonumber
\frac{\partial I}{\partial t} &=&  \tilde a  E-\gamma_I I+\nabla_x ({D_I}\nabla_x I)\\
\nonumber
\frac{\partial R}{\partial t} &=&  \tilde \gamma_E  E+ \gamma_I I+\nabla_x ({D_R}\nabla_x R),
\end{eqnarray}
coupled with system \eqref{eq:diffuse2} for the non commuters. Clearly, the latter model also retains the multiscale spatial structure and allows for different situations of urban centers, with intense small-scale movements, connected by major road networks.

%In the sequel we will describe the details of the numerical method used to discretize the above systems. To simplify notations we will describe the schemes in the case of the MK-SIR model, being the extension to the MK-SEIR model and more structured models strightforward.  

\begin{remark}~
In the kinetic transport models describing the dynamics of commuters we assumed the disease transmission independent from the velocity of individuals. On the other hand, if one considers a velocity dependent transmission rate where interactions are not homogeneous with respect to the speed of motion the term $F(g,I_T)$ in \eqref{eq:incf} can be replaced by $F(g,\IO)(x,t)+K(g,f_I)(x,v,t)$ where
\[
K(g,f_I)(x,v,t)=\frac{g(x,v,t)}{2\pi(1+k I_T)}\int_{\mathbb{S}^1} \beta(x,v,v_*) f^p_I(x,v_*,t)\,dv_*.
\]
We refer to \cite{bellomo2020multiscale,Deli,PuSa,RS} for other recent approaches in this direction based on Boltzmann-type equations. Note, however, that this additional level of description can be hardly connected with the few experimental data at disposal. 
\end{remark}

In the next section, we discuss several numerical examples based on the models \eqref{eq:kineticc}-\eqref{eq:diffuse} and \eqref{eq:kineticc2}-\eqref{eq:diffuse2}. It is important to underline that the discretization of the resulting multiscale system of PDEs is not trivial and therefore requires the construction of a specific numerical method able to correctly describe the transition from a convective to a diffusive regime in realistic geometries. Although this is an important aspect of the present contribution, in order to make the presentation more readable all the details concerning the numerical scheme and its validation in terms of accuracy are reported in a separate Appendix.

%------------------------------------------------------------------------------------
%------------------------------------------------------------------------------------

%------------------------------------------------------------------------------------
%\subsection{Asymptotic-preserving schemes}

%------------------------------------------------------------------------------------
%------------------------------------------------------------------------------------
\section{Numerical examples}
In this section, we present several numerical experiments in order to validate the proposed models. First, the kinetic model \eqref{eq:kineticc} is solved in a spatially heterogeneous environment, involving a space-dependent contact rate. Second, the influence of the commuters propagation and of the hyperbolic and parabolic regimes is tested for the the MK-SIR model \eqref{eq:kineticc}-\eqref{eq:diffuse} in a simple setting involving three abstract urban areas and related connections. The last example concerns a realistic application: we simulate the first ten days of the COVID-19 outbreak around March 2020 in the Italian region of Emilia-Romagna. For this test we employ the MK-SEIR model \eqref{eq:kineticc2}-\eqref{eq:diffuse2} and we show its capability to reproduce correctly a complex epidemic scenario. Unstructured meshes composed by triangles or polygonal cells are used to pave the computational domain, while a fully second order in space and time asymptotic preserving discretization is adopted. We refer to Appendix A.1-A.4 for the details and the properties of the numerical scheme used (A.1 for the velocity approximation, A.2 for the finite volume space discretization on unstructured grids, A.3 for the time integration and A.4 for the numerical convergence analysis). In the sequel, if not stated otherwise, in the models we assume that propagation speeds as well as relaxation times are space-dependent but not changing among the populations. Finally, a total number of $M=8$ discrete velocity directions accordingly to a Gaussian quadrature has been used in all simulations (see Appendix A.1 for details).

%------------------------------------------------------------------------------------

%------------------------------------------------------------------------------------
\subsection{Prototype test cases for the MK-SIR model}
\subsubsection{Test 1. Commuter dynamics in epidemic heterogeneous environments}
In this test, we analyse the role played by the epidemic parameters when they are chosen space dependent. To further simplify the dynamic here we ignore the non commuters population and consider only the dynamic of the commuters given by the kinetic system \eqref{eq:kineticc}.

\begin{figure}[tbp]
	\begin{center}
		\begin{tabular}{ccc} 
			\includegraphics[width=0.32\textwidth]{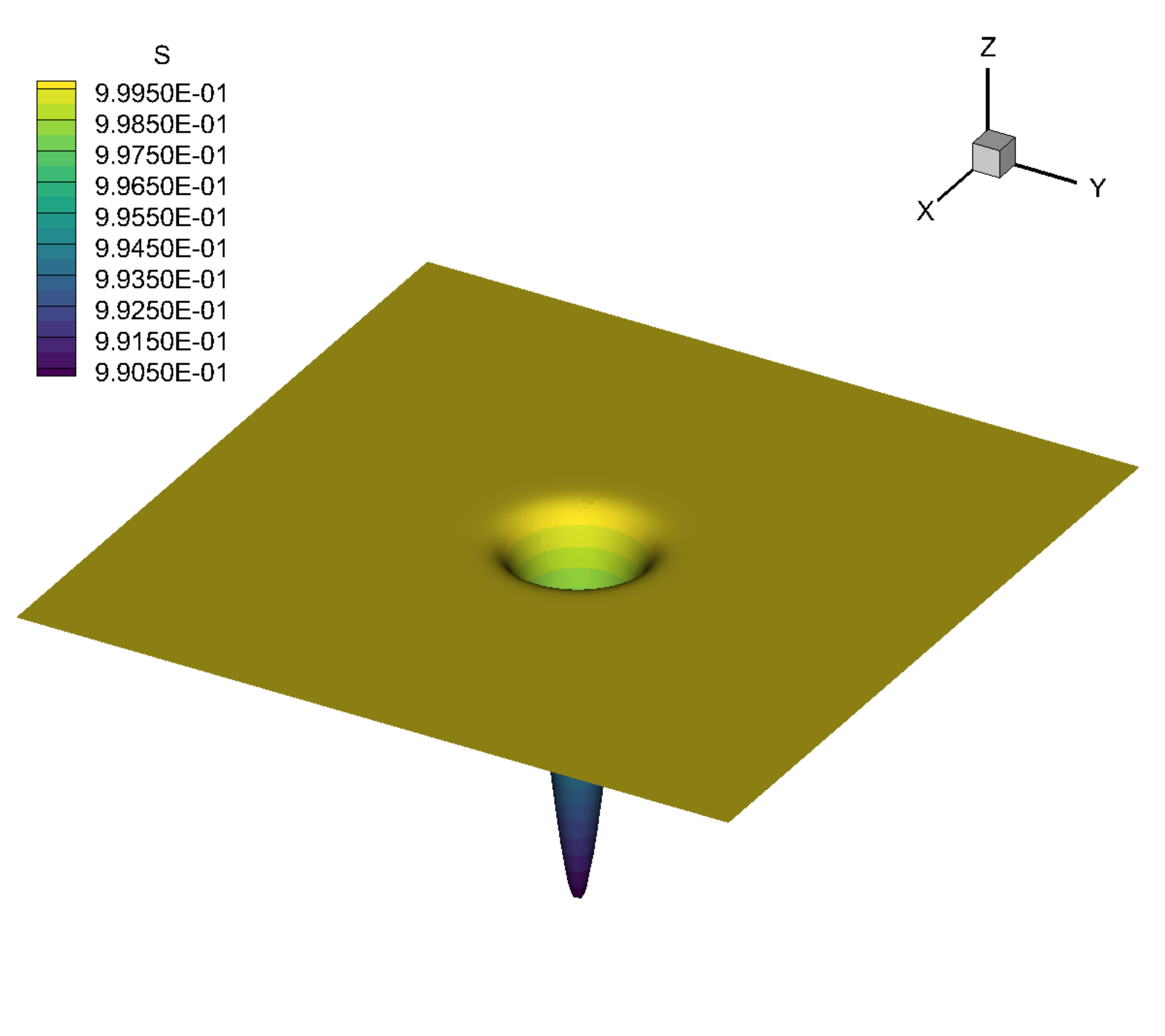} & 
			\includegraphics[width=0.32\textwidth]{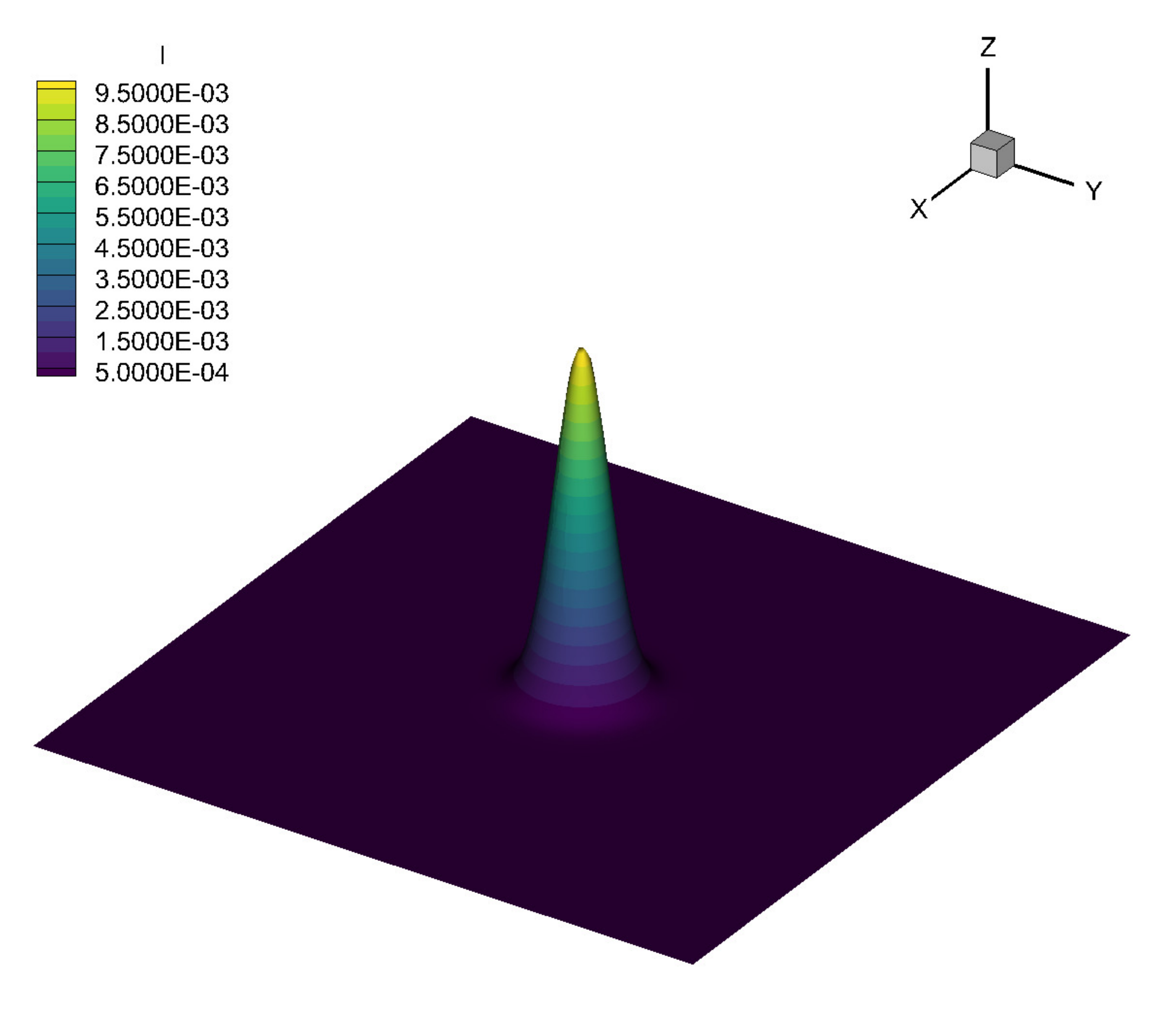} &
			\includegraphics[width=0.32\textwidth]{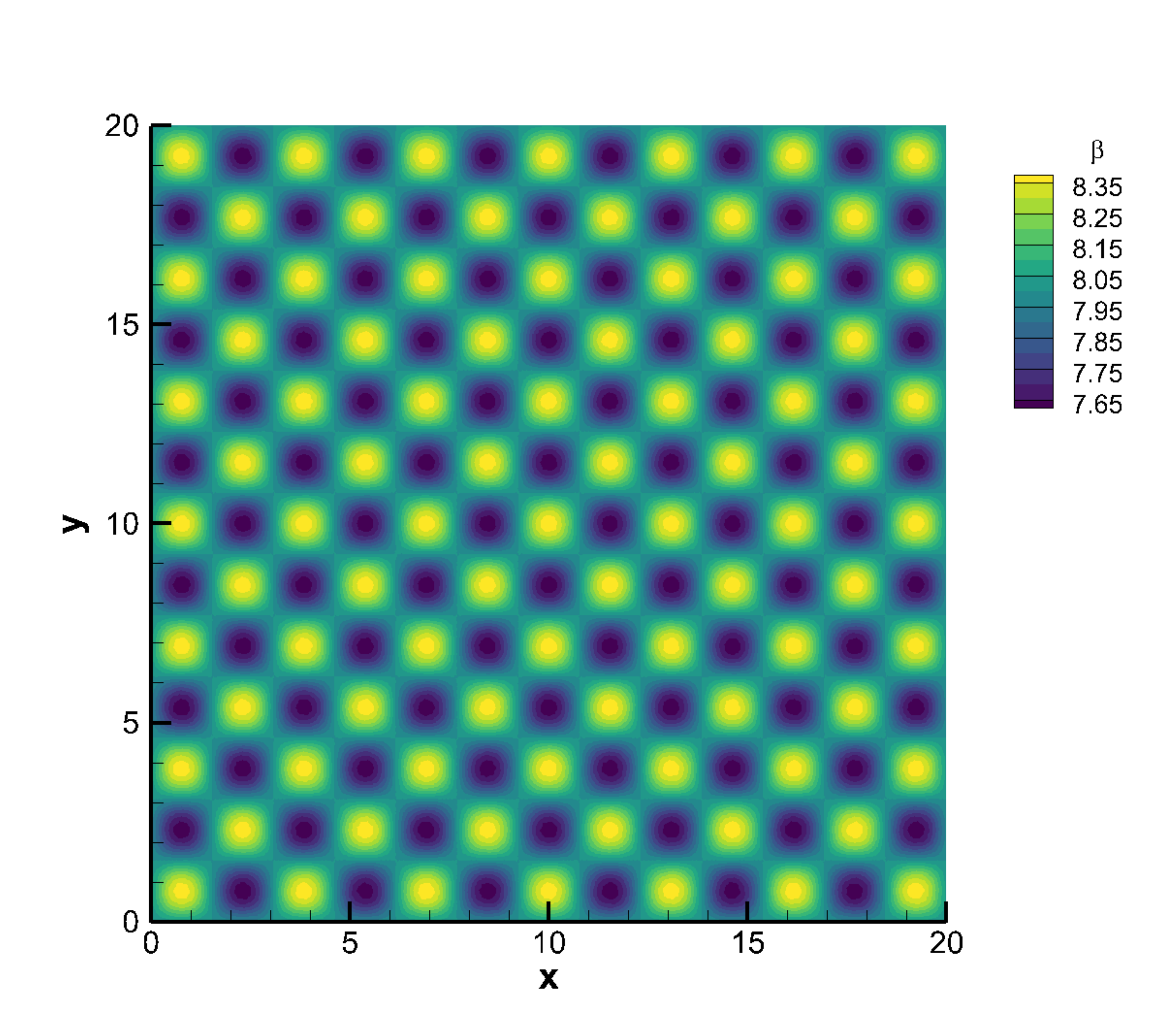}\\
		\end{tabular} 
	\end{center}
		\caption{Test 1. Initial condition for $S$ (left), $I$ (middle) and $\beta$ (right).}
		\label{fig.test2_IC}
\end{figure}

\begin{figure}[tbp]
	\begin{center}
		\begin{tabular}{cc} 
			\includegraphics[width=0.4\textwidth]{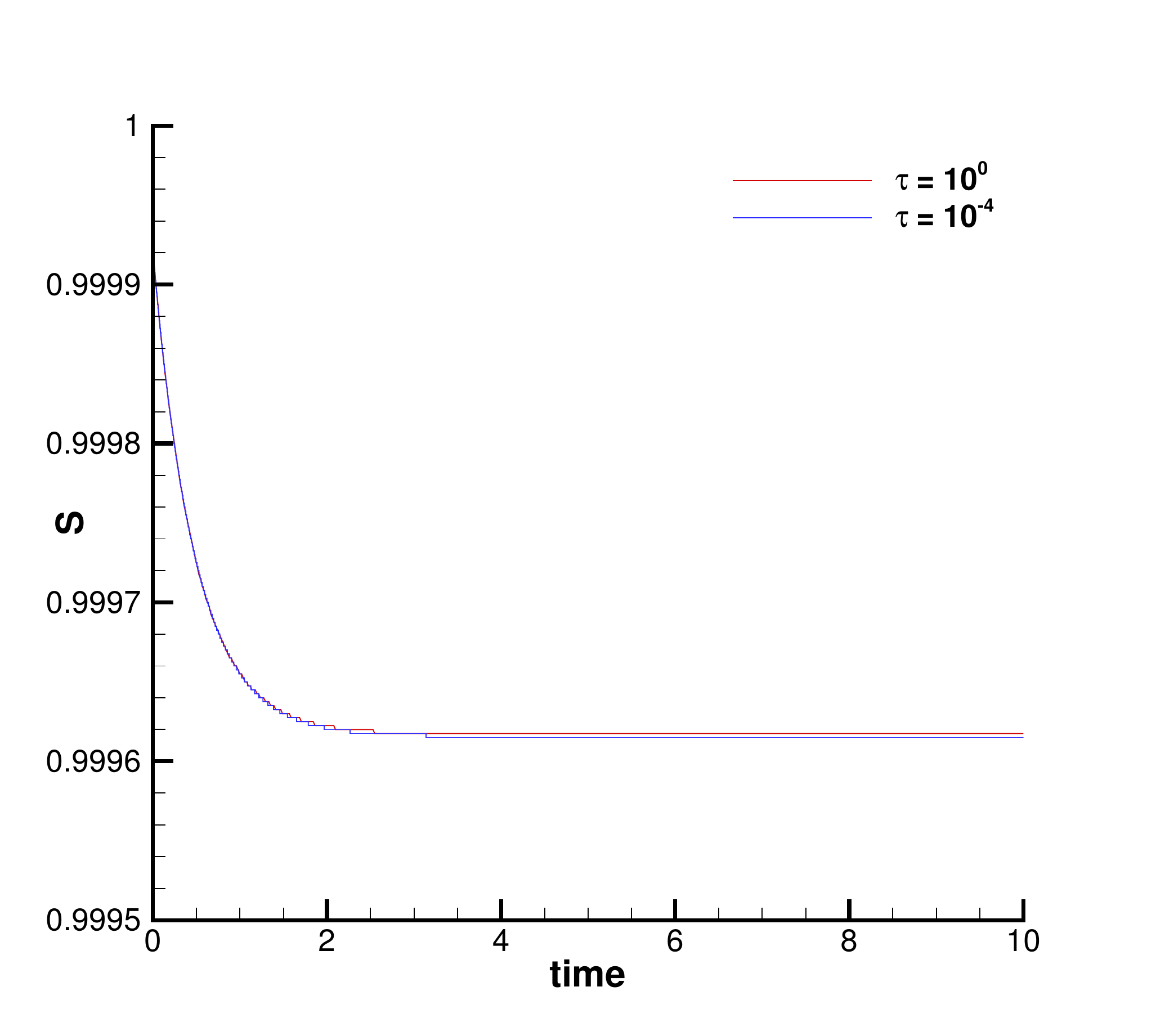} & 
			\includegraphics[width=0.4\textwidth]{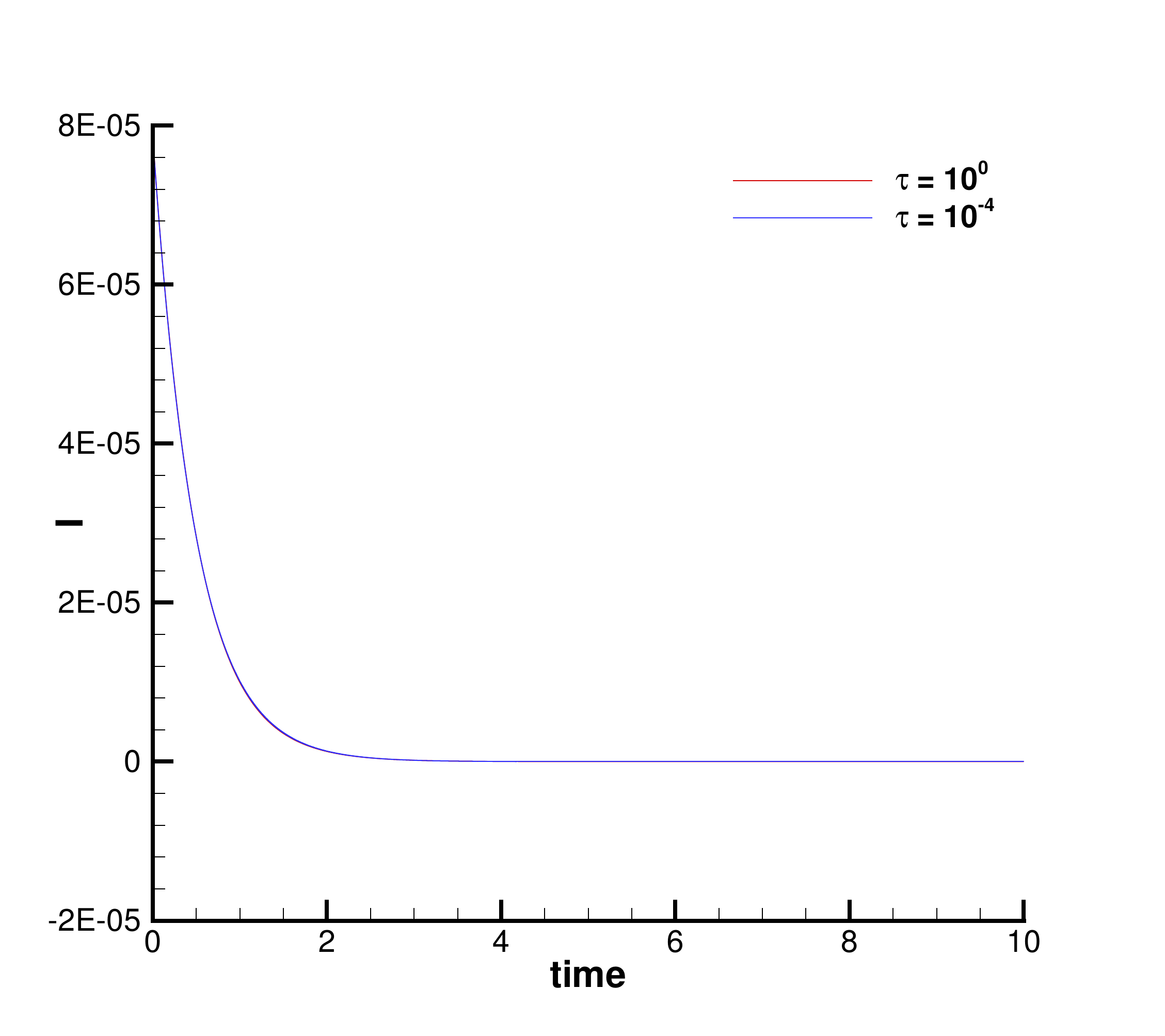} \\
		\end{tabular} 
	\end{center}
		\caption{Test 1. Evolution of $S$ (left) and $I$ (right) for a reproduction number $R_0<1$ and relaxation times $\tau=1.0$ with  $\lambda^2=1.0$ (hyperbolic regime, blue line) and $\tau=10^{-4}$ with $\lambda^2=10^4$ (parabolic regime, red line).}
		\label{fig.test2_R0small_SIR}
\end{figure}

\begin{figure}[tbp]
	\begin{center}
		\begin{tabular}{cc} 
			\includegraphics[width=0.33\textwidth]{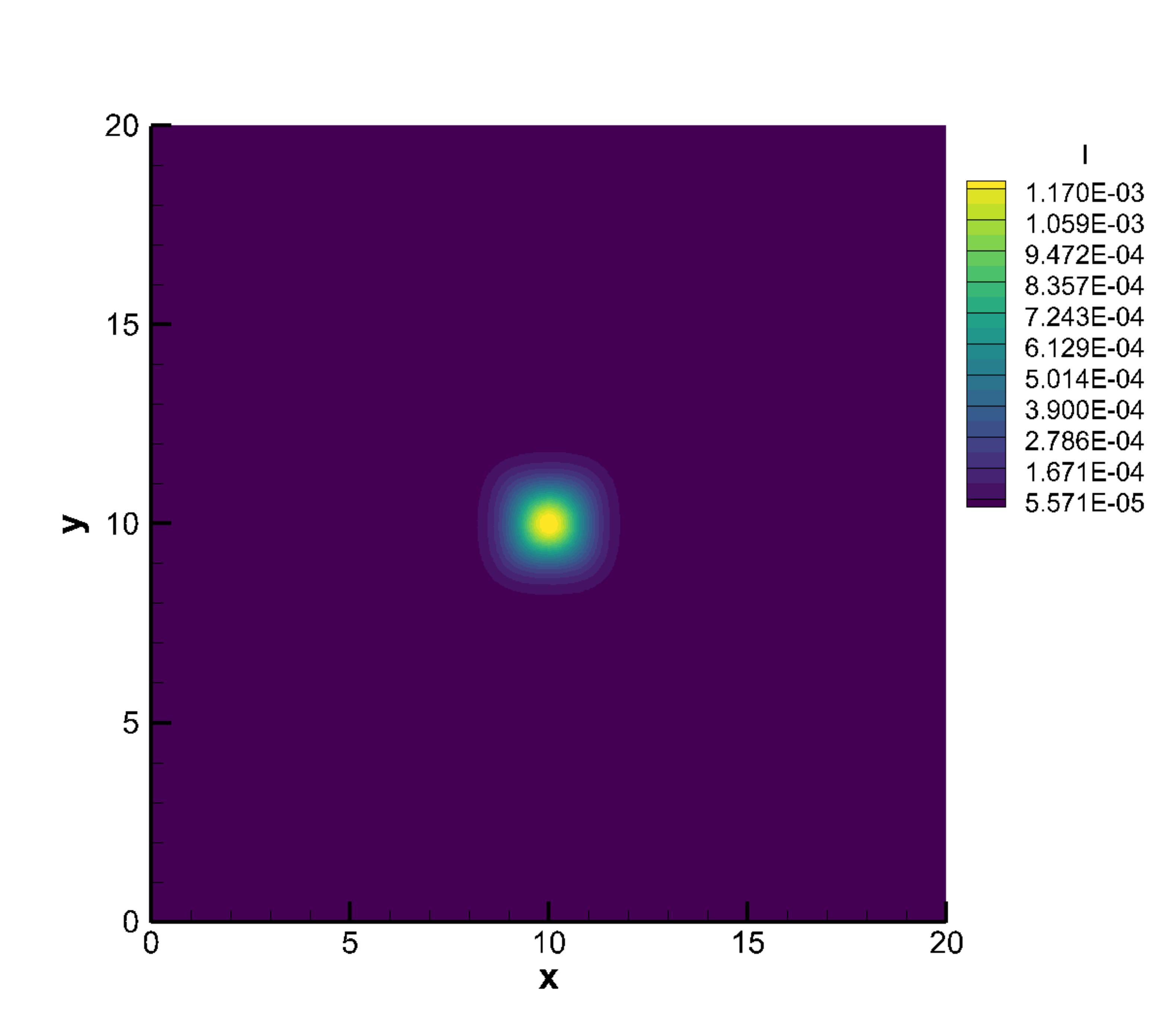} & 
			\includegraphics[width=0.33\textwidth]{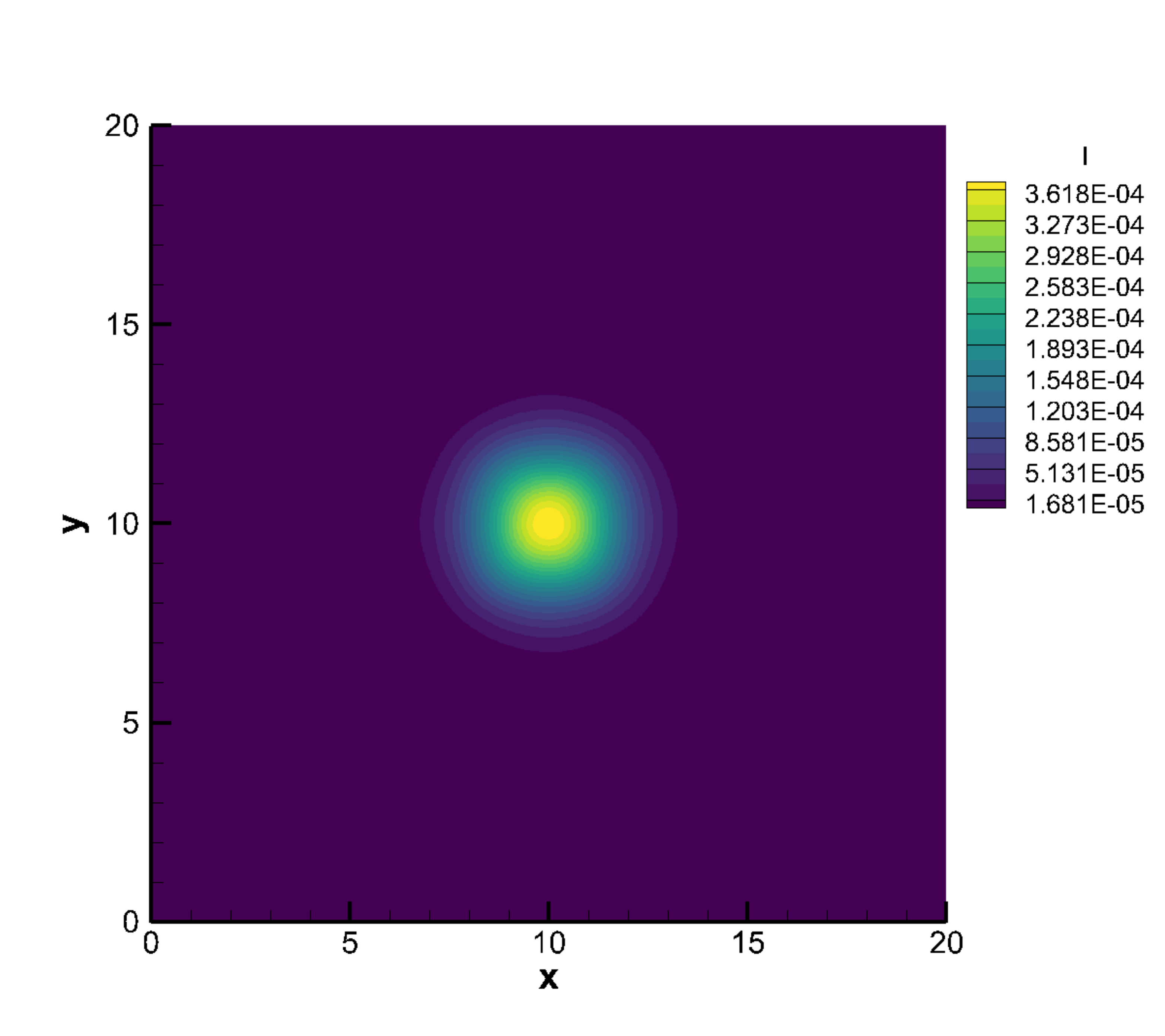} \\
			\includegraphics[width=0.33\textwidth]{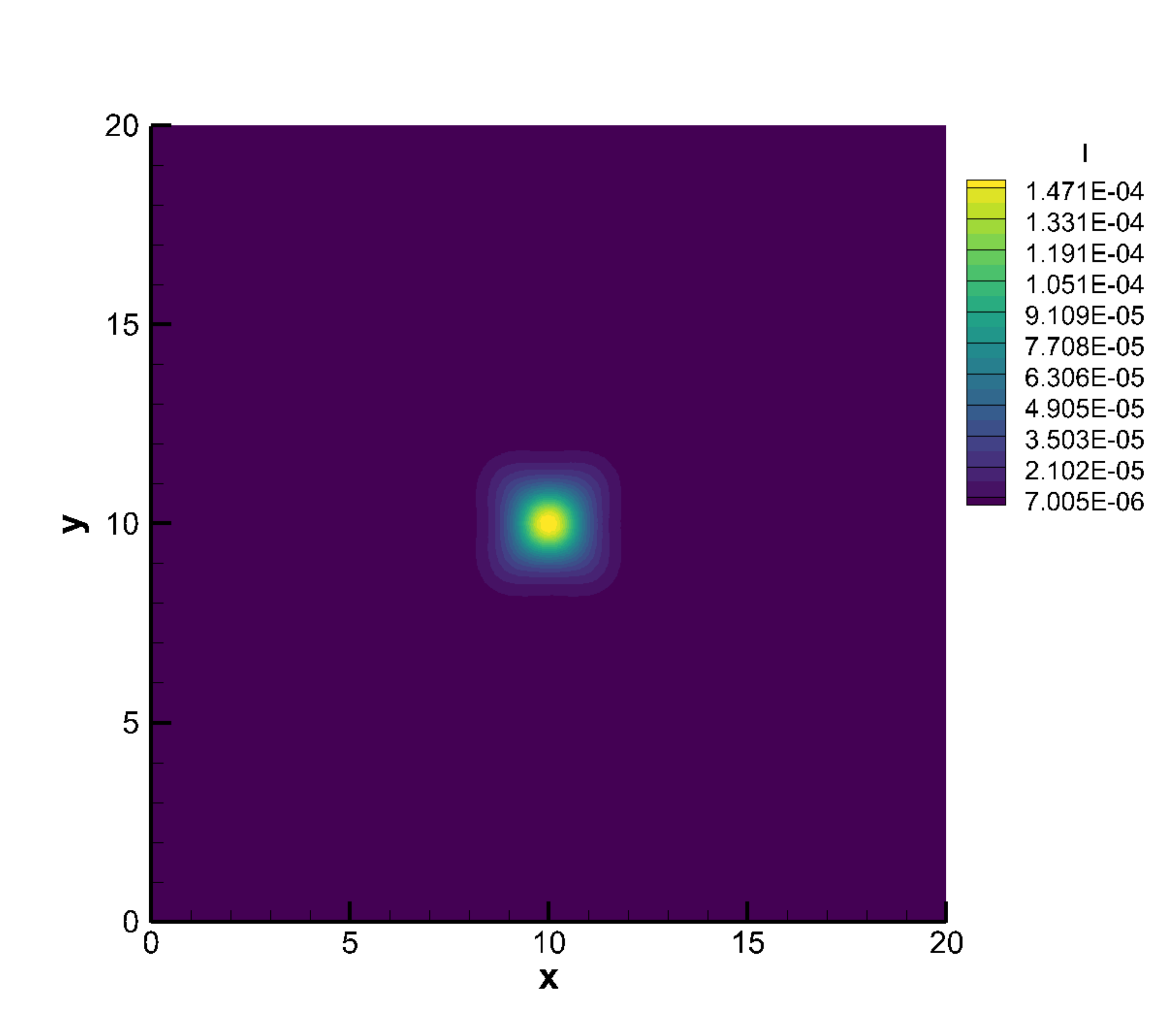} & 
			\includegraphics[width=0.33\textwidth]{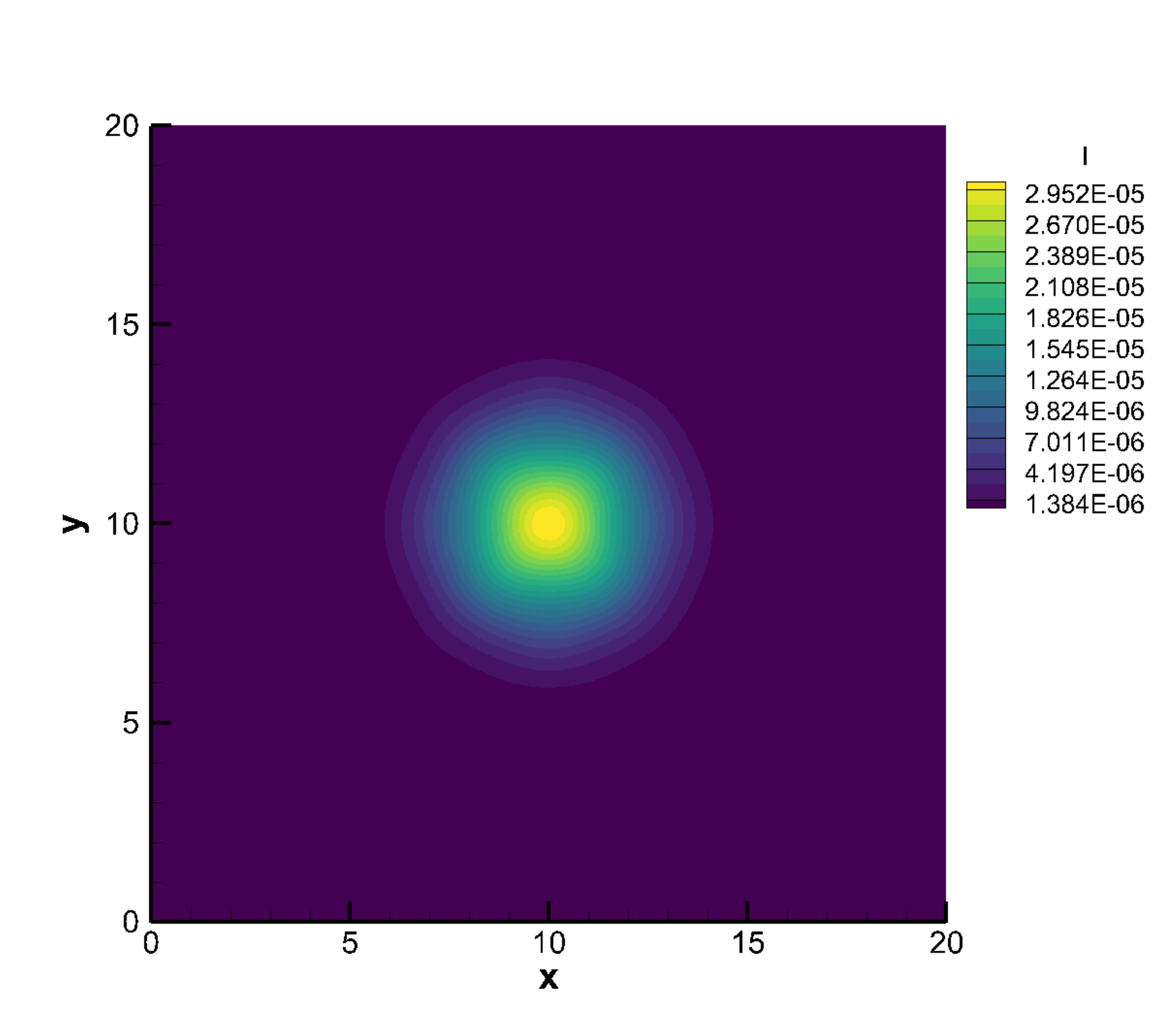} \\
		\end{tabular} 
	\end{center}
		\caption{Test 1. Time evolution of the infected $I$ for a reproduction number $R_0<1$. Left image shows the results for a relaxation time $\tau=1.0$ with  $\lambda^2=1.0$, i.e. a hyperbolic regime. Right image shows the results for $\tau=10^{-4}$ with $\lambda^2=10^4$, i.e. a parabolic regime. Numerical results at output times $t=1$ and $t=2$ (from top to bottom).}
		\label{fig.test2_R0small_SI}
\end{figure}

\begin{figure}[tbp]
	\begin{center}
		\begin{tabular}{cc} 
			\includegraphics[width=0.4\textwidth]{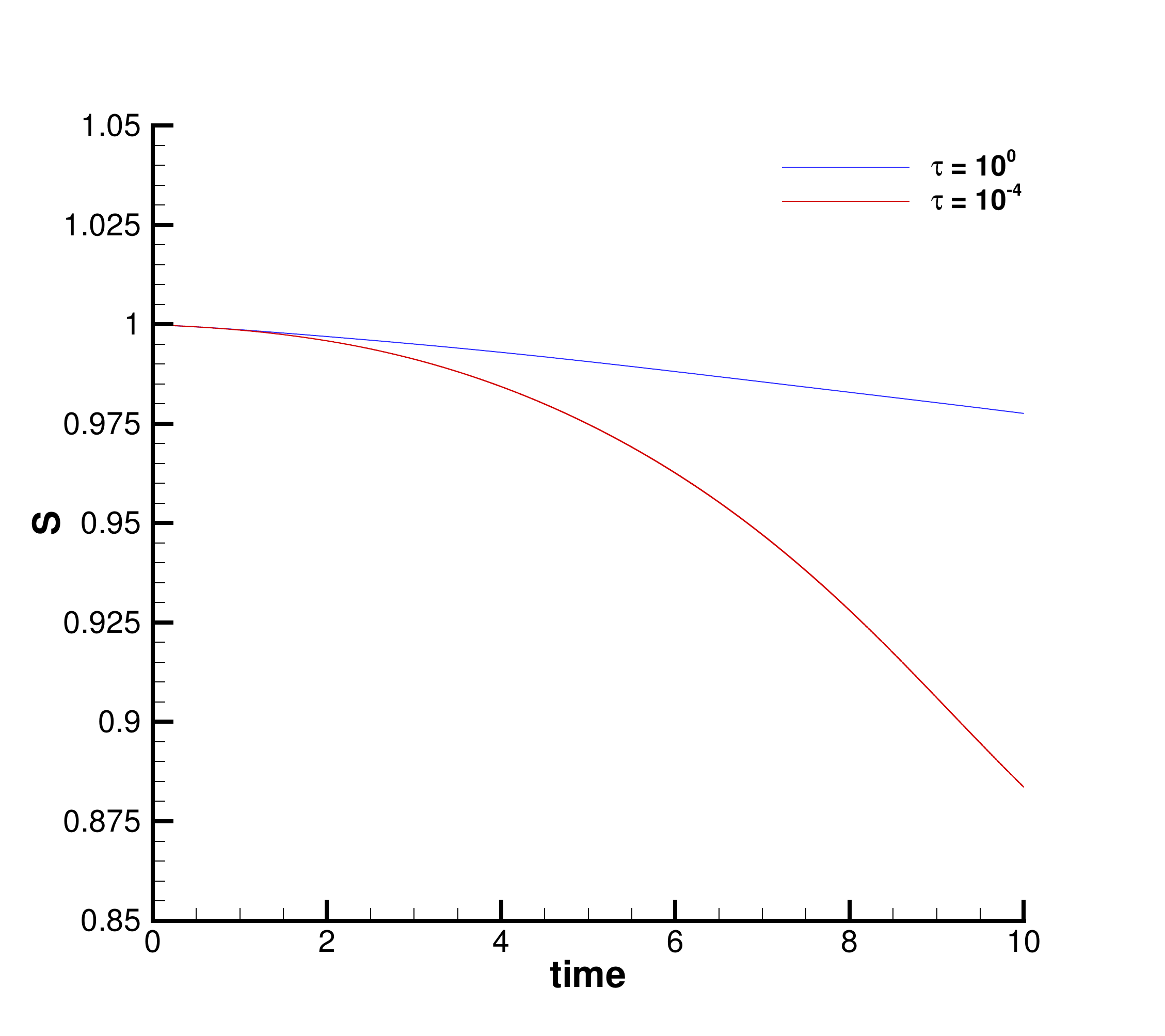} & 
			\includegraphics[width=0.4\textwidth]{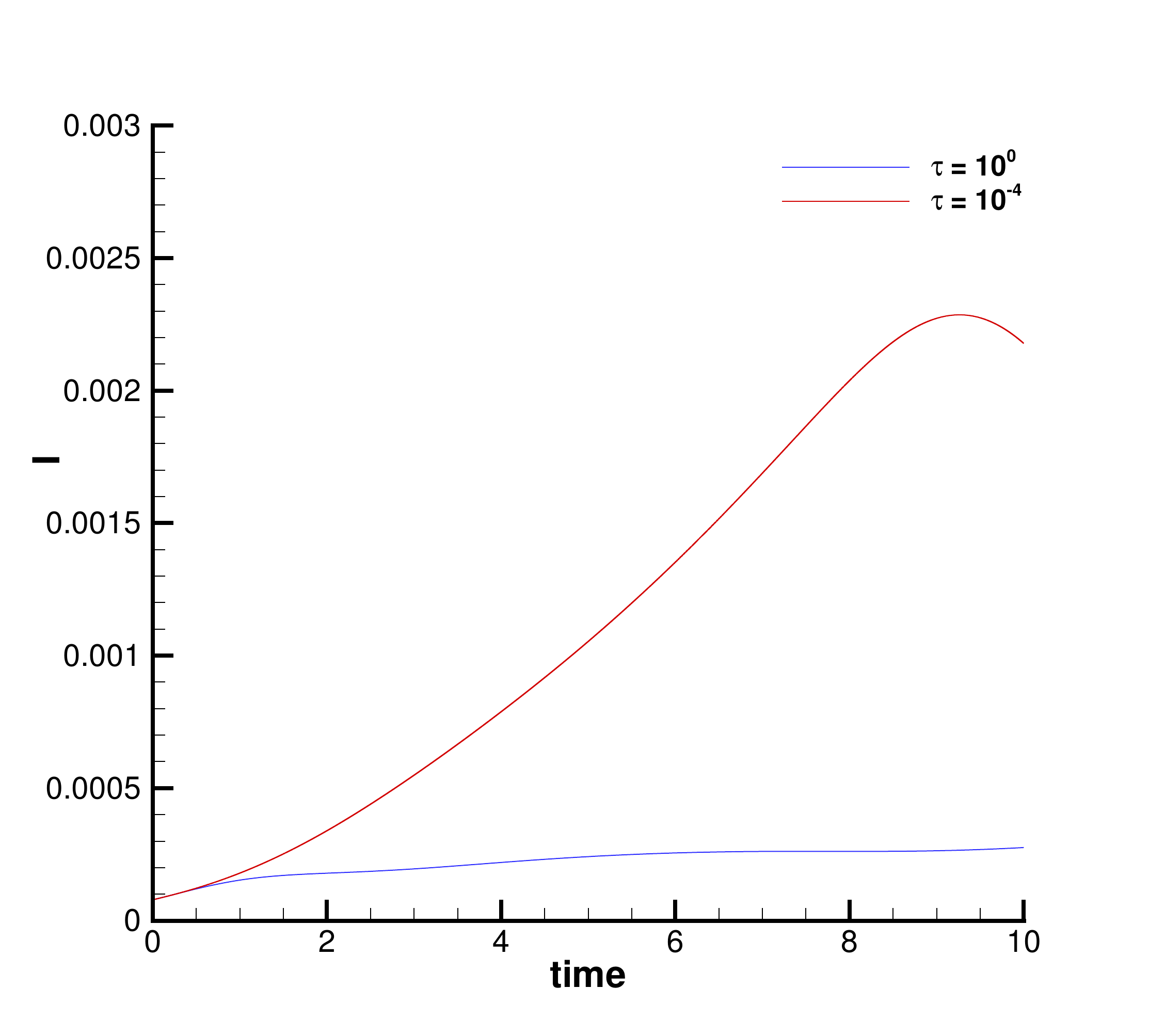} \\
		\end{tabular} 
	\end{center}
		\caption{Test 1. Evolution of $S$ (left) and $I$ (right) for a reproduction number $R_0>1$ and relaxation times $\tau=1.0$ with  $\lambda^2=1.0$ (blue line, hyperbolic regime) and $\tau=10^{-4}$ with $\lambda^2=10^4$ (red line, parabolic regime).}
		\label{fig.test2_R0big_SIR}
\end{figure}

\begin{figure}[!htbp]
	\begin{center}
		\begin{tabular}{cc} 
			\includegraphics[width=0.33\textwidth]{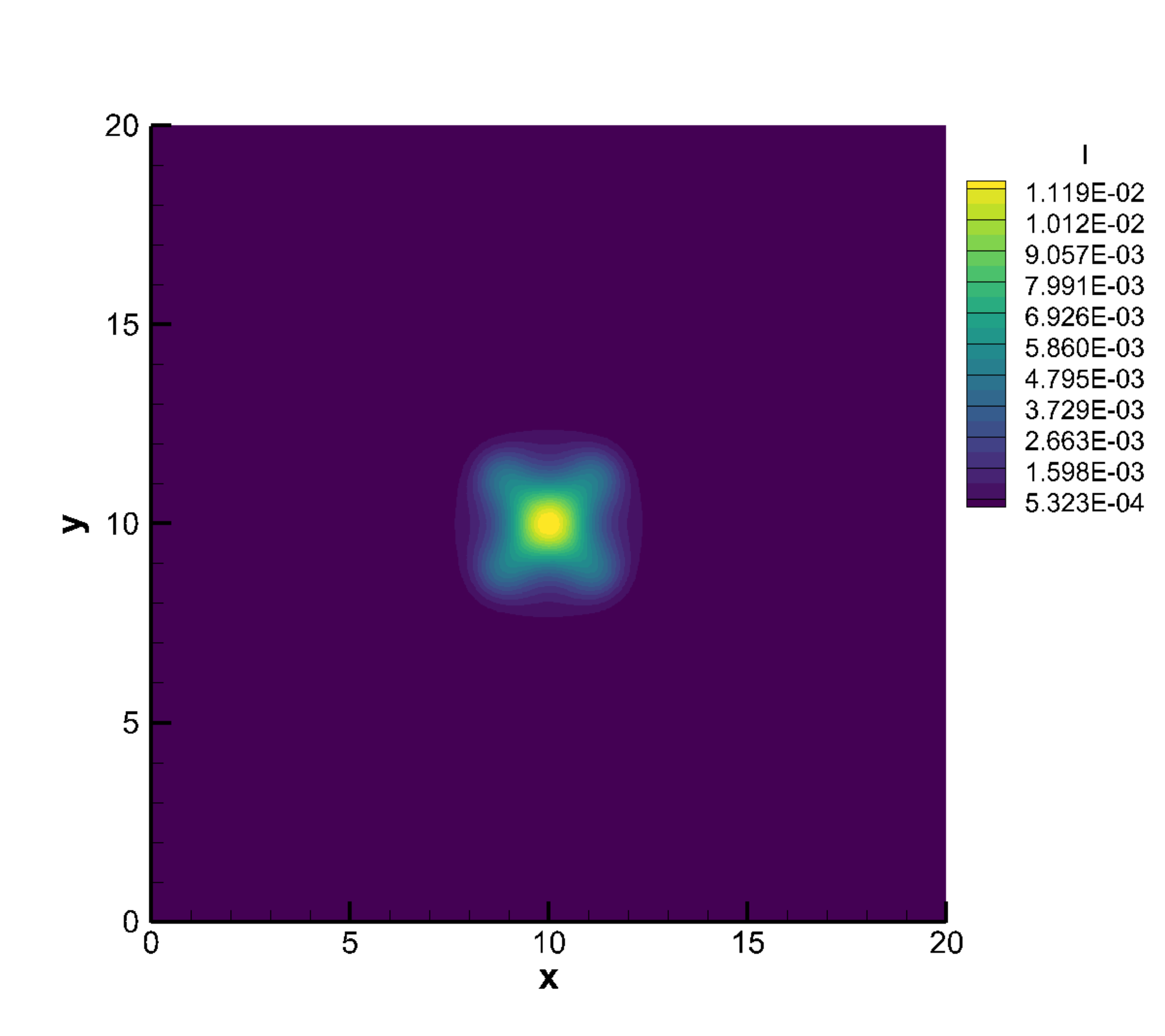} & 
			\includegraphics[width=0.33\textwidth]{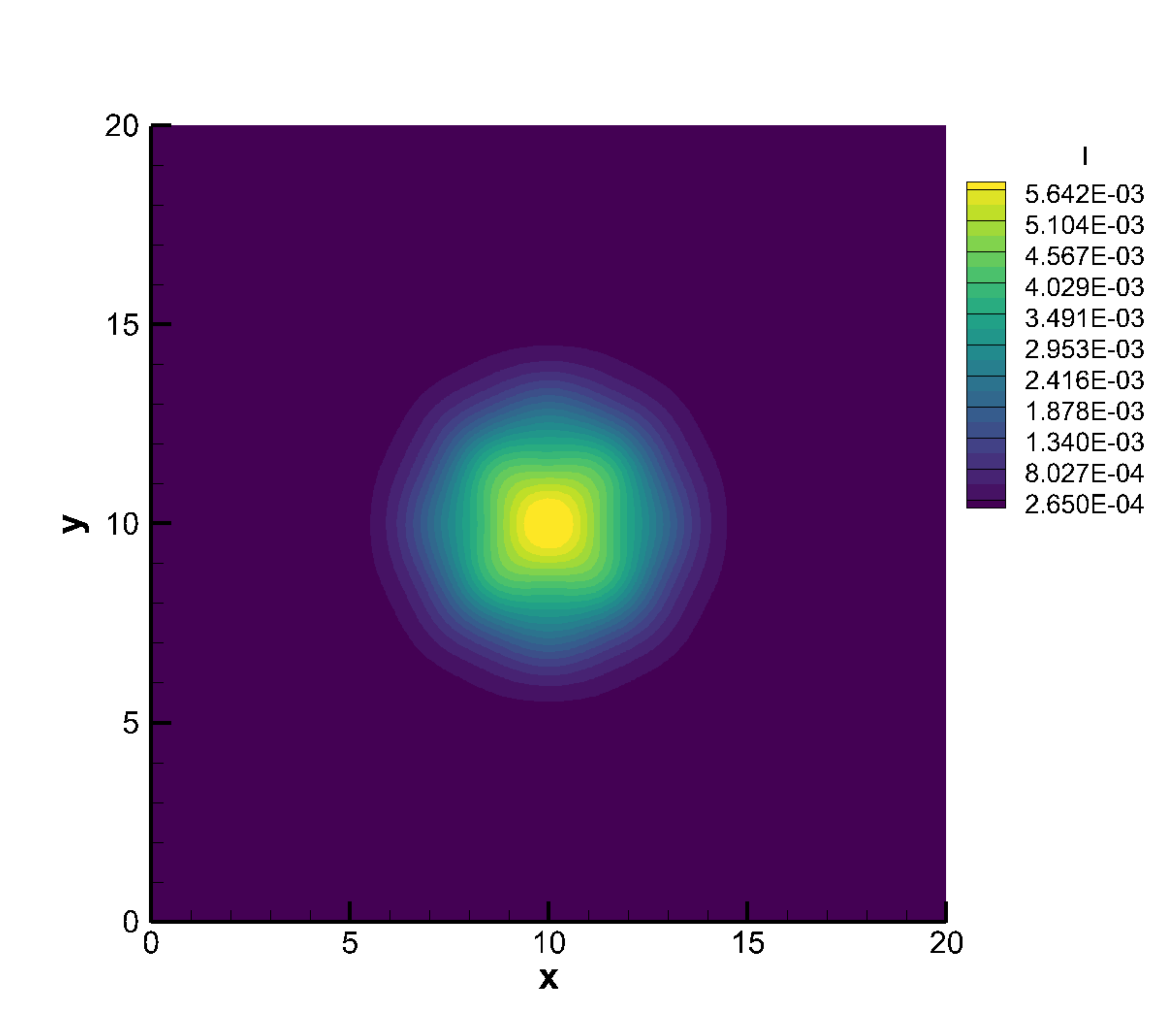} \\
			\includegraphics[width=0.33\textwidth]{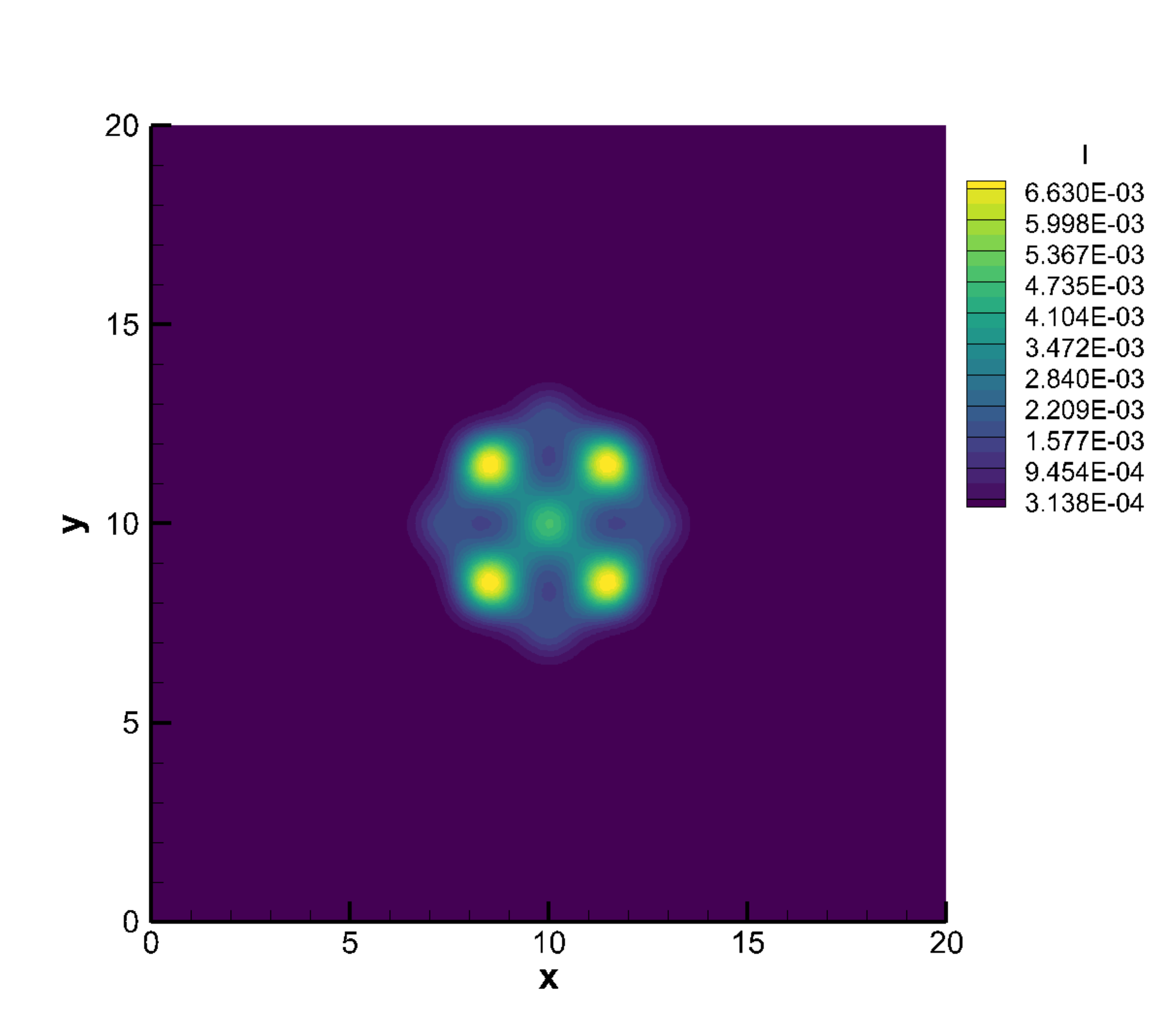} & 
			\includegraphics[width=0.33\textwidth]{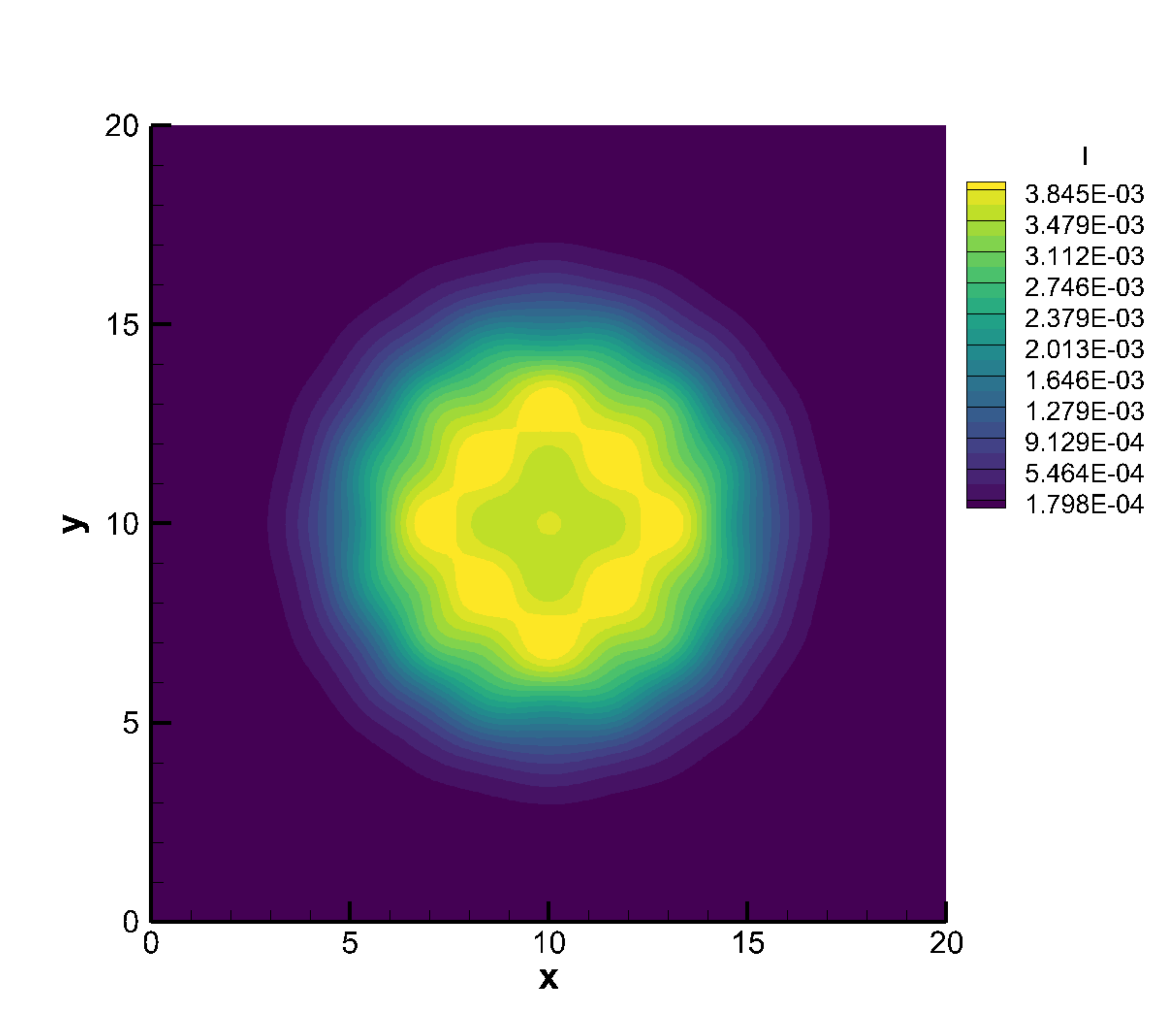} \\
			\includegraphics[width=0.33\textwidth]{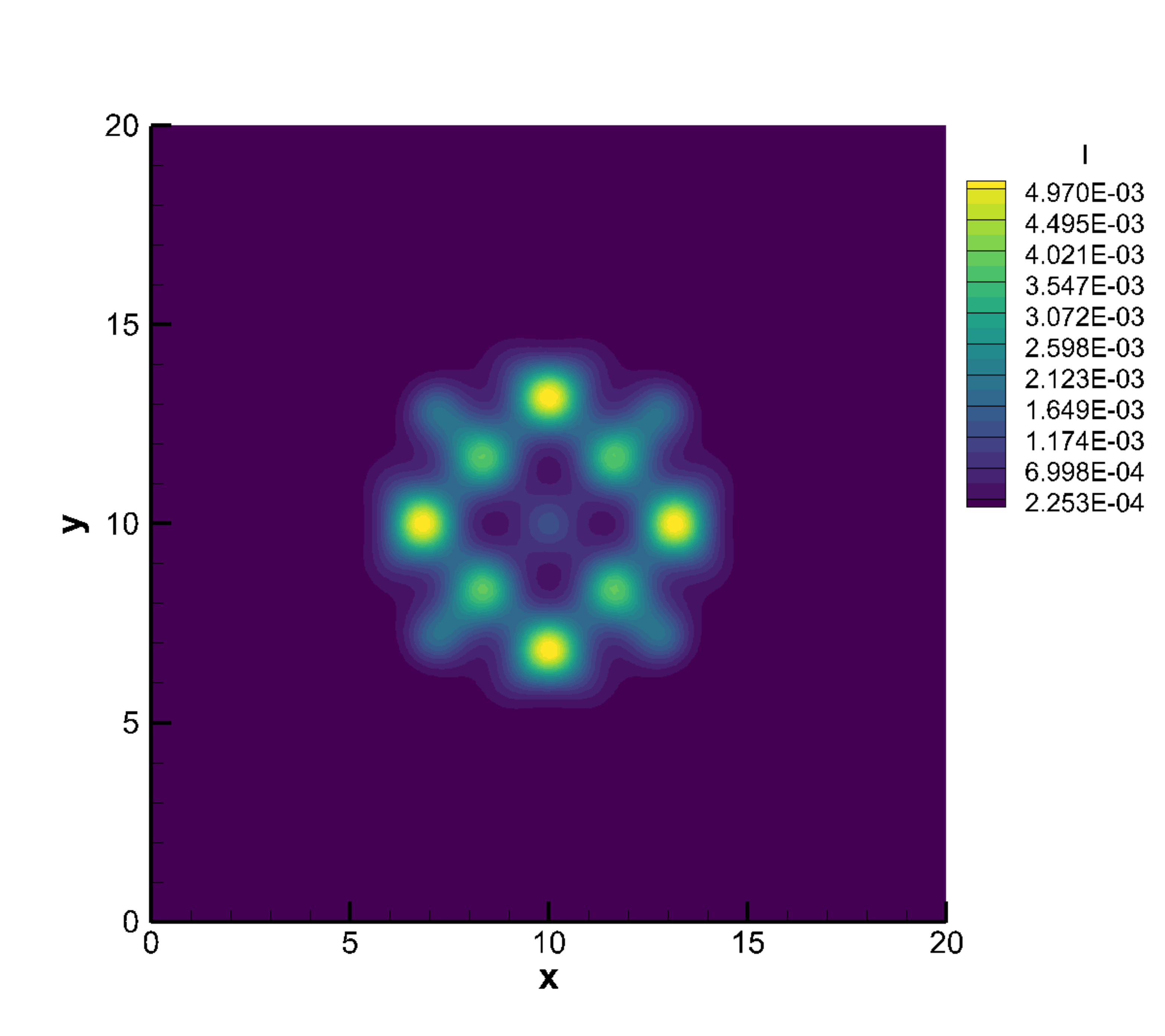} & 
			\includegraphics[width=0.33\textwidth]{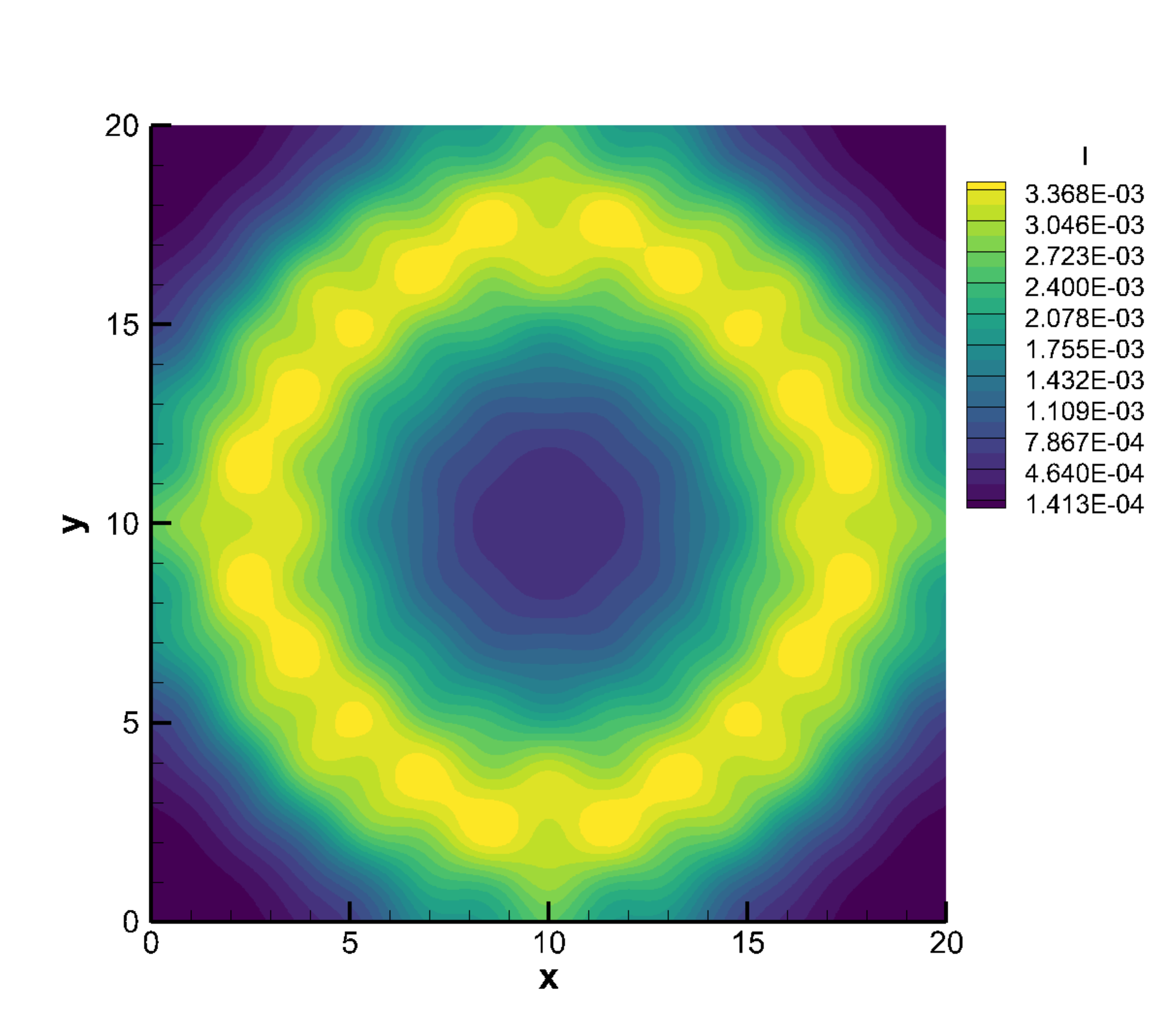} \\
			\includegraphics[width=0.33\textwidth]{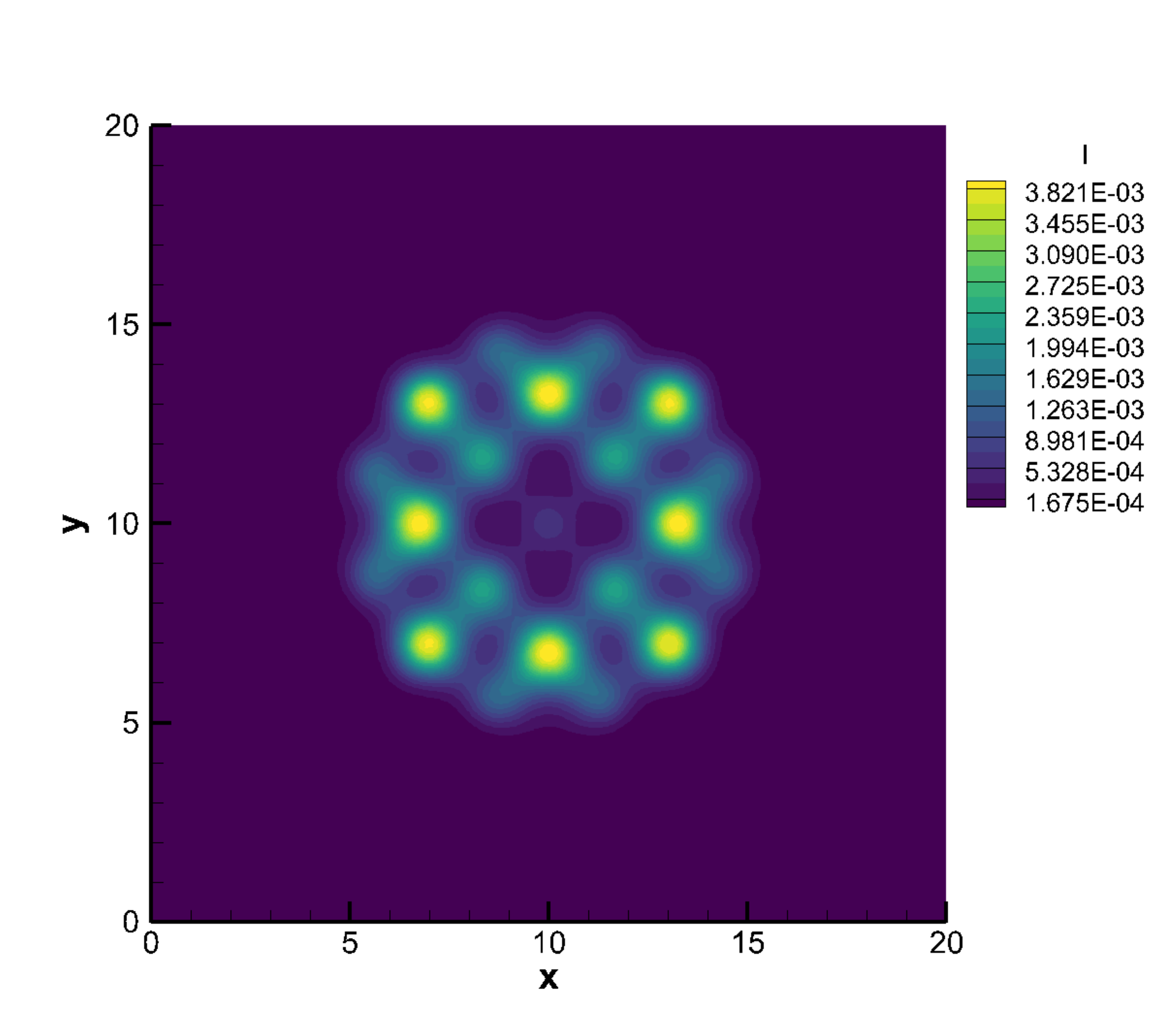} & 
			\includegraphics[width=0.33\textwidth]{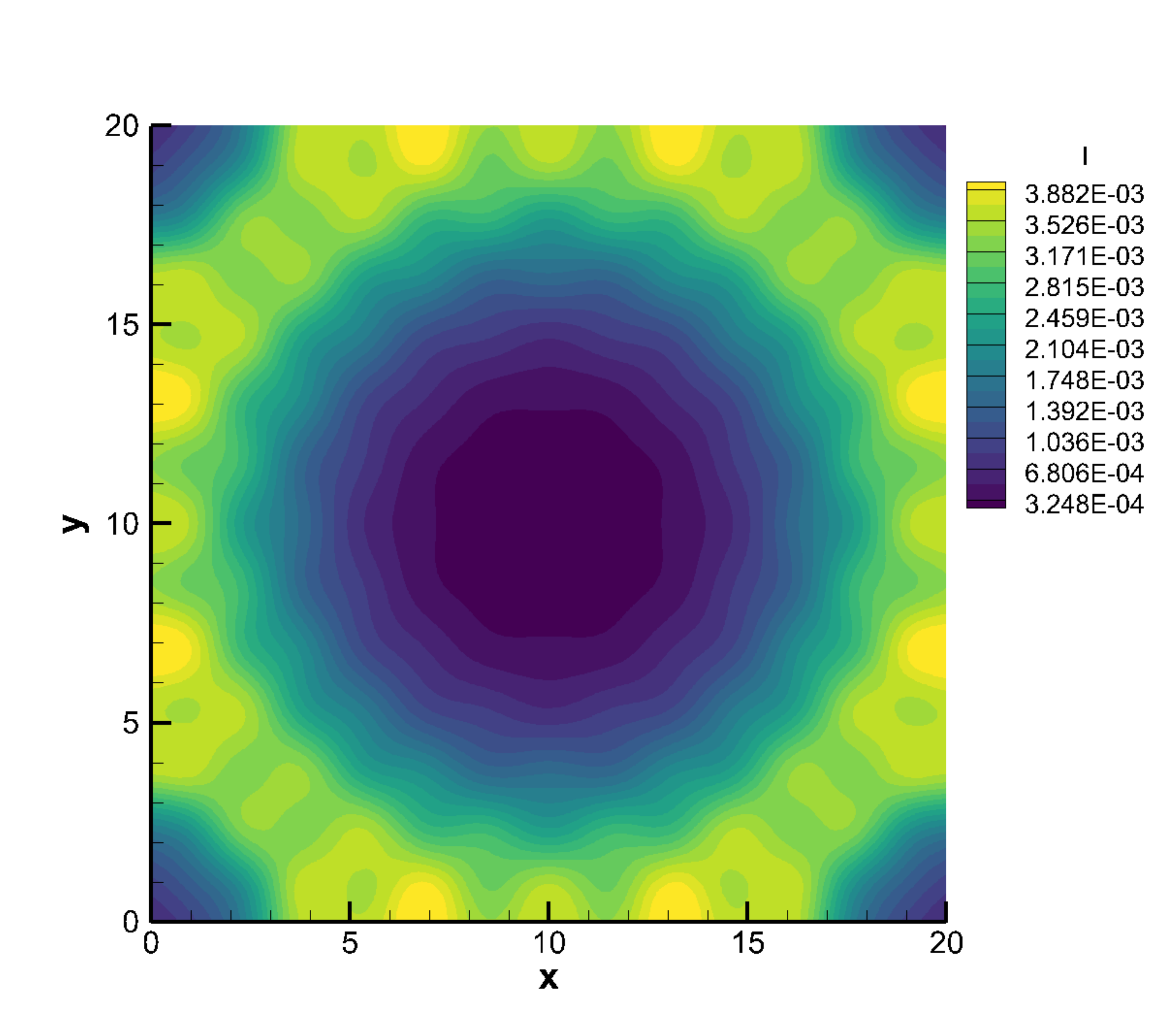} \\
		\end{tabular} 
	\end{center}
			\caption{Test 1. Time evolution of the infected $I$ for a reproduction number $R_0>1$. Left image shows the results for a relaxation time $\tau=1.0$ with  $\lambda^2=1.0$, i.e. a hyperbolic regime. Right image shows the results for $\tau=10^{-4}$ with $\lambda^2=10^4$, i.e. a parabolic regime. Numerical results at output times $t=2$, $t=4$, $t=8$ and $t=10$ (from top to bottom).}
		\label{fig.test2_R0big_I}
\end{figure}

In particular, we focus on the contact rate $\beta=\beta(x)$ and we move along the lines of a test firstly proposed in \cite{Wang2020}. The distribution of the population is initially assigned as follows within the computational domain $\Omega=[0;20]^2$:
\begin{equation}
S=1-I, \qquad I = 0.01 \cdot e^{-(x-10)^2-(y-10)^2}, \qquad R=0.
\end{equation}
Zero-flux boundaries are adopted and the fluxes associated to each population are initially set to zero. The recovery rate is set to $\gamma=10$, whereas the contact rate is given by
\begin{equation}
\beta = \tilde{\beta} \left[ 1 + 0.05 \sin\left( \frac{13\pi x}{20} \right) \sin\left( \frac{13\pi y}{20} \right) \right].
\label{eqn.beta0-test2}
\end{equation}
The initial condition for this test is shown in Figure \ref{fig.test2_IC} for a computational mesh composed of $N_E=15672$ polygonal control volumes. 
By choosing $\tilde{\beta}=8$ an initial reproduction number of $R_0=0.808$ is obtained. This means that the infection is not able to start spreading being $R_0<1$. In contrast to this first choice, we consider a second possibility by fixing $\tilde{\beta}=10$ which leads to $R_0=1.111$. In this second case we expect the infection to expand in the host population and the epidemic disease starts spreading. These two scenarios are considered in two different regimes, namely a hyperbolic configuration with $\tau=1.0$ and $\lambda^2=1.0$ as well as a diffusive setting with $\tau=10^{-4}$ and $\lambda^2=10^4$ to appreciate the differences in modeling the migration of individuals in a given region. The final time of the simulation is $t_f=10$. The time evolution of the %susceptible $S$ and 
infected $I$ population is shown in Figure \ref{fig.test2_R0small_SIR} for both hyperbolic and parabolic configurations with $R_0<1$, while Figure \ref{fig.test2_R0small_SI} depicts the distribution %of $S$ and 
$I$ at different output times. The same kind of plots are found in Figures \ref{fig.test2_R0big_SIR} and
\ref{fig.test2_R0big_I} in the case $R_0>1$, clearly capturing the spreading of the disease. Different dynamics of the infectious spread are noticed according to the hyperbolic or parabolic regime, highlighting the different behavior of the model which is one of the relevant aspects of our approach. This is even more evident when considering the scenario with $R_0>1$: the propagations speeds in the parabolic case are much faster than the ones in the hyperbolic regime, as depicted in Figure \ref{fig.test2_R0big_I}. Finally, let us observe that the patterns generated by the infective population in the hyperbolic limit, correctly resemble the spatial distribution of the contact rate given by \eqref{eqn.beta0-test2}. Contrarily, the parabolic limit leads to a more homogeneous pattern, which tends to become independent of the variable contact rate. 
%This latter is a less physical behavior that we aim to avoid in a realistic setting such as the one shown in the last test case.

% R0<1 --------------------------------------------------

%\begin{figure}[tbp]
%	\begin{center}
%		\begin{tabular}{cc} 
%			\includegraphics[width=0.40\textwidth]{test2_R0small-hyperbolic-I-t1} & 
%			\includegraphics[width=0.40\textwidth]{test2_R0small-parabolic-I-t1} \\
%			\includegraphics[width=0.40\textwidth]{test2_R0small-hyperbolic-I-t2} & 
%			\includegraphics[width=0.40\textwidth]{test2_R0small-parabolic-I-t2} \\
%		\end{tabular} 
%		\caption{Test 1. Evolution of $I$ for reproduction number $R_0<1$ and relaxation times $\tau=1.0$ with  $\lambda^2=1.0$ (left) and $\tau=10^{-4}$ with $\lambda^2=10^4$ (right). Numerical results at output times $t=1$ and $t=2$ (from top to bottom).}
%		\label{fig.test2_R0small_I}
%	\end{center}
%\end{figure}

% R0>1 --------------------------------------------------

%\begin{figure}[!htbp]
%	\begin{center}
%		\begin{tabular}{cc} 
%			\includegraphics[width=0.33\textwidth]{test2_R0big-hyperbolic-S-t2} & 
%			\includegraphics[width=0.33\textwidth]{test2_R0big-parabolic-S-t2} \\
%			\includegraphics[width=0.33\textwidth]{test2_R0big-hyperbolic-S-t4} & 
%			\includegraphics[width=0.33\textwidth]{test2_R0big-parabolic-S-t4} \\
%			\includegraphics[width=0.33\textwidth]{test2_R0big-hyperbolic-S-t8} & 
%			\includegraphics[width=0.33\textwidth]{test2_R0big-parabolic-S-t8} \\
%			\includegraphics[width=0.33\textwidth]{test2_R0big-hyperbolic-S-t10} & 
%			\includegraphics[width=0.33\textwidth]{test2_R0big-parabolic-S-t10} \\
%		\end{tabular} 
%	\end{center}
%		\caption{Test 1. Evolution of $S$ for reproduction number $R_0>1$ and relaxation times $\tau=1.0$ with  $\lambda^2=1.0$ (left, hyperbolic regime) and $\tau=10^{-4}$ with $\lambda^2=10^4$ (right, parabolic regime). Numerical results at output times $t=2$, $t=4$, $t=8$ and $t=10$ (from top to bottom).}
%		\label{fig.test2_R0big_S}
%\end{figure}

%------------------------------------------------------------------------------------
\subsubsection{Test 2. Commuter and non commuter interactions in disease spread}
\begin{figure}[!htbp]
	\begin{center}
		\begin{tabular}{ccc} 
			\includegraphics[width=0.32\textwidth]{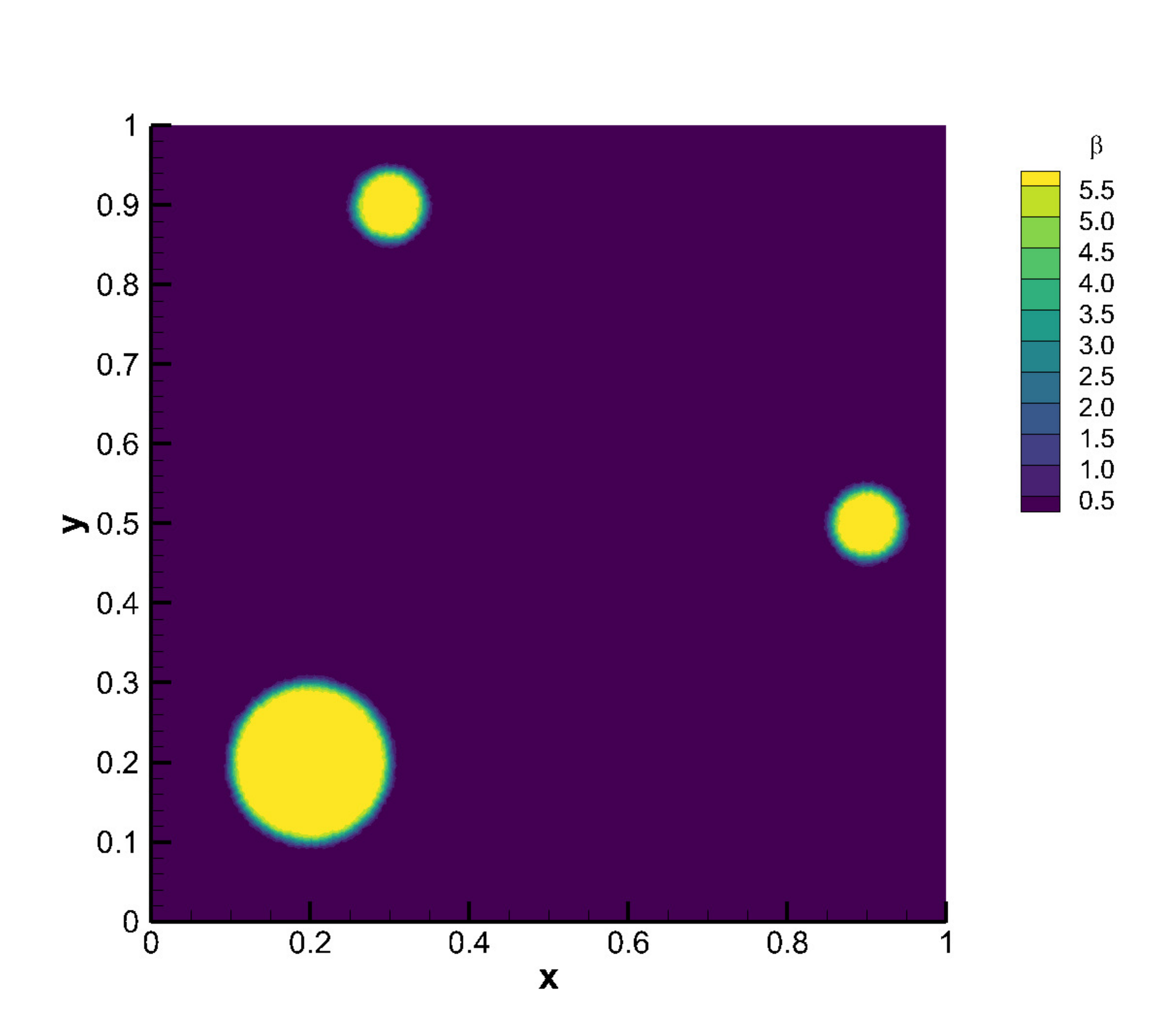} & 
			\includegraphics[width=0.32\textwidth]{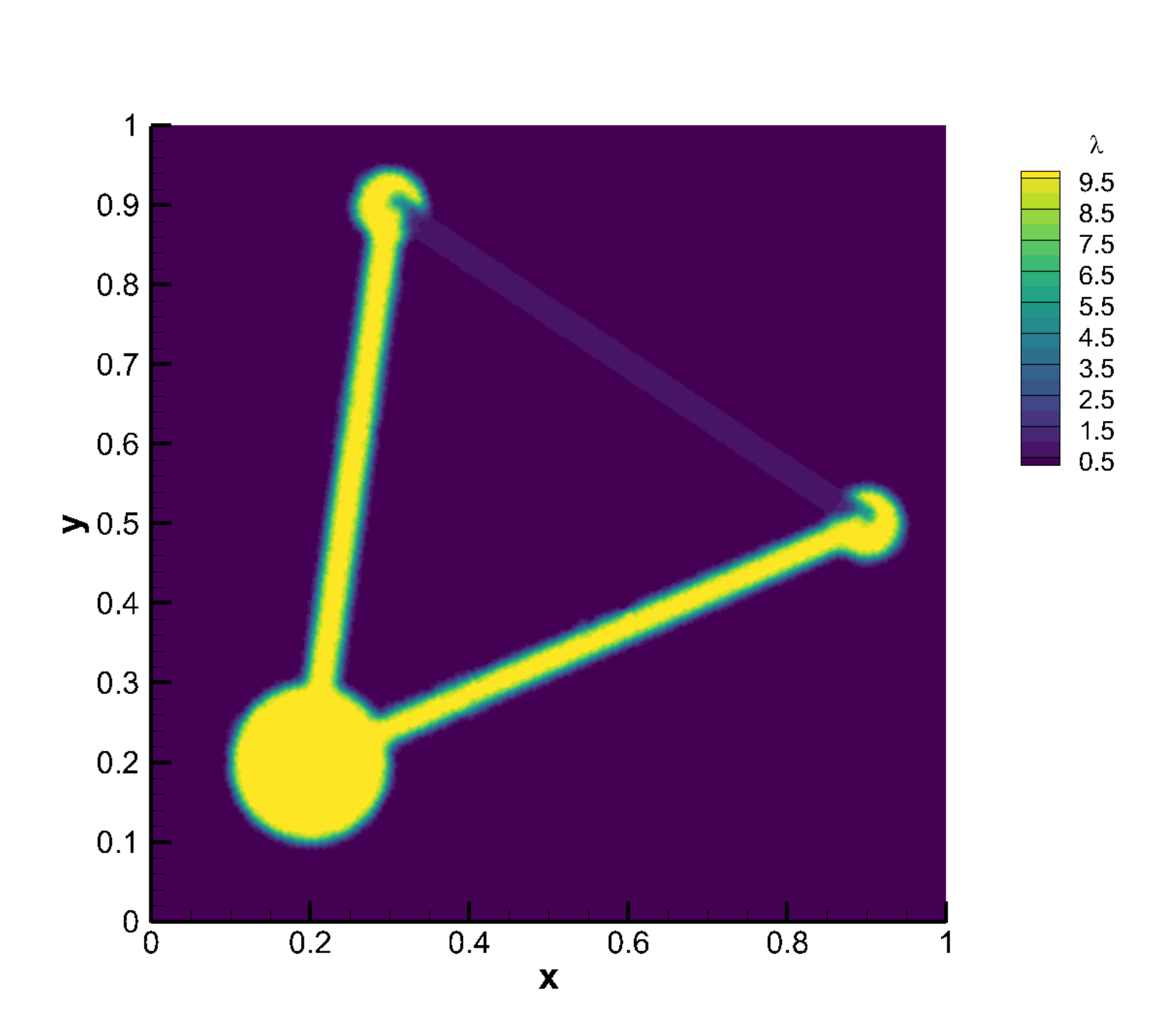} &
			\includegraphics[width=0.32\textwidth]{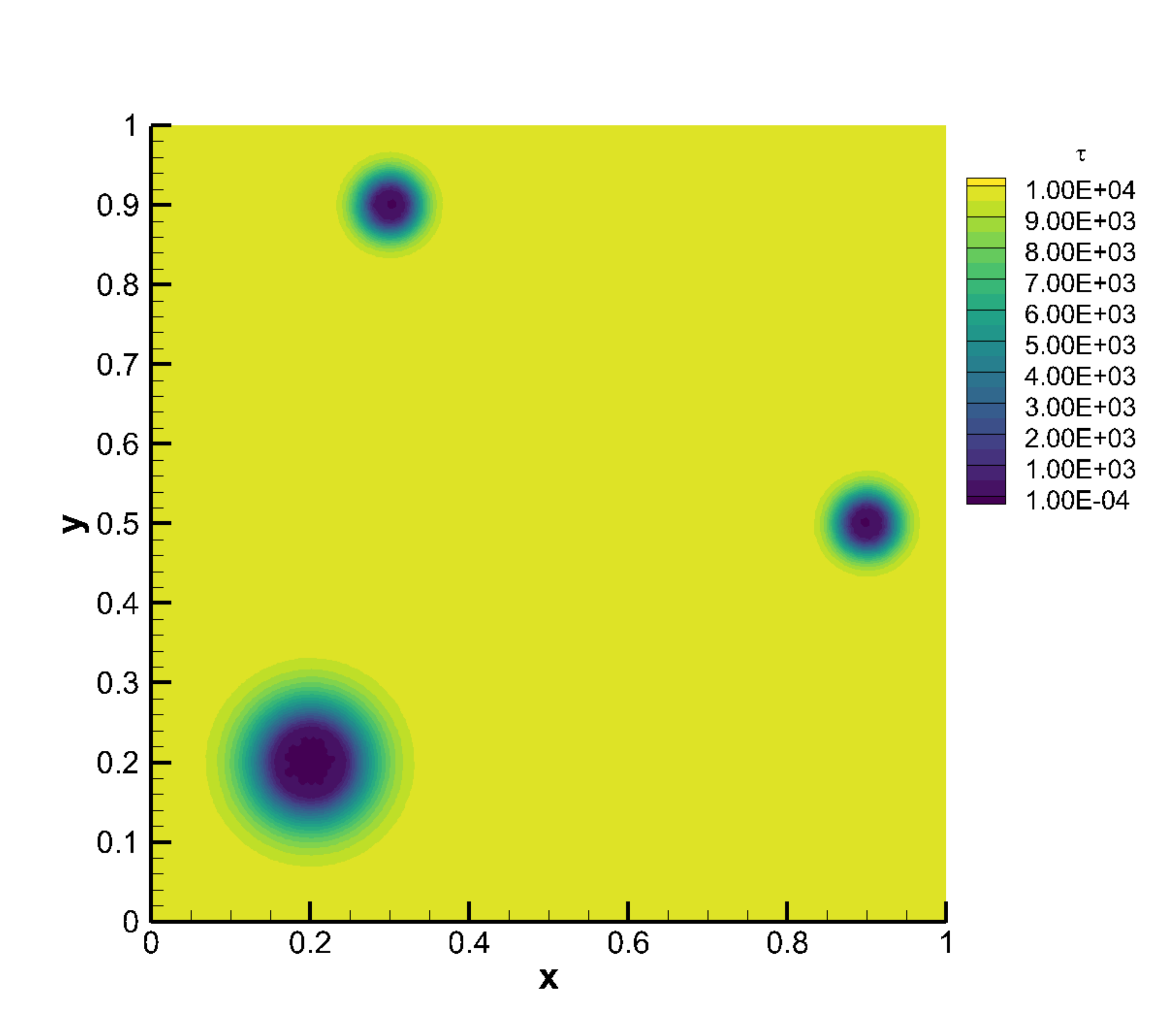} \\   
		\end{tabular} 
	\end{center}
		\caption{Test 2. Initial condition for $\beta$ (left), $\lambda$ (middle) and $\tau$ (right) for the hyperbolic regime. }
		\label{fig.test1_IC}
\end{figure}

\begin{figure}[!htbp]
	\begin{center}
	\begin{tabular}{ccc} 
			\includegraphics[width=0.32\textwidth]{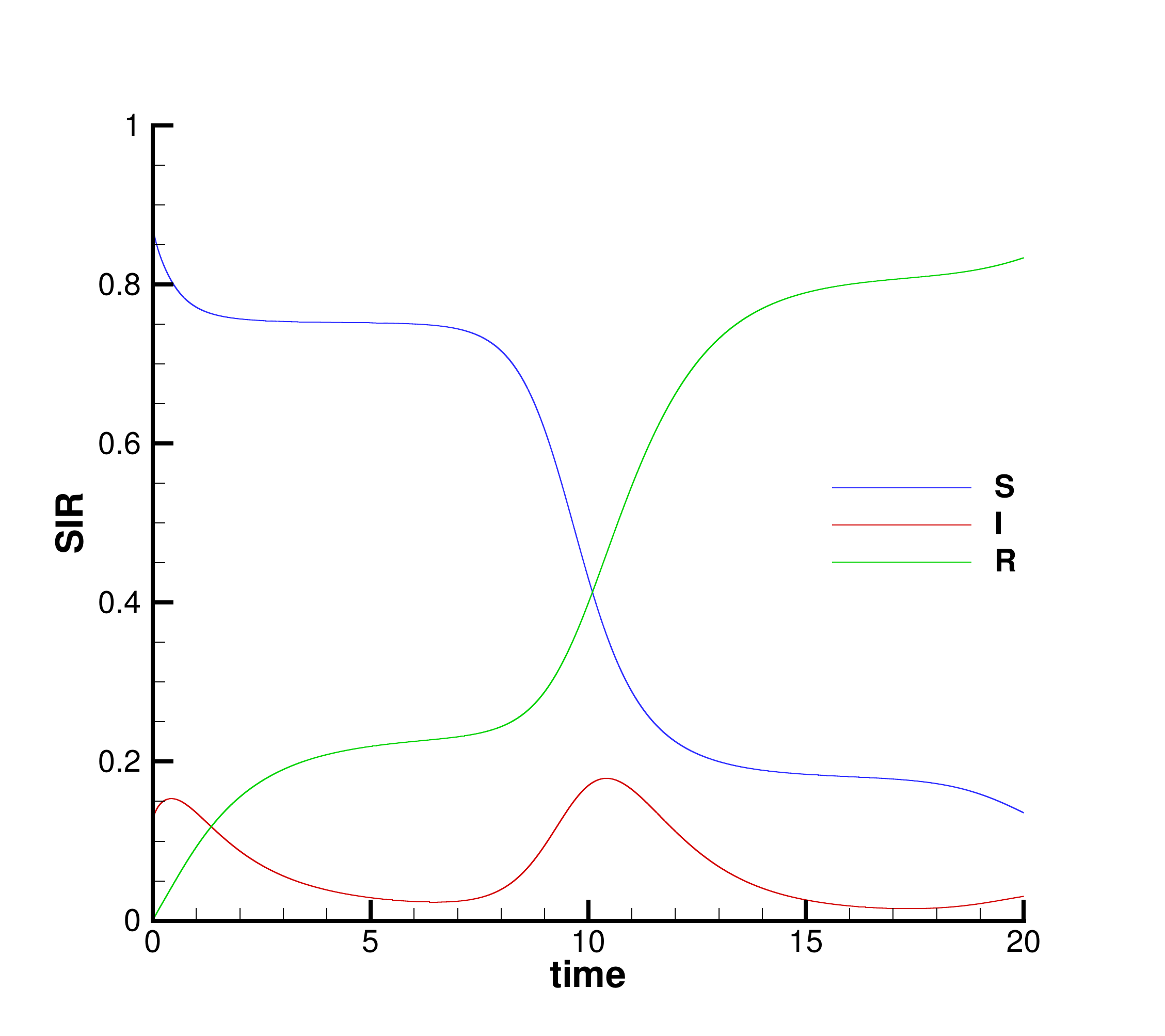} & 
			\includegraphics[width=0.32\textwidth]{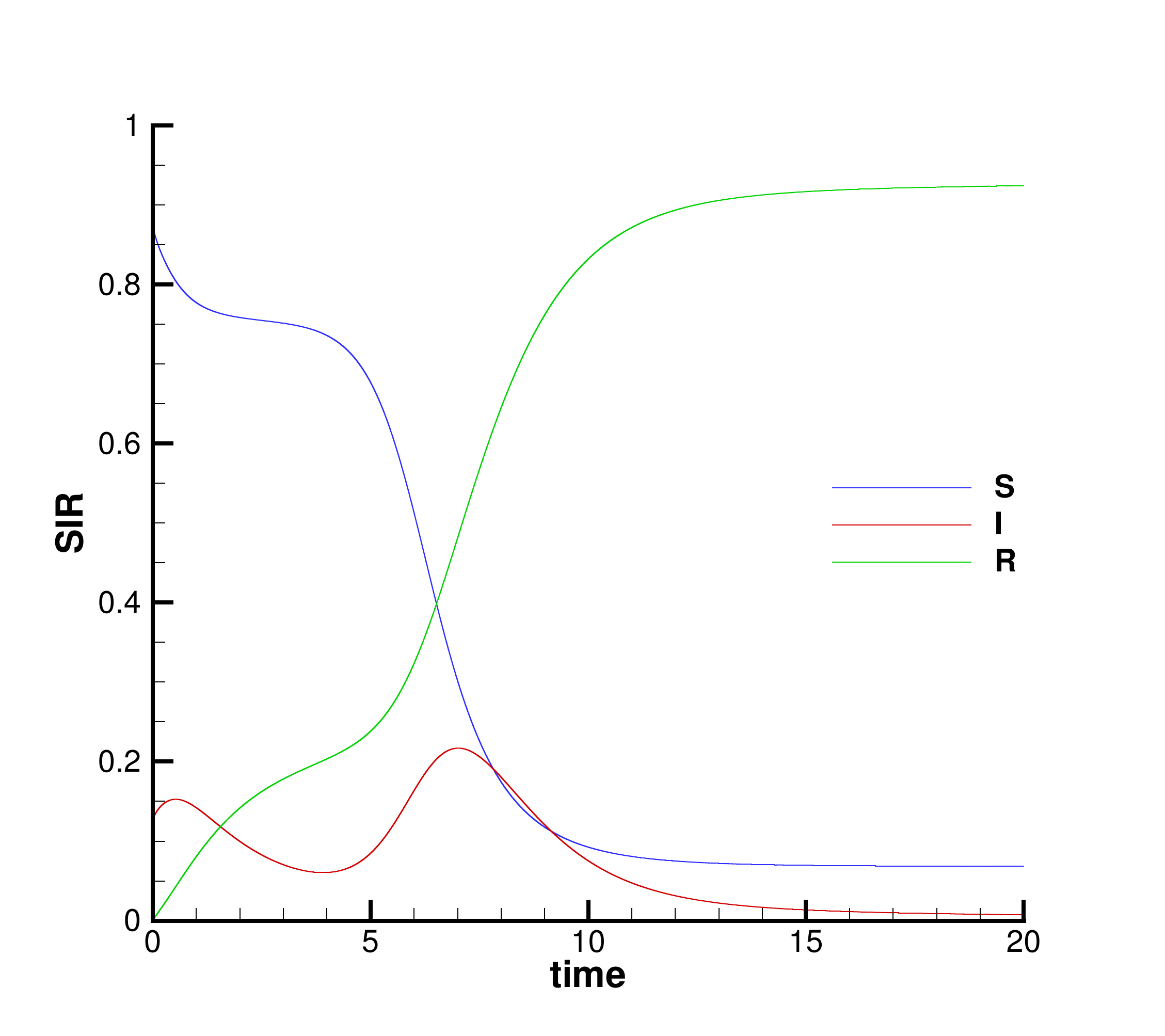} &
			\includegraphics[width=0.32\textwidth]{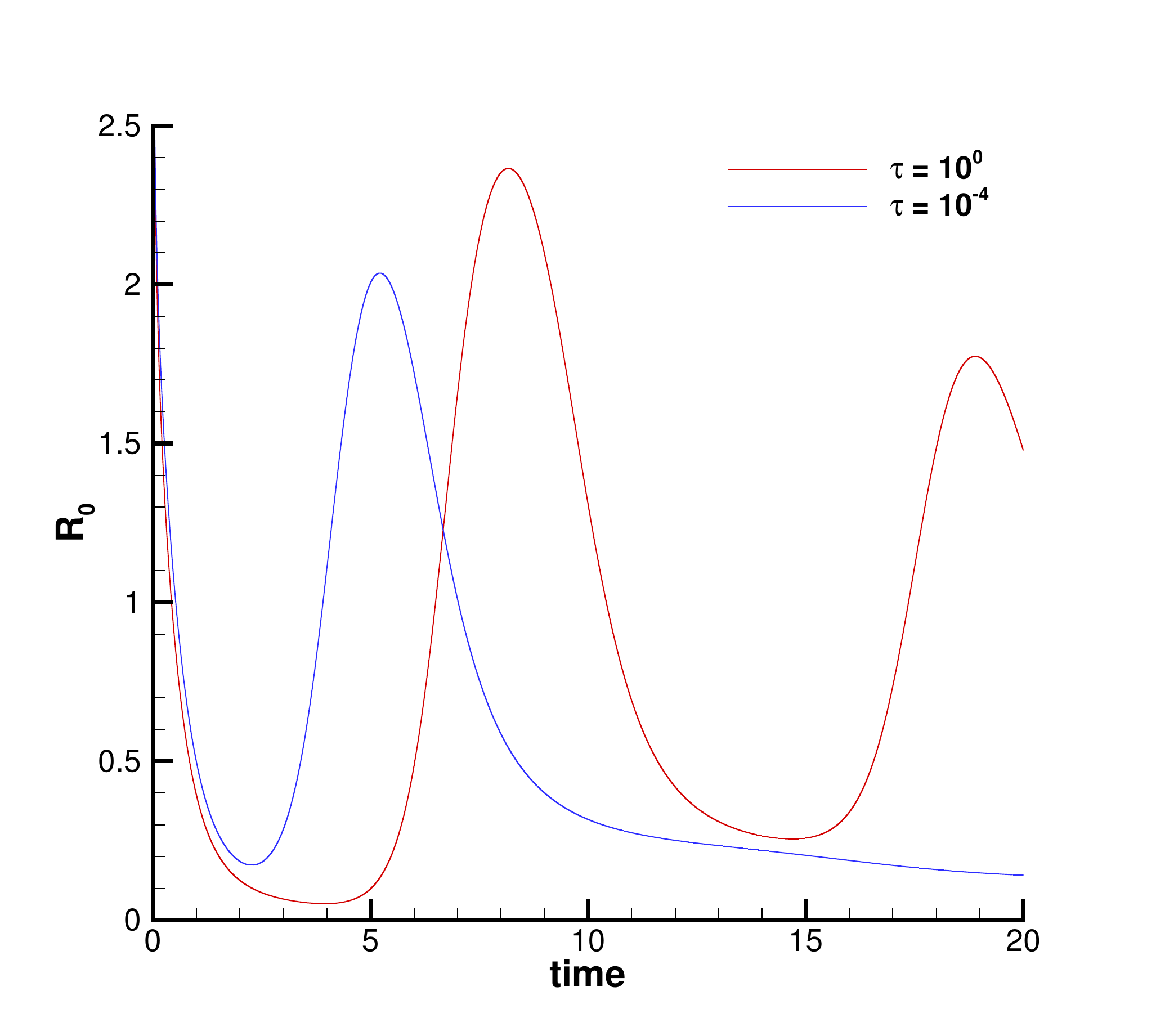} \\  	
			\end{tabular}		
	\end{center}
			\caption{Test 2. Time evolution of the susceptible ($S$), infected ($I$) and recovered ($R$) population in the hyperbolic regime with $\tau=10^{4}$ (left) and parabolic regime with $\tau=10^{-4}$ (middle). Hyperbolic velocities $\lambda^2=10^2$ and SIR parameters $\beta=6$ and $\gamma=1$. Right: evolution of index $R_0$ for hyperbolic (red line) and parabolic (blue line) regime.}
		\label{fig.test1_SIR}
\end{figure}

\begin{figure}[!htbp]
	\begin{center}
		\begin{tabular}{cc} 
			\includegraphics[width=0.33\textwidth]{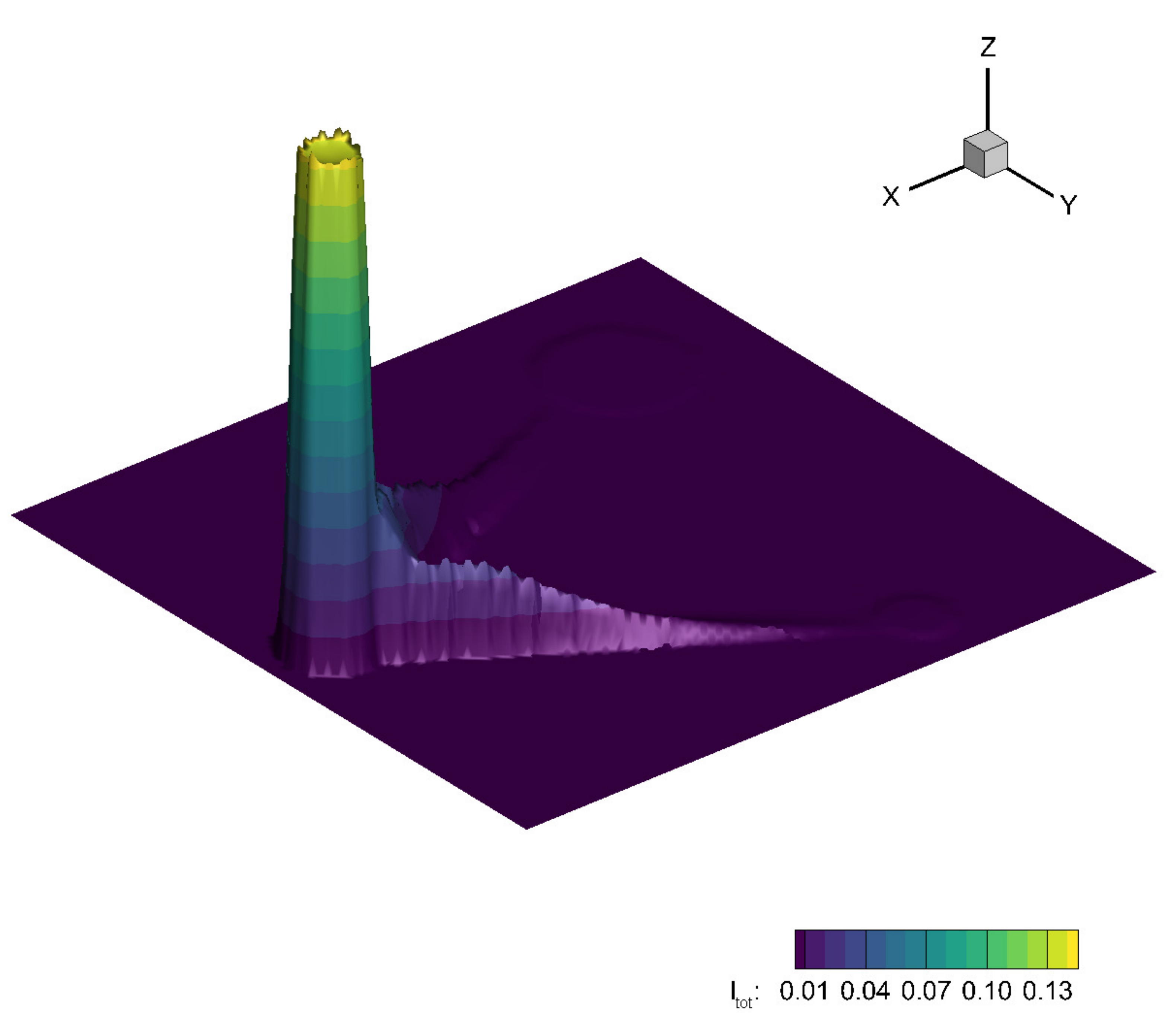} & 
			\includegraphics[width=0.33\textwidth]{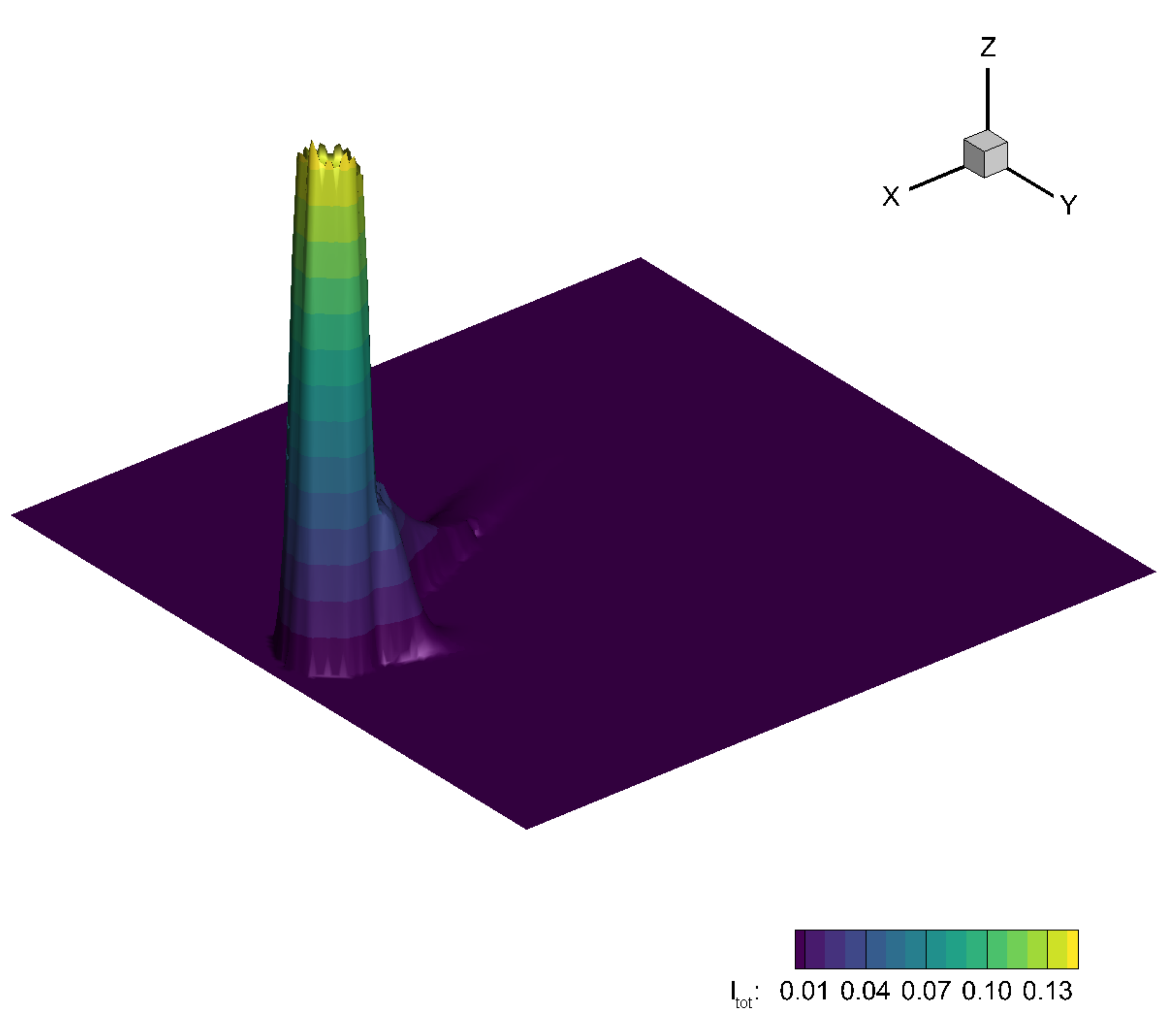} \\
			\includegraphics[width=0.33\textwidth]{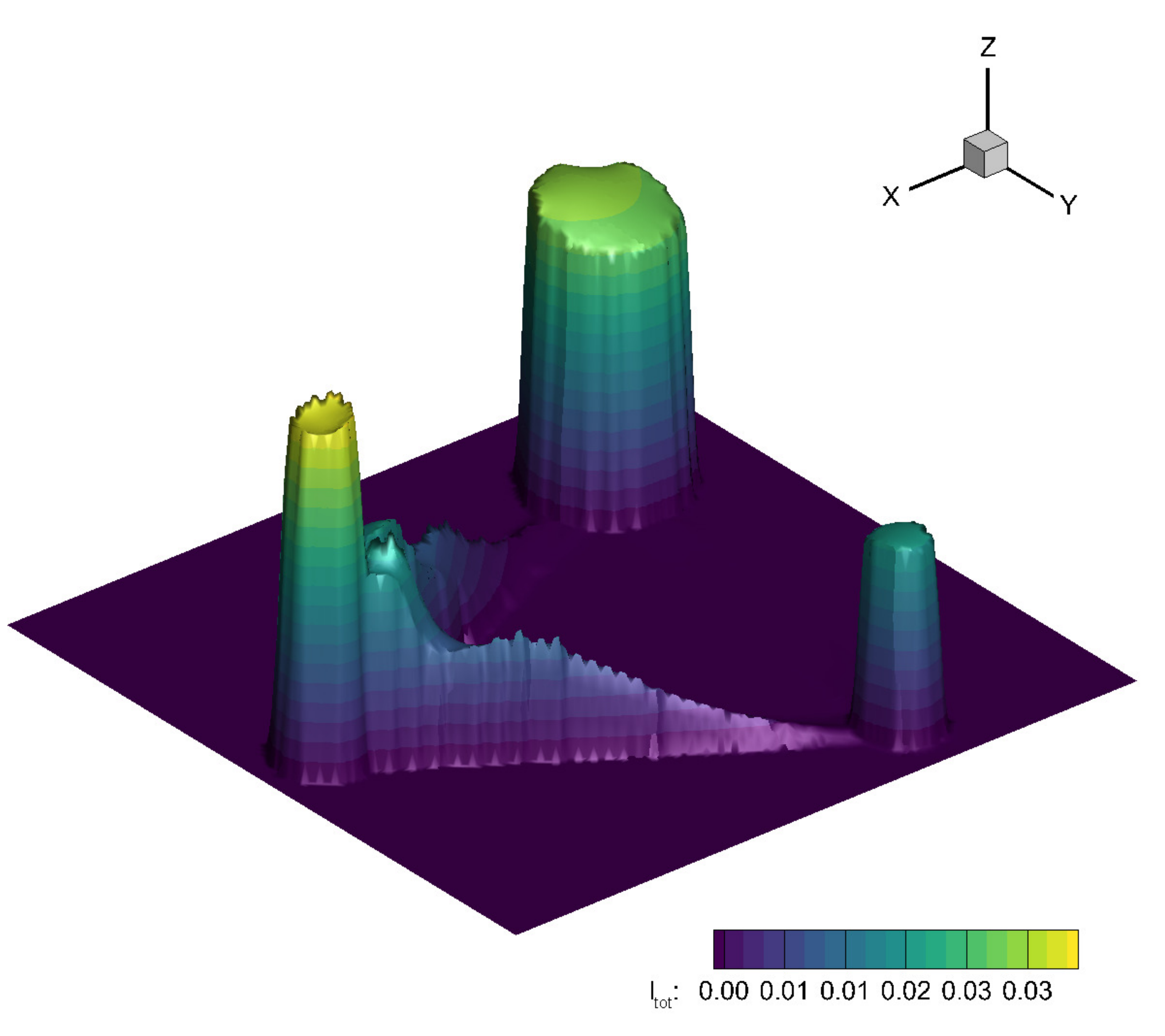} & 
			\includegraphics[width=0.33\textwidth]{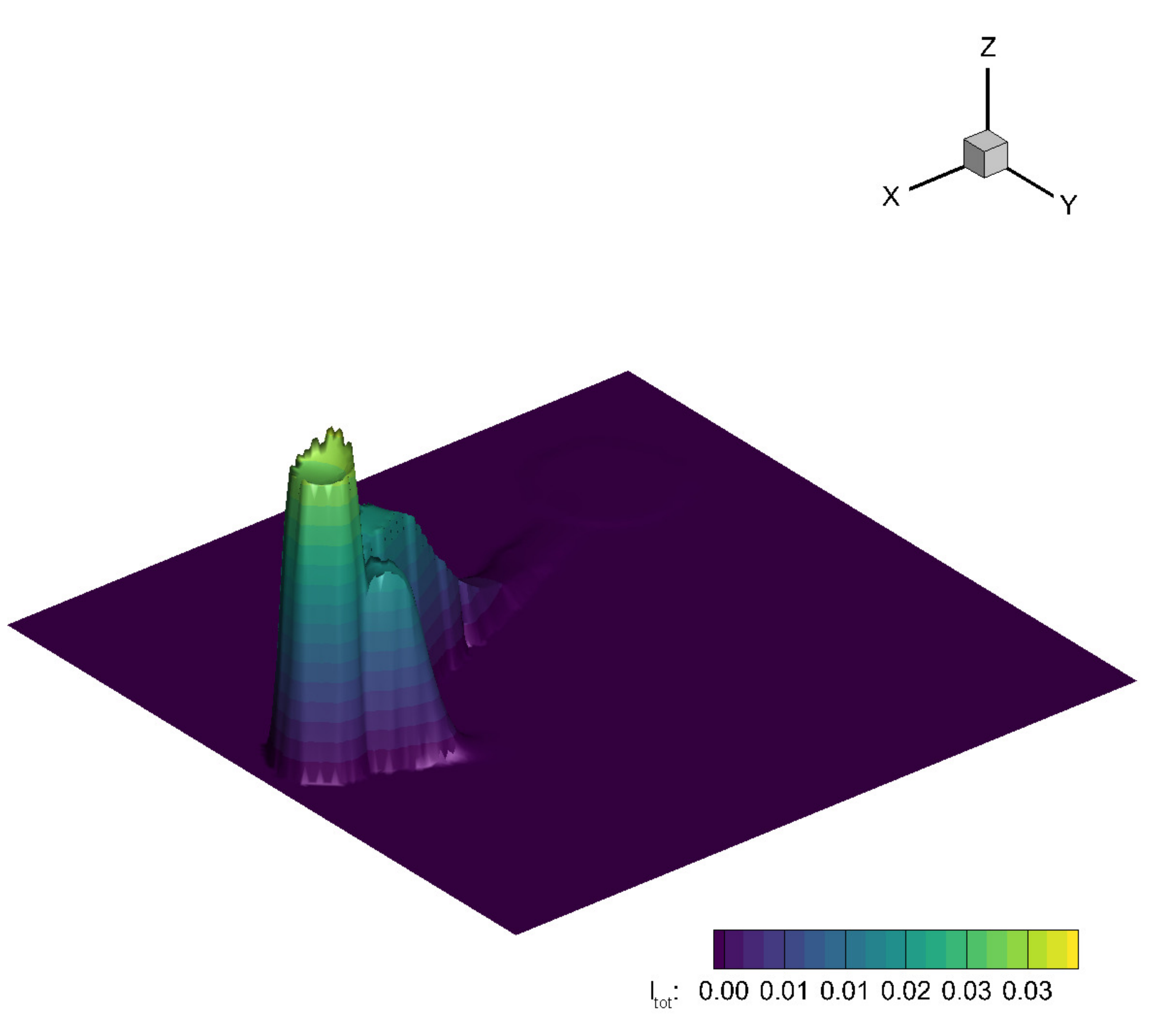} \\
			\includegraphics[width=0.33\textwidth]{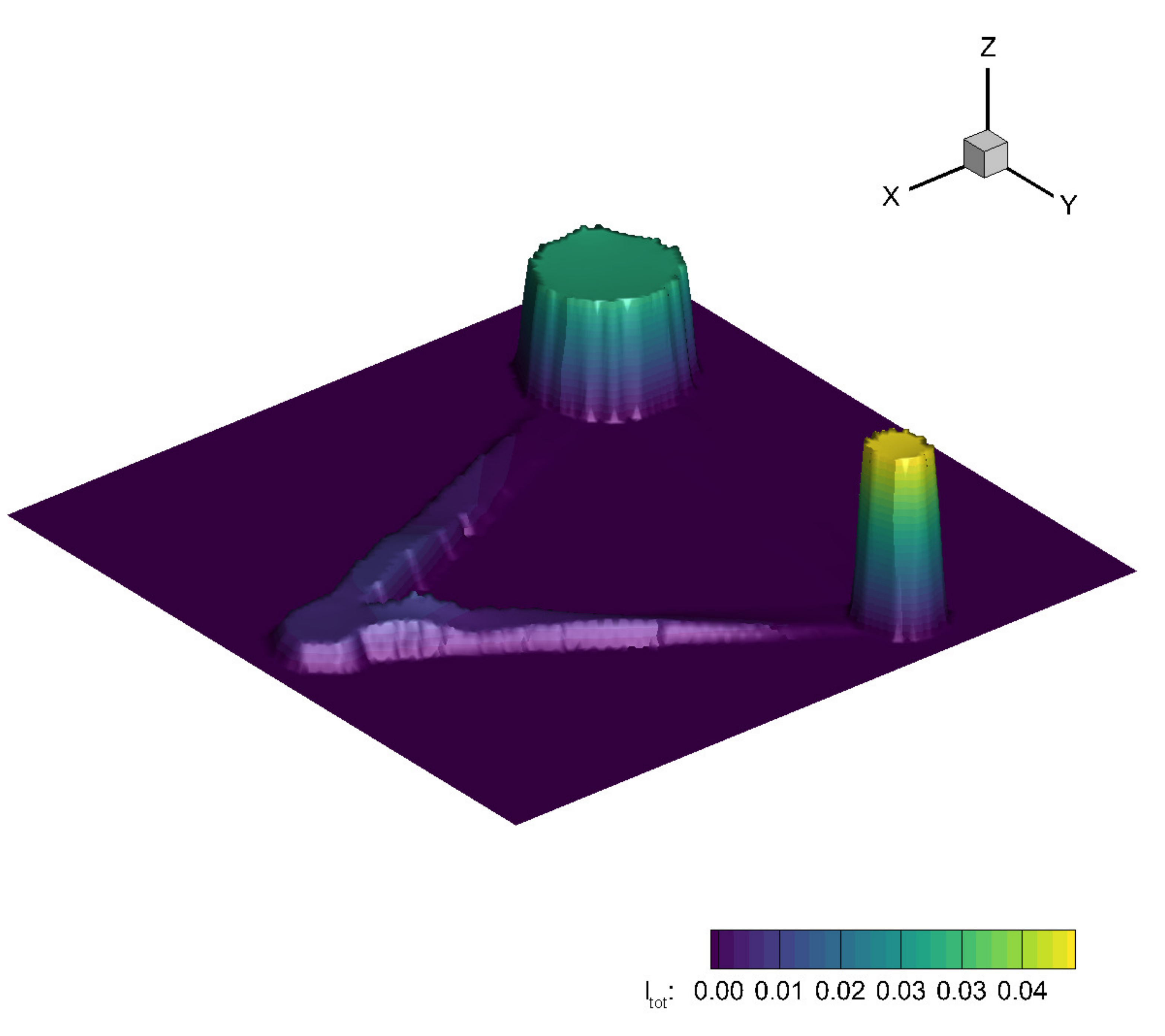} & 
			\includegraphics[width=0.33\textwidth]{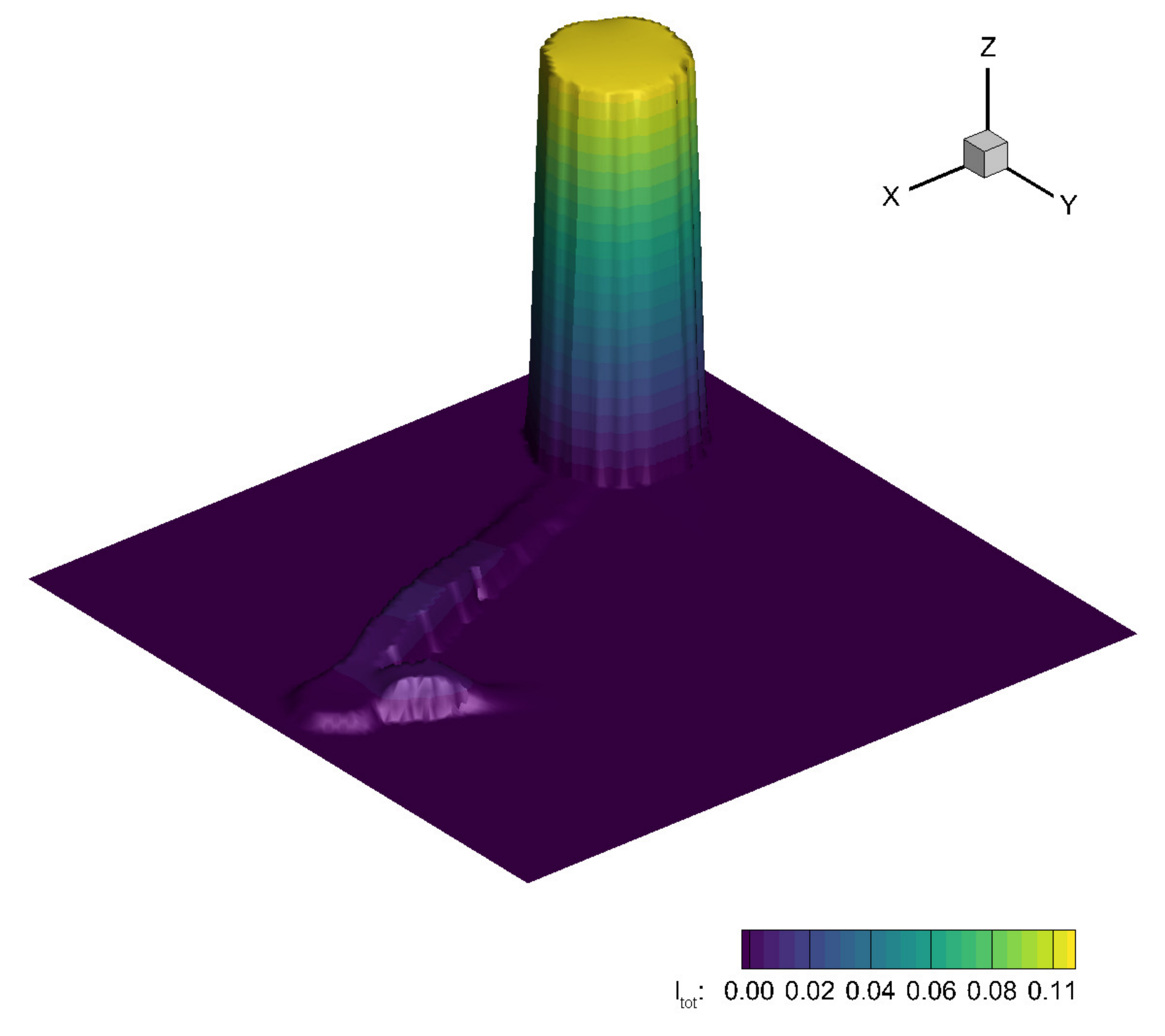} \\
			\includegraphics[width=0.33\textwidth]{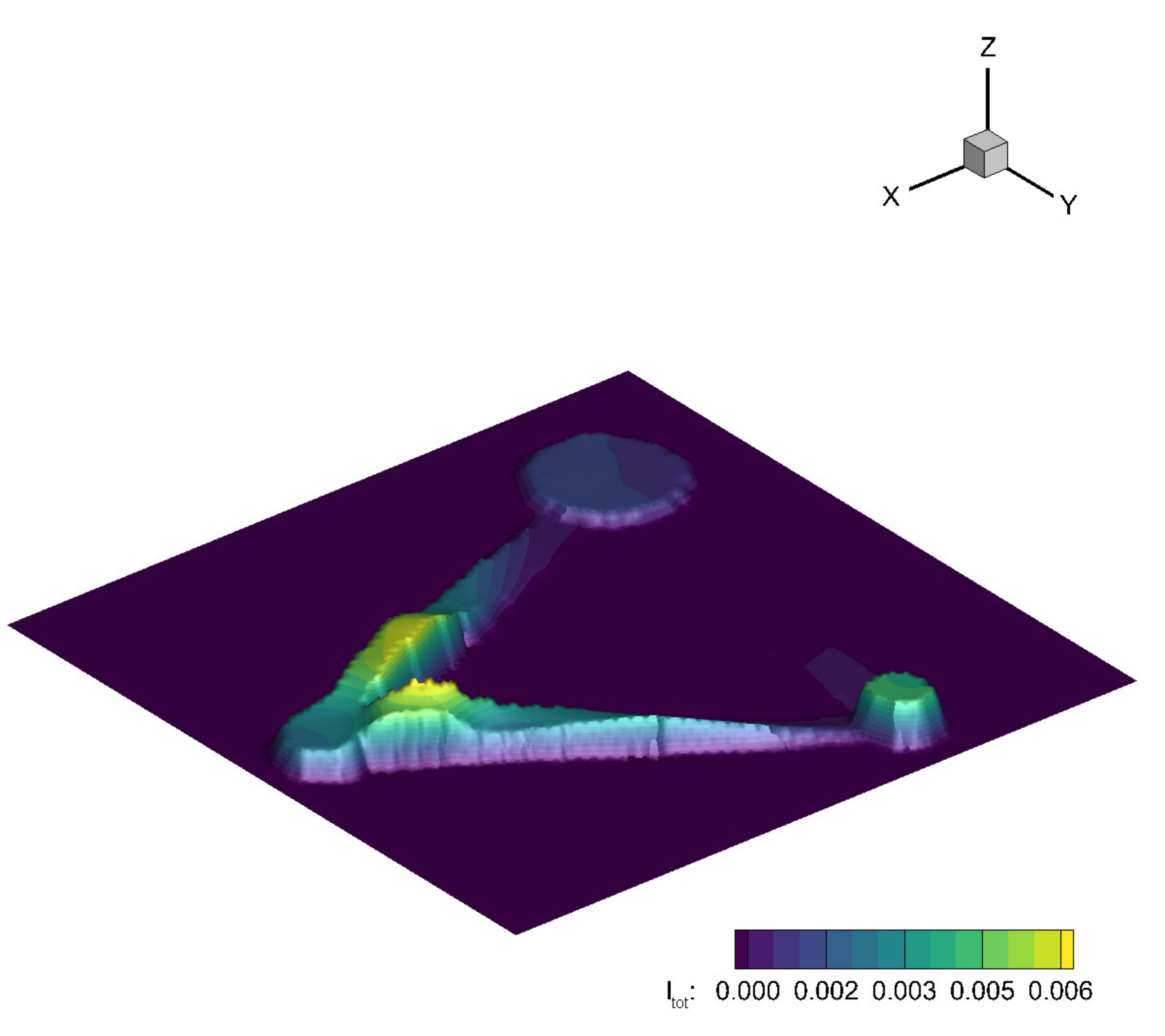} & 
			\includegraphics[width=0.33\textwidth]{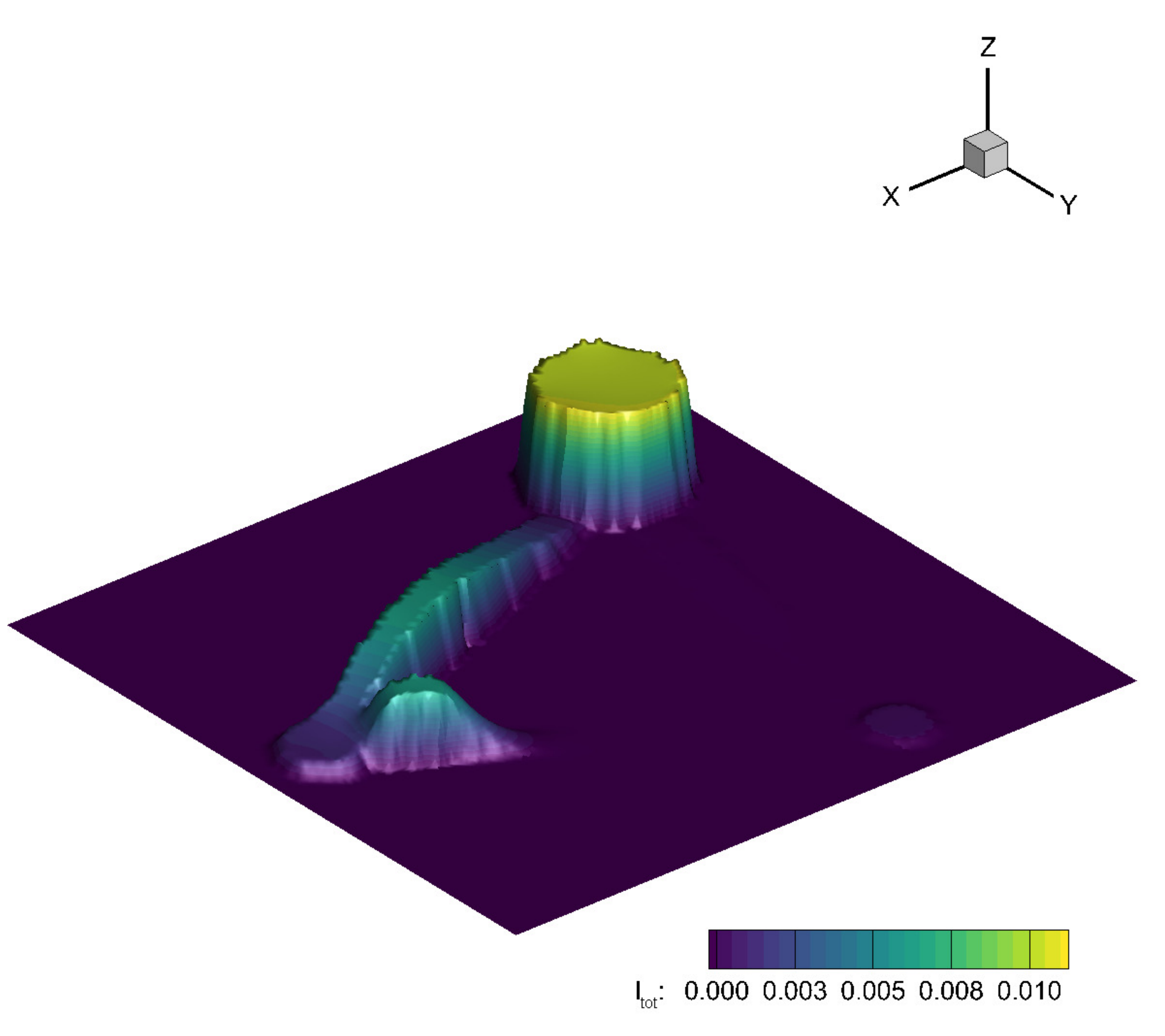} \\
		\end{tabular} 
	\end{center}
		\caption{Test 2. Distribution of total infected population $I_T=I+\IO$ in hyperbolic regime with $\tau=10^{4}$ (left) and parabolic regime with $\tau=10^{-4}$ (right). Hyperbolic velocities $\lambda^2=10^2$ and SIR parameters $\beta=6$ and $\gamma=1$. Output at times $t=2.5$, $t=5$, $t=10$ and $t=15$ (from top to bottom).}
\label{fig.test1_I3D}
\end{figure}

The MK-SIR model \eqref{eq:kineticc}-\eqref{eq:diffuse} was employed to simulate the spread of a disease in an environment that aims to reproduce three urban areas of different sizes connected by paths to which we assigned different propagation speeds. The computational domain is given by $\Omega=[0;1]^2$ with zero-flux boundaries set everywhere and it is discretized with a total number of polygonal cells $N_E=10112$. The urban areas $A$, $B$ and $C$ are expanded around the centers of coordinates
\begin{equation}
\x_A = (x_A,y_A) = (0.2,0.2), \quad \x_B= (x_B,y_B) =(0.9,0.5), \quad \x_C= (x_C,y_C) =(0.3,0.9).
\end{equation} 
The propagation speed is $\lambda^2=10^2$ along the straight paths $A$-$B$ and $A$-$C$, while it drops to $\lambda^2=10^{-2}$ along the connection $B$-$C$, thus we expect a lower level of individuals moving between the urban areas $B$ and $C$. Each path has a width of $h=0.04$ and in the remaining part of the computational domain we set $\lambda^2=10^{-12}$, this means that almost all individuals move or in the urban region in a diffusive regime or along the connecting paths. The recovery rate is set to $\gamma=1$ and the contact rate is space-dependent, that is
\begin{equation}
\beta = \left\{ \begin{array}{lll}
6 & \textnormal{if} & (S_T+I_T)>0 \\
0 & \textnormal{if} & (S_T+I_T)=0,
\end{array} \right.
\end{equation}
with the initial distribution of the populations given by
\begin{eqnarray}
S_T &=& \max \left[ 0, \, -100 (x-x_A)^2-100(y-y_A)^2\right] + \max \left[ 0, \, 1-500(x-x_B)^2-500(y-y_B)^2\right] \nonumber \\
&+& \max \left[ 0, \, 1-500(x-x_C)^2-500(y-y_C)^2 \right], \nonumber \\
I_T &=& \max \left[ 0, \, 1-500(x-x_B)^2-500(y-y_B)^2 \right], \nonumber \\
R_T &=& 0.
\end{eqnarray}
We furthermore impose that the commuters are 
\begin{equation}
S = 0.01 \cdot S_T, \qquad I = 0.8 \cdot I_T, \qquad R = R_T,
\end{equation}
thus, for the sake of conservation, leading to $\SO=S_T-S$ and $\IO=I_T-I$ for the non commuters. Finally, the relaxation time is also spatially heterogeneous. Let us define an auxiliary variable $\tilde{\tau}$ which writes
\begin{equation}
\tilde{\tau} = \tau_r \, \left[ e^{-\frac{1}{2} \frac{\left( (x-x_A)^2+(y-y_A)^2\right) }{s_1^2} } + e^{-\frac{1}{2} \frac{\left( (x-x_B)^2+(y-y_B)^2\right) }{s_2^2} } 
+ \, e^{-\frac{1}{2} \frac{\left( (x-x_C)^2+(y-y_C)^2\right) }{s_2^2} } \right],
\end{equation}
with $s_1=0.05$ and $s_2=0.025$. The parameter $\tau_r=10^{4}$ is set to define a hyperbolic regime, while $\tau_r=10^{-4}$ allows the diffusive system to be recovered. The relaxation time within each urban area is always defined as $\tau_0=10^{-4}$, therefore the initial distribution of the space-dependent relaxation time is given by
\begin{equation}
\tau = \max \left[ \tau_0, \, \tau_r - \frac{3}{2} \tilde{\tau} \right].
\end{equation}  
Figure \ref{fig.test1_IC} depicts the initial distribution for $\beta$, $\lambda$ and $\tau$ in the hyperbolic regime configuration. 
The final time of the simulation is $t_f=20$ and the evolution of the total susceptible, infected and recovered is shown in Figure \ref{fig.test1_SIR} for both the hyperbolic and the diffusive regimes for the populations of commuters. It is interesting to observe the appearance of a second wave of infected due to the spatial dynamics in the region. Additionally, the dynamic induced by the diffusive regime is faster compared to the propagation of the disease in the hyperbolic configuration. This is also evident from the time evolution of the current reproduction number $R_0$ reported in Figure \ref{fig.test1_SIR}. A comparison between hyperbolic and parabolic numerical results is shown at different output times for %both  $S_T=S+\SO$ and  
$I_T=I+\IO$ in Figure %\ref{fig.test1_S3D} and 
\ref{fig.test1_I3D}%, respectively
. Let us observe how the relaxation time drastically determines the dynamic, while the same set of governing equations is maintained. In fact it is clear from Figure \ref{fig.test1_SIR} that when a parabolic regime is chosen the epidemic spreads at much faster pace.
%\begin{figure}[!htbp]
%	\begin{center}
%		\begin{tabular}{cc} 
%			\includegraphics[width=0.33\textwidth]{test1_S3D_t25-hyp} & 
%			\includegraphics[width=0.33\textwidth]{test1_S3D_t25-par} \\
%			\includegraphics[width=0.33\textwidth]{test1_S3D_t50-hyp} & 
%			\includegraphics[width=0.33\textwidth]{test1_S3D_t50-par} \\
%			\includegraphics[width=0.33\textwidth]{test1_S3D_t100-hyp} & 
%			\includegraphics[width=0.33\textwidth]{test1_S3D_t100-par} \\
%			\includegraphics[width=0.33\textwidth]{test1_S3D_t150-hyp} & 
%			\includegraphics[width=0.33\textwidth]{test1_S3D_t150-par} \\
%		\end{tabular} 
%	\end{center}
%		\caption{Test 2. Hyperbolic velocities $\lambda^2=10^2$ and SIR parameters $\beta=6$ and $\gamma=1$. Distribution of total susceptible population $S_T=S+\SO$ in hyperbolic regime with $\tau=10^{4}$ (left) and parabolic regime with $\tau=10^{-4}$ (right). Output at times $t=2.5$, $t=5$, $t=10$ and $t=15$ (from top to bottom).}
%		\label{fig.test1_S3D}
%\end{figure}

\begin{figure}[!htbp]
	\begin{center}
		\begin{tabular}{ccc} 
			\includegraphics[width=0.32\textwidth]{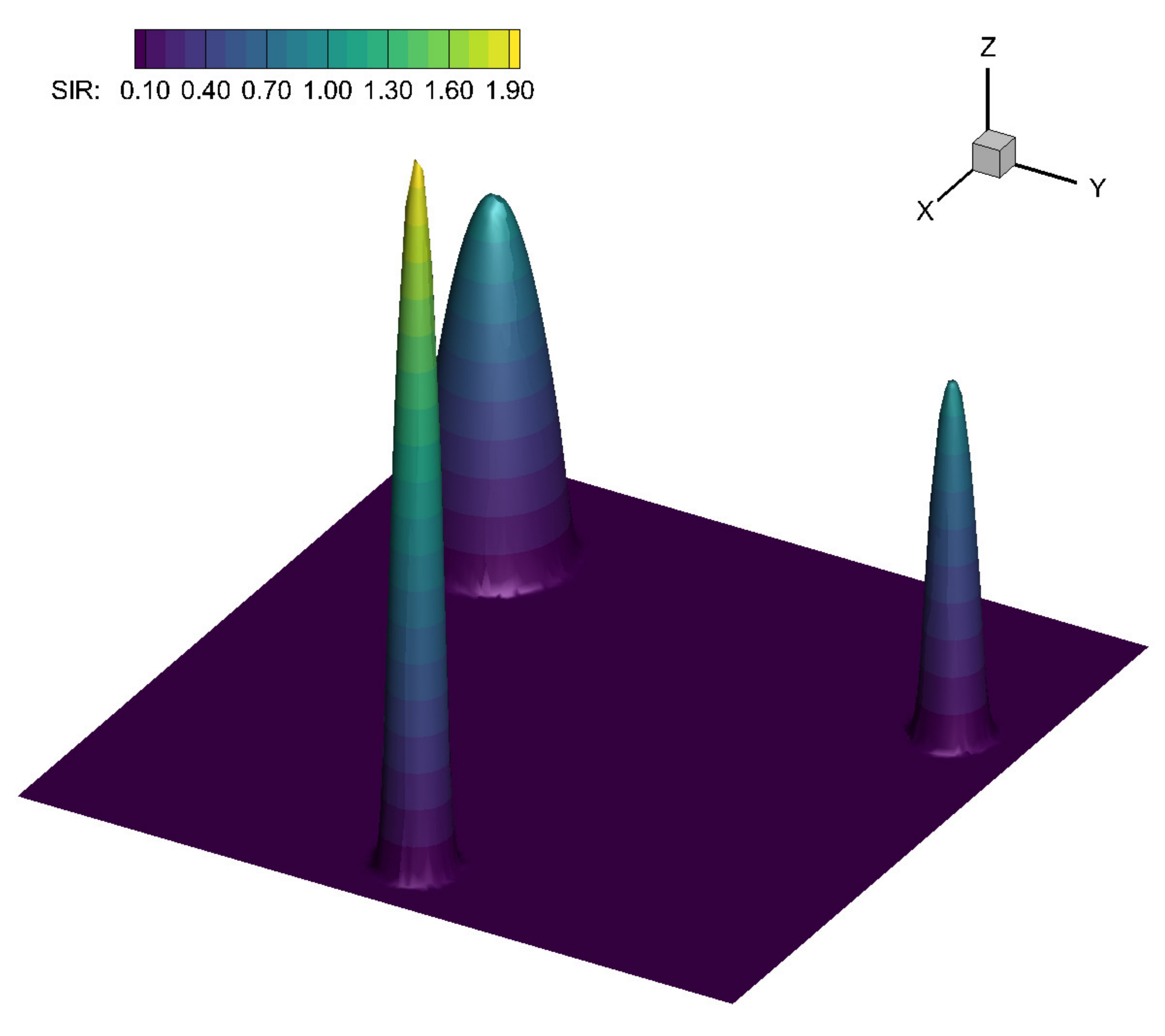} & 
			\includegraphics[width=0.32\textwidth]{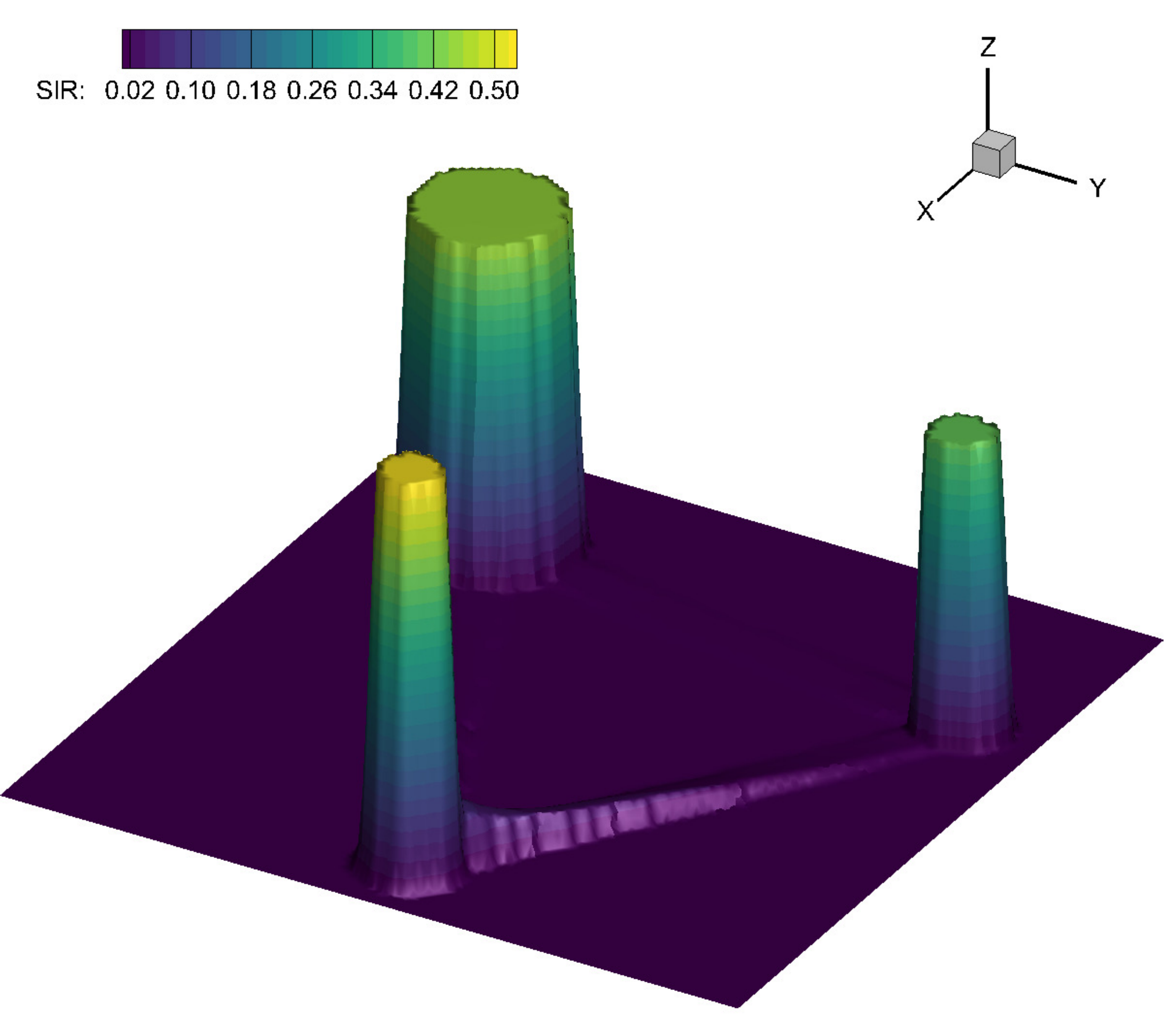} &
			\includegraphics[width=0.32\textwidth]{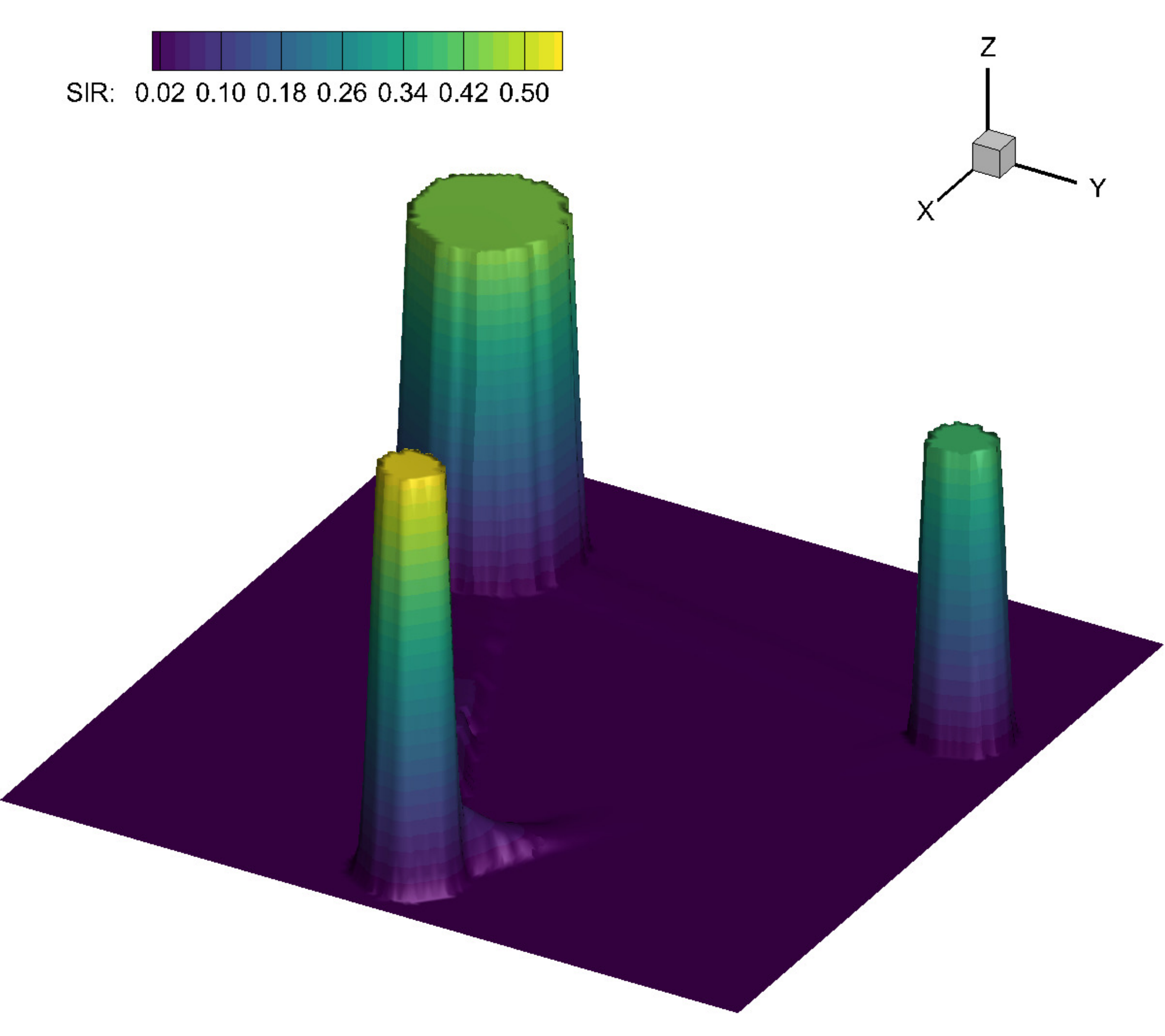} \\
		\end{tabular} 
	\end{center}
		\caption{Test 2. Distribution of total population $(S+\SO)+ (I+\IO) + (R+\RO)$ at initial time (left), at time $t=15$ in hyperbolic regime (middle) and parabolic regime (right). Hyperbolic velocities $\lambda^2=10^2$ and SIR parameters $\beta=6$ and $\gamma=1$. }
			\label{fig.test1_SIR3D}
\end{figure}

Finally, the distribution of the total population is depicted in Figure \ref{fig.test1_SIR3D}, demonstrating that, as expected, the majority the population still remains at their initial locations. Indeed, only the fraction of commuters moves to other locations. 
%a at final time, which is coherent with real world population dynamics. Indeed, only a few people over the entire population of an urban area move to other locations and do not come back. This feature of the model is achieved by the introduction of the background population $\SO$ that is only affected by diffusive phenomena restricted to a region around the city and not by hyperbolic propagation over the entire space.

%------------------------------------------------------------------------------------
\subsection{Application to the spatial spread of COVID-19 in a realistic geographical scenario}
As last numerical example, we apply the novel MK-SEIR model \eqref{eq:kineticc2}-\eqref{eq:diffuse2} to a realistic geographical and epidemic setting. Specifically, we consider the COVID-19 outbreak which took place in a region of northern Italy, Emilia-Romagna, in the first ten days of March 2020. The setup of this simulation requires the knowledge of (i) the computational domain, (ii) the location of the main cities and the boundary definition of the provinces within the region and (iii) the spatial distribution of the population as well as of the infected people at the initial time of the computation. Here, the MK-SEIR model is adopted, thus the initial number of the exposed population (which includes also the asymptomatic) must be estimated and therefore it is affected by uncertainty. We leave, however, the analysis of the influence of such uncertain data to further study and in this example limit ourselves to a deterministic setting in agreement with observations. 

\begin{figure}[!htbp]
	\begin{center}
		\begin{tabular}{cc} 
			\includegraphics[width=0.45\textwidth]{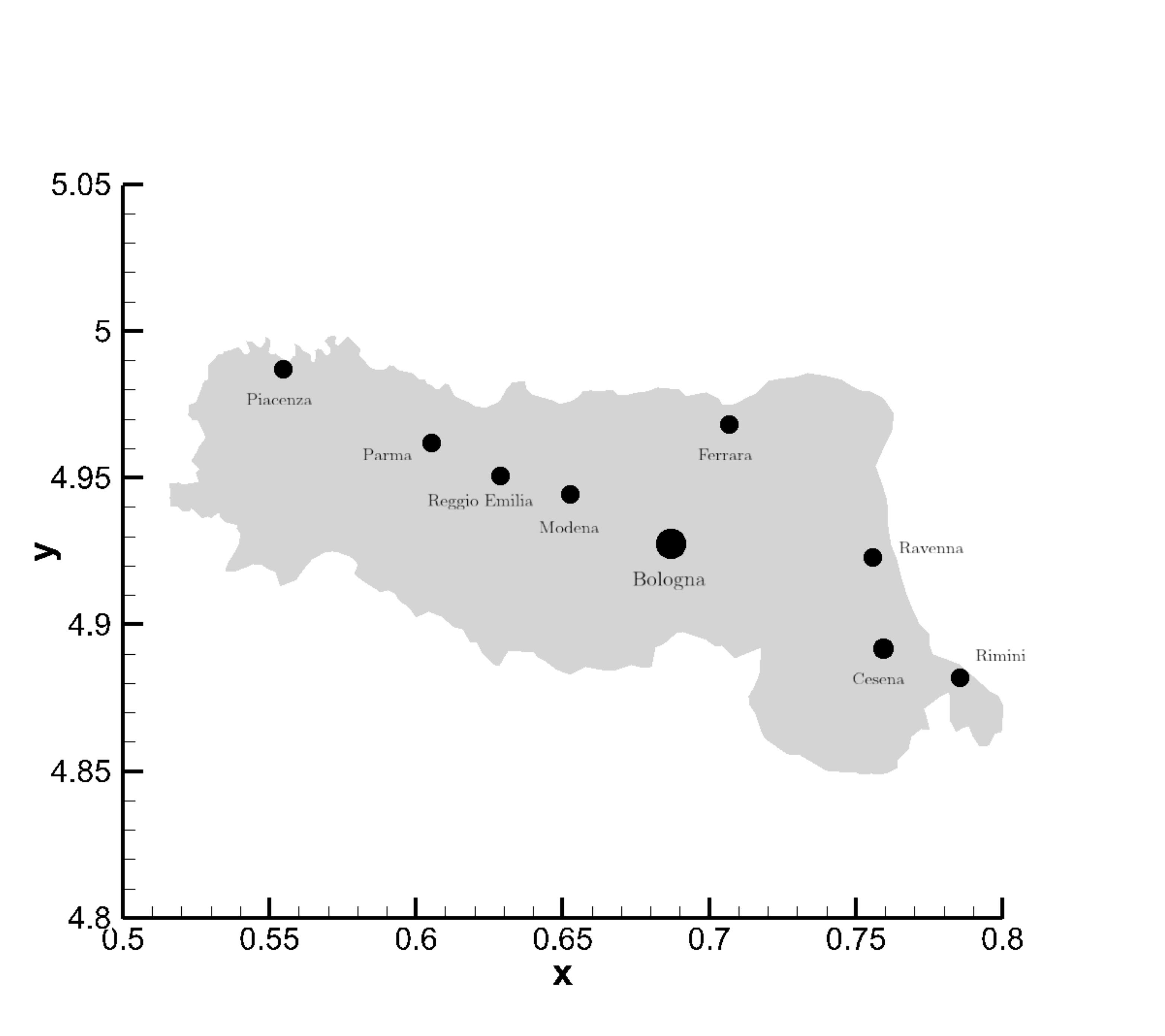} & 
			\includegraphics[width=0.45\textwidth]{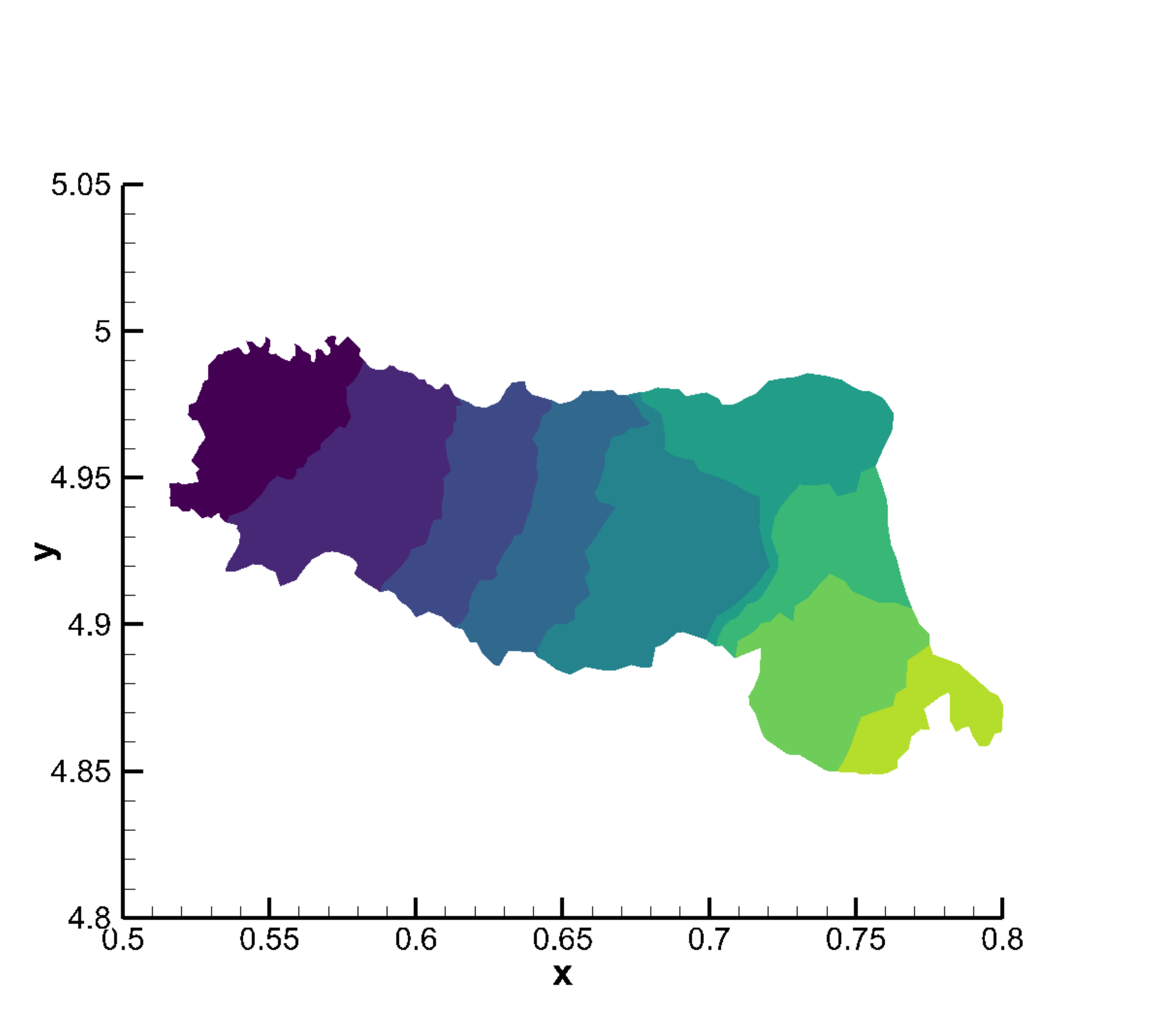} \\
			\includegraphics[width=0.45\textwidth]{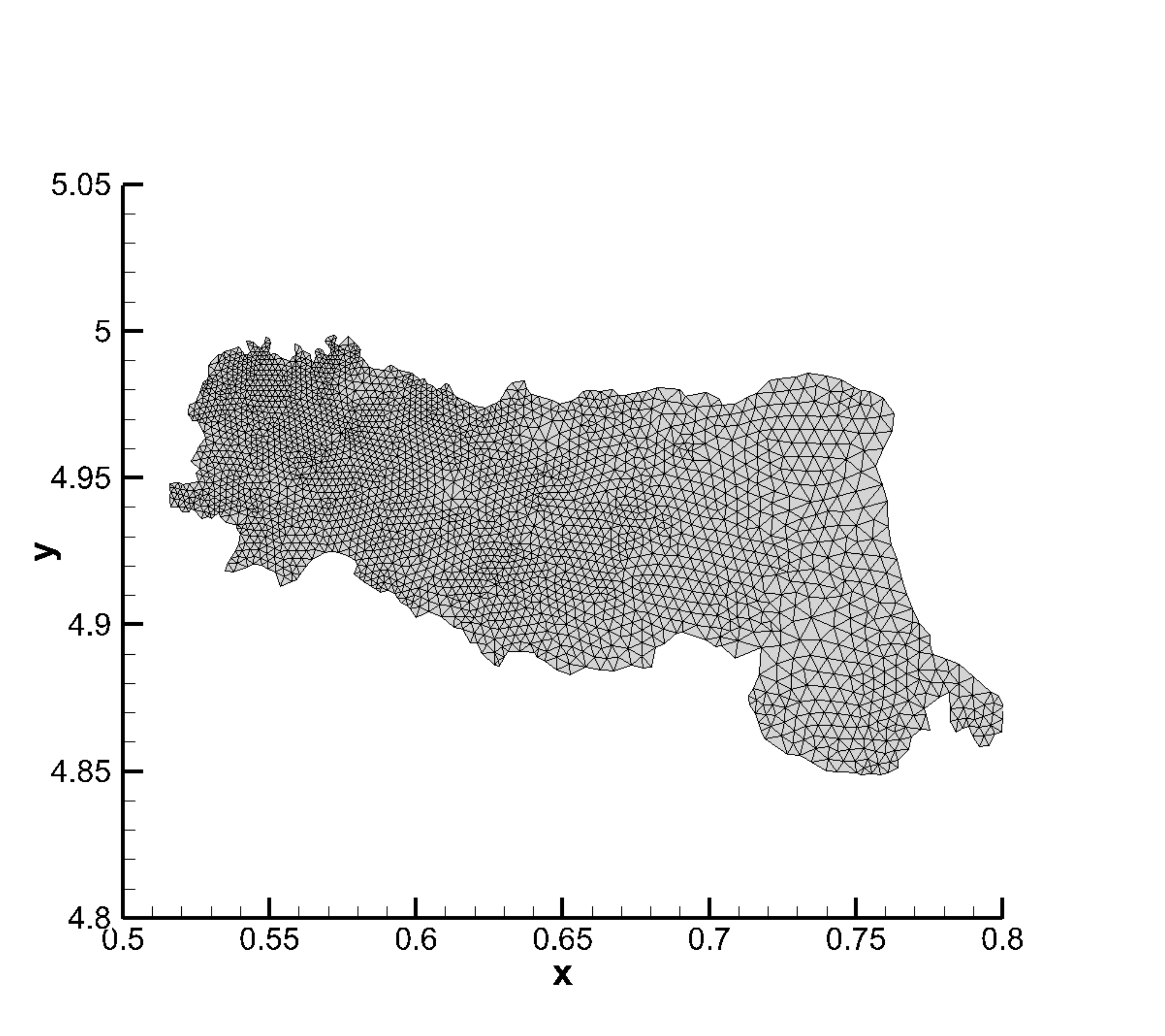} & 
			\includegraphics[width=0.45\textwidth]{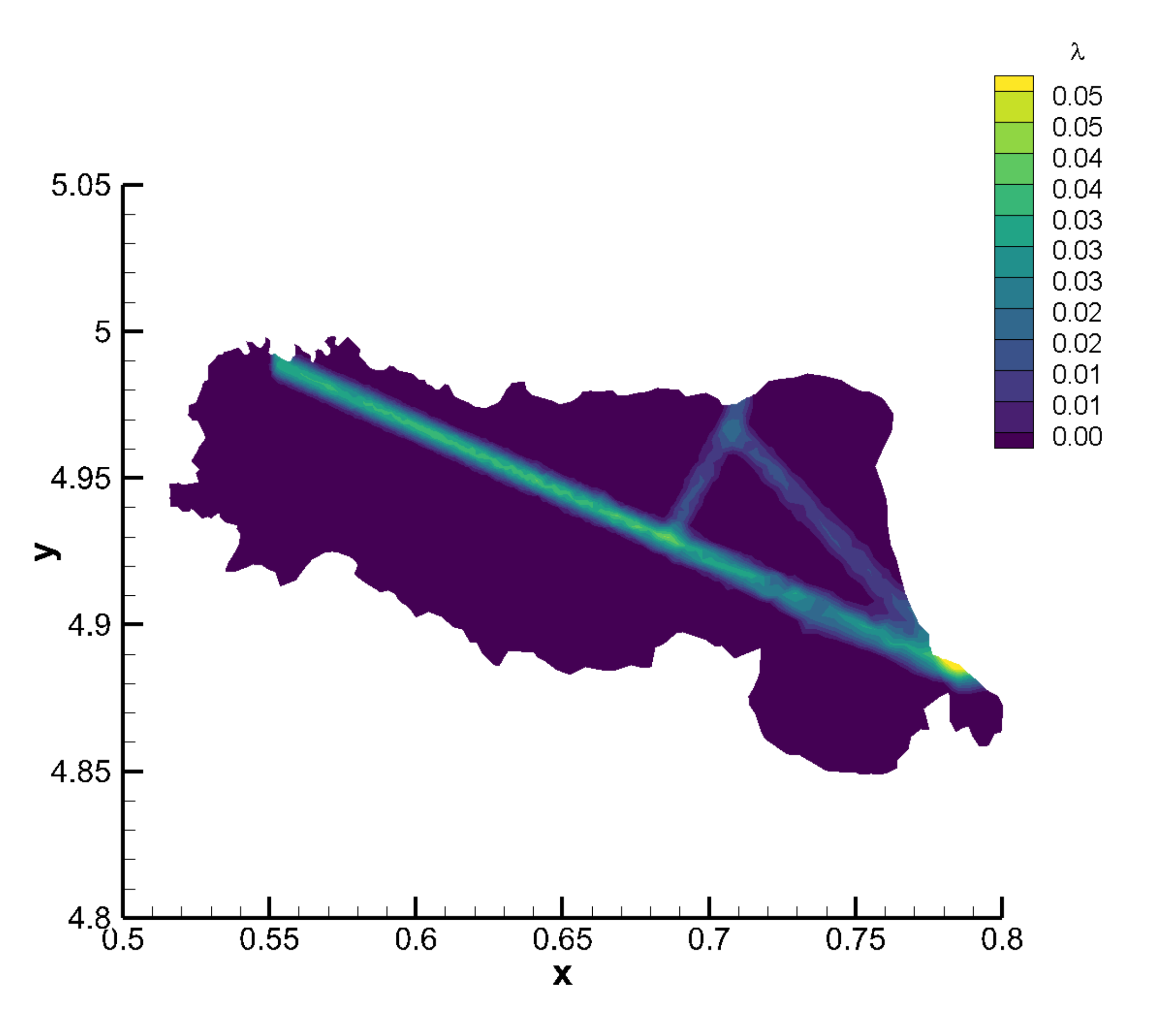} \\			
		\end{tabular} 
		\end{center}
		\caption{Top: map of the Emilia-Romagna region with main cities (left) and associated provinces (right). Bottom: unstructured computational mesh used to discretize the region (left) and initial condition for $\lambda$ (right).}
			\label{fig.test3_ini}
\end{figure}

The computational domain is defined in terms of the boundary that limits the region of Emilia-Romagna. This can be found in \cite{istat} as a list of georeferenced points in the ED50/UTM Zone 32N reference coordinate system. In order to avoid ill-conditioned reconstruction matrices and other related problems arising while dealing with big numbers in finite arithmetics, all coordinates are rescaled by a factor of $\alpha_x=10^{6}$. The computational grid is composed of a total number of triangular control volumes $N_E=5057$ and zero-flux boundary conditions are imposed everywhere, thus assuming that no exchange of population is present with the surrounding regions. The region is then subdivided into a total number of $\mathcal{N}_c=9$ provinces with the associated main cities, as depicted in Figure \ref{fig.test3_ini}.
 
The initial distribution of a generic population $f(x,y)$ is assigned to each main city denoted with subscript $c$ as a multivariate Gaussian function with the variance being the radius of the urban area, that is
\begin{equation}
f(x,y) = \frac{1}{2\pi r_c} \, e^{ -\frac{(x-x_c)^2+(y-y_c)^2}{2 r_c^2} } \, f_c,
\label{eq:ini_f}
\end{equation}
where $f_c$ is the number of the individuals in the population associated to the generic city. The radius of the city refers to the circular area of each province proportional to the order of magnitude of the population of the biggest city (Bologna). In this way, the integral over the computational domain of the initial population distribution exactly retrieves the quantity $f_c$ in \eqref{eq:ini_f}. The initial settings for each province of the Emilia-Romagna region are taken from \cite{prot_civile} and reported in Table \ref{tab.ER_city_ini} with $P_T$ representing the total inhabitants of each province of the region. 

\begin{table}[!htbp] 
	\caption{Initial data: name of the city, radius of the city $r_c$ measured in km, total population $P_T$, total number of infected ($I_T$) and exposed ($I_T$) people on 1/03/2020 in the region Emilia-Romagna.} 
	\renewcommand{\arraystretch}{1.05}
	\begin{center} 	
			\begin{tabular}{l|cccc} 				
				\multicolumn{1}{c|}{City}			& Radius ($r_c$) & Population ($P_T$) & Infected ($I_T$) & Exposed ($E_T$) \\
				& [km] & [person] & [person] & [person] \\
				\hline
				Piacenza		& $4.5$ & $2.87 \cdot 10^{5}$  & $174$ & $696$ \\
				Parma			& $6$   & $4.54 \cdot 10^{5}$  & $59$  & $236$ \\
				Reggio Emilia	& $6$   & $5.32 \cdot 10^{5}$  & $7$   & $28$ \\
				Modena			& $9$   & $7.07 \cdot 10^{5}$  & $24$  & $96$ \\
				Bologna			& $15$  & $10.18 \cdot 10^{5}$ & $2$   & $8$ \\
				Ferrara			& $4.5$ & $3.49 \cdot 10^{5}$  & $0$   & $0$ \\	 
				Ravenna			& $4.5$ & $3.89 \cdot 10^{5}$  & $2$   & $8$ \\
				Cesena			& $4.5$ & $3.98 \cdot 10^{5}$  & $1$   & $4$ \\
		        Rimini			& $3$   & $3.40 \cdot 10^{5}$  & $16$  & $64$ \\	
				\hline
			\end{tabular}
	\end{center}
	\label{tab.ER_city_ini}
\end{table}
The total infected $I_T$ and exposed $E_T$ populations are initialized with \eqref{eq:ini_f} according to the values reported in Table \ref{tab.ER_city_ini}, while the total susceptible population is initially prescribed by setting $S_T=P_T-I_T-E_T$. The part of the population which involves the commuters is assigned according to the mobility data of the region that can be found in \cite{mobilityER}. In particular, the total number of individuals that commute among the provinces of the region are counted, so that the commuters matrix reported in Table \ref{tab.ER_mobility-abs} can be constructed. We only consider the paths connecting each province to its direct neighbor provinces, thus people crossing more than one province in one day are not taken into account. The relative number $\mathcal{C}$ of all commuters within each province normalized by the population of the city of origin is then easily evaluated and reported in the last column of Table \ref{tab.ER_mobility-abs} as a percentage. These percentages are then adopted for setting the commuters populations $S$, $E$, $I$ and $R$. Specifically, each control volume is assigned the total percentage of commuters referred to the province where it is located. For instance, in the province of Ferrara we set
\begin{equation}
S = 0.068 \cdot S_T, \qquad E = 0.068 \cdot E_T, \qquad I = 0.068 \cdot I_T.
\end{equation}

\begin{table}[!htbp] 
	\caption{Matrix of commuters among the provinces of the region Emilia-Romagna. Departure provinces are listed on the first left column, while arrival provinces are reported in columns (PC=Piacenza, PR=Parma, RE=Reggio Emilia, MO=Modena, BO=Bologna, FE=Ferrara, RA=Ravenna, FC=Cesena, RN=Rimini).  } 
	\begin{small}
	\begin{center} 	
		\begin{tabular}{c|ccccccccc||c} 				
			\multicolumn{1}{c|}{From \textbackslash To}	& PC & PR & RE & MO & BO & FE & RA & FC & RN & $\mathcal{C}$ [\%] \\
			\hline
			PC		& -      & $4178$ 	& - & - & - & - & - & - & - & 1.45 \\
			PR			& $1707$ & -      	& $5142$ & - & - & - & - & - & - & 1.51 \\
			RE	& -      & $8969$ 	& - & $19841$ & - & - & - & - & - & 5.42 \\
			MO		& -      & -      	& $11488$ & - & $13034$ & $1173$ & - & - & - & 3.63 \\
			BO		& -      & - 		& - & $6842$ & - & $5983$ & $3887$ & - & - & 1.64 \\
			FE		& -      & - 		& - & $2682$ & $16865$ & - & $2610$ & - & - & 6.80 \\	 
			RA		& -      & - 		& - & - & $9808$ & $1016$ & - & $9211$ & - & 5.14\\
			FC		& -      & - 		& - & - & & - & $6646$ & - & $6944$ & 3.41\\
			RN		& -      & - 		& - & - & & - & - & $6075$ & - & 1.79\\	
			\hline
		\end{tabular}
	\end{center}
    \end{small}
	\label{tab.ER_mobility-abs}
\end{table}
Then, the static populations $(\SO,\EO,\IO)$ which are not moving away from the cities are computed relying on conservation principles. All values of populations are normalized by a factor of $\alpha_p=10^5$ and the recovered population is initially set to zero. 
The contact rates have been estimated from available data on the entire Emilia-Romagna region using the corresponding zero-order SEIR model and are set to $\beta_I=\beta_E=3.7 \cdot 10^{-3}$, while the recovery rate and the incubation period are taken from established values in the literature \cite{Veneziani2021} and read $\gamma_I=\gamma_E=1/12$ and $a=1/7$, respectively. The incidence functions \eqref{eq:incf}-\eqref{eq:incfE} are adopted with $p=1$, $\kappa_I=6 \cdot 10^{7}$ and $\kappa_E=2.3 \cdot 10^{4}$. The values $\kappa_I$ and $\kappa_E$ have been chosen in agreement with the zero-order model at a regional level. The higher value of $\kappa_I$ allows to mimic the effect of quarantines and social distancing with respect to infected people. 
We finally set $\sigma=\zeta=0.25$ according to \cite{Gatto}, hence we estimate the initial number of exposed individuals, including asymptomatic, as $E_T=4 \, I_T$, thus obtaining an initial reproduction number $R_0\approx 2.3$ that is in accordance with available literature \cite{FBK2020}.  

The propagation speed has been selected in order to match the overall epidemic data in the different provinces, and has been fixed to $\lambda^2=6.25 \cdot 10^{-4}$ for susceptible, exposed and recovered populations. The same speed is prescribed along the main connections of the region, namely the highway Piacenza-Rimini, the connection between Ferrara and Rimini, passing through Ravenna, as well as the minor path joining Ferrara to Bologna, as shown in the last panel of Figure \ref{fig.test3_ini}. 
A value of $\lambda^2=10^{-12}$ is set in the rest of the computational domain and each connection band has a width of $h=1.5$ km. The propagation speed of the infected population is set to zero, i.e. $\lambda_I^2=0$, meaning that all individuals who have been detected as infected are not allowed to move according to the governement restrictions. However, the infected people, even if limited by quarantines and social distancing, can still contribute to the spread the disease via the diffusion process at the urban scale (mimicking for instance the still possible infections happening at the family level). The relaxation time is $\tau_r=10^4$ so that the model recovers a hyperbolic regime in the entire region, apart from the main cities, where a parabolic setting is prescribed in order to correctly capture the diffusive behavior of the disease spreading which typically occurs in highly urbanized zones. Therefore, the relaxation time $\tau$ is  prescribed as follows:
\begin{equation}
\tau = \tau_r + (\tau_0-\tau_r) \, \sum \limits_{c=1}^{\mathcal{N}_c} e^{-\frac{1}{2} \frac{\left( (x-x_c)^2+(y-y_c)^2\right) }{r_c^2} },
\end{equation}
with $(x_c,y_c)$ representing the coordinates of a generic city center and $r_c$ the associated radius defined in Table \ref{tab.ER_city_ini}. The diffusion relaxation time is chosen to be $\tau_0=10^{-4}$. 

\begin{figure}[htbp]
	\begin{center}
		\begin{tabular}{cc} 
			\includegraphics[width=0.45\textwidth]{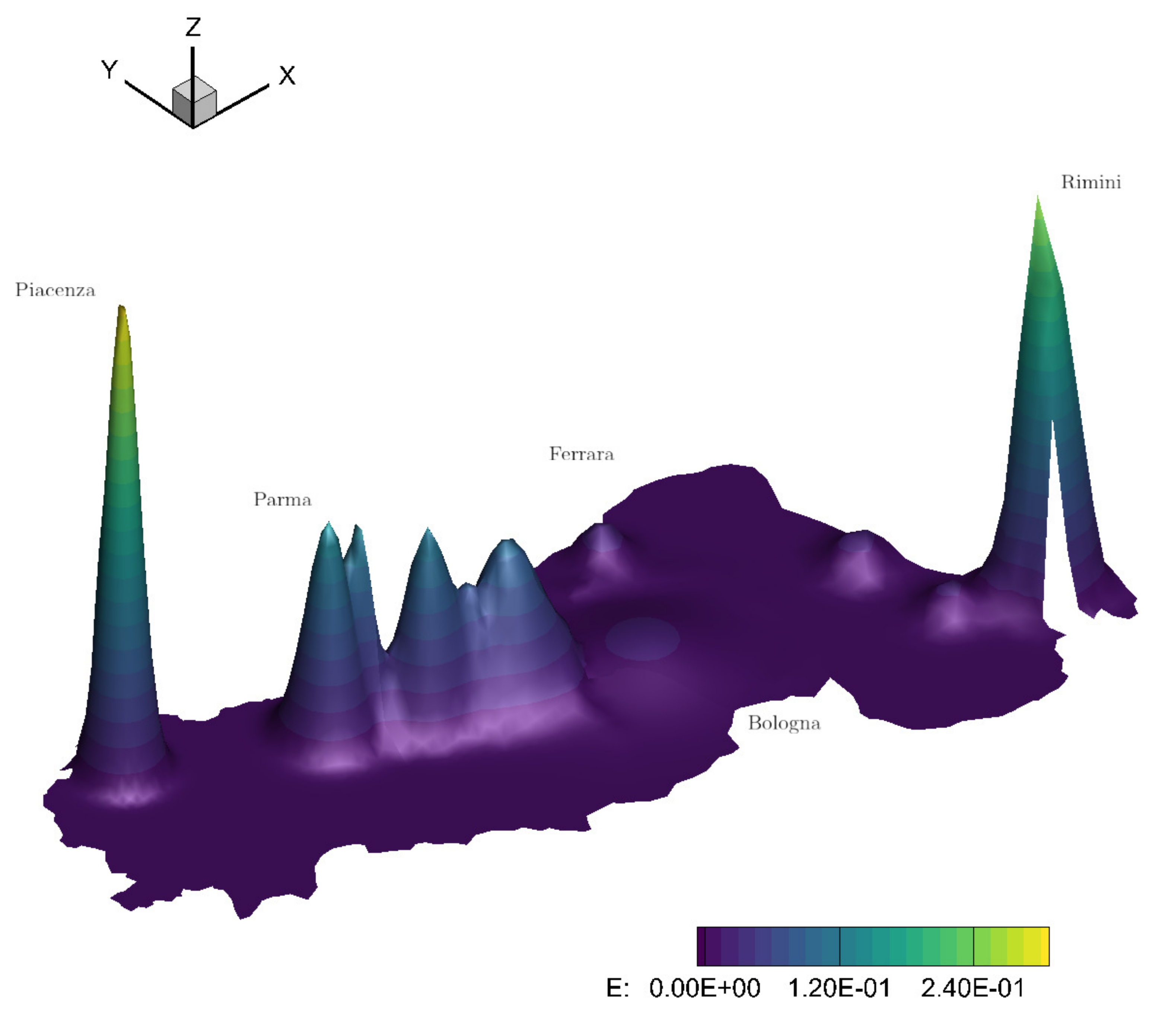} & 
			\includegraphics[width=0.45\textwidth]{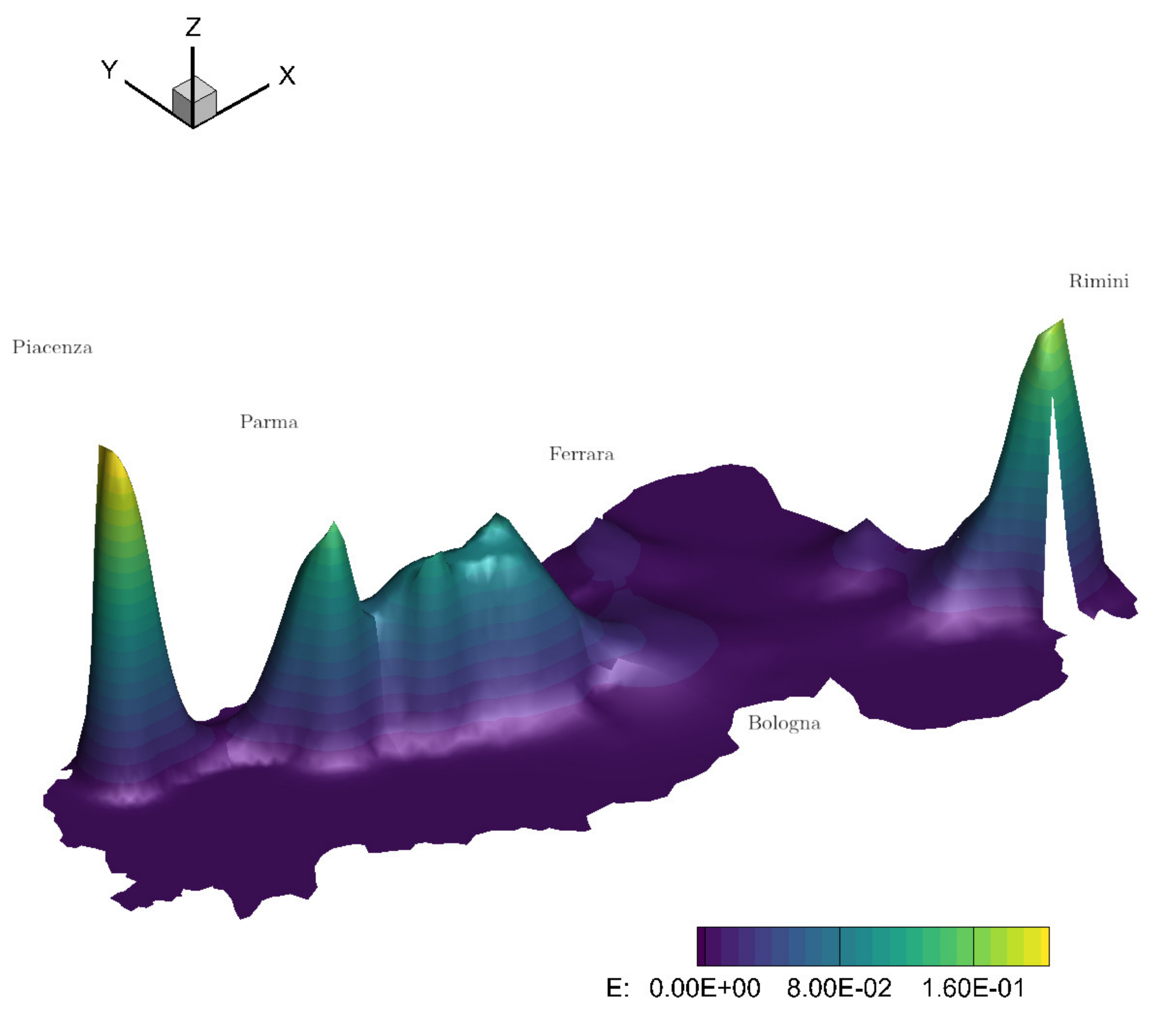} \\
			\includegraphics[width=0.45\textwidth]{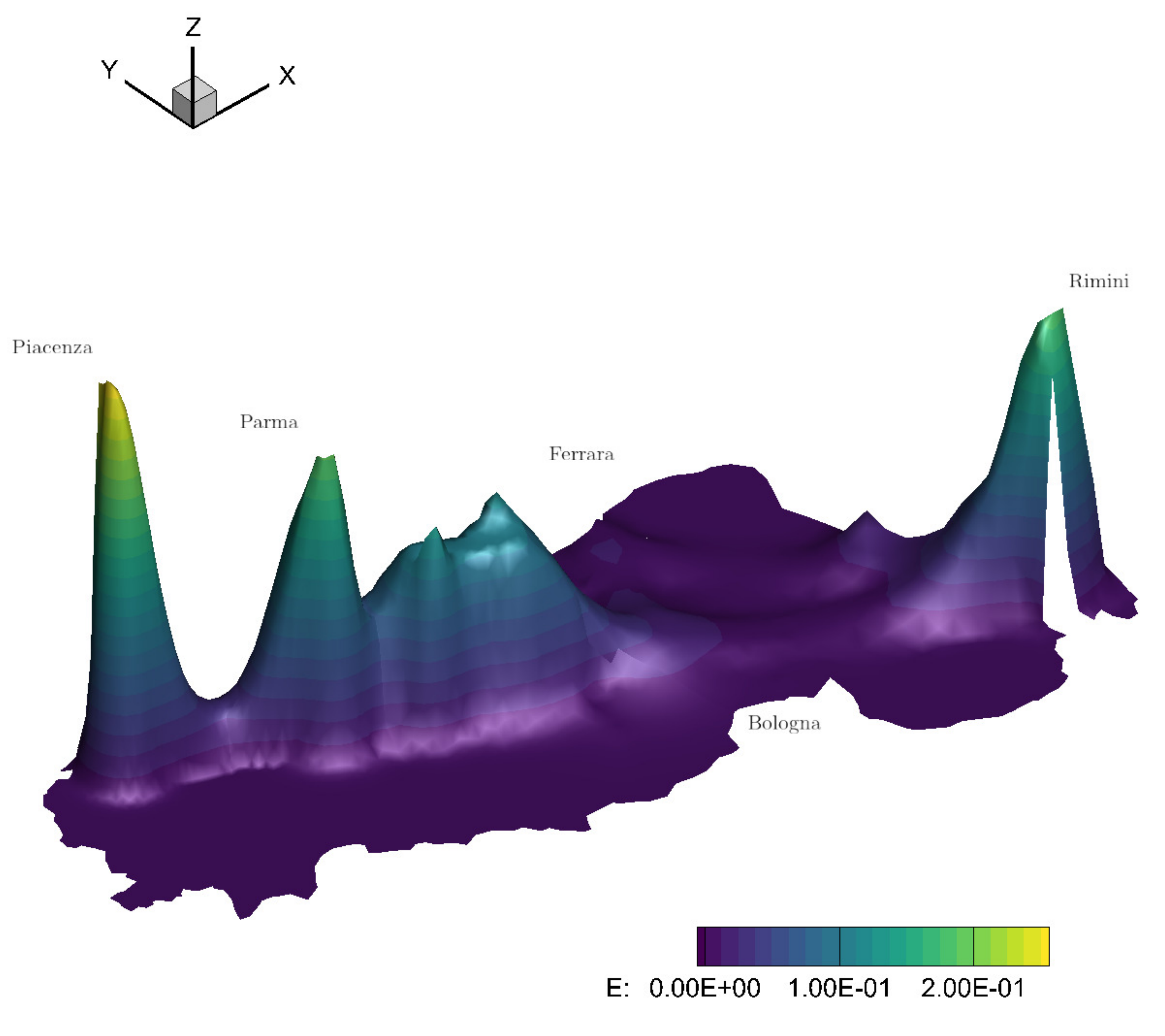} & 
            \includegraphics[width=0.45\textwidth]{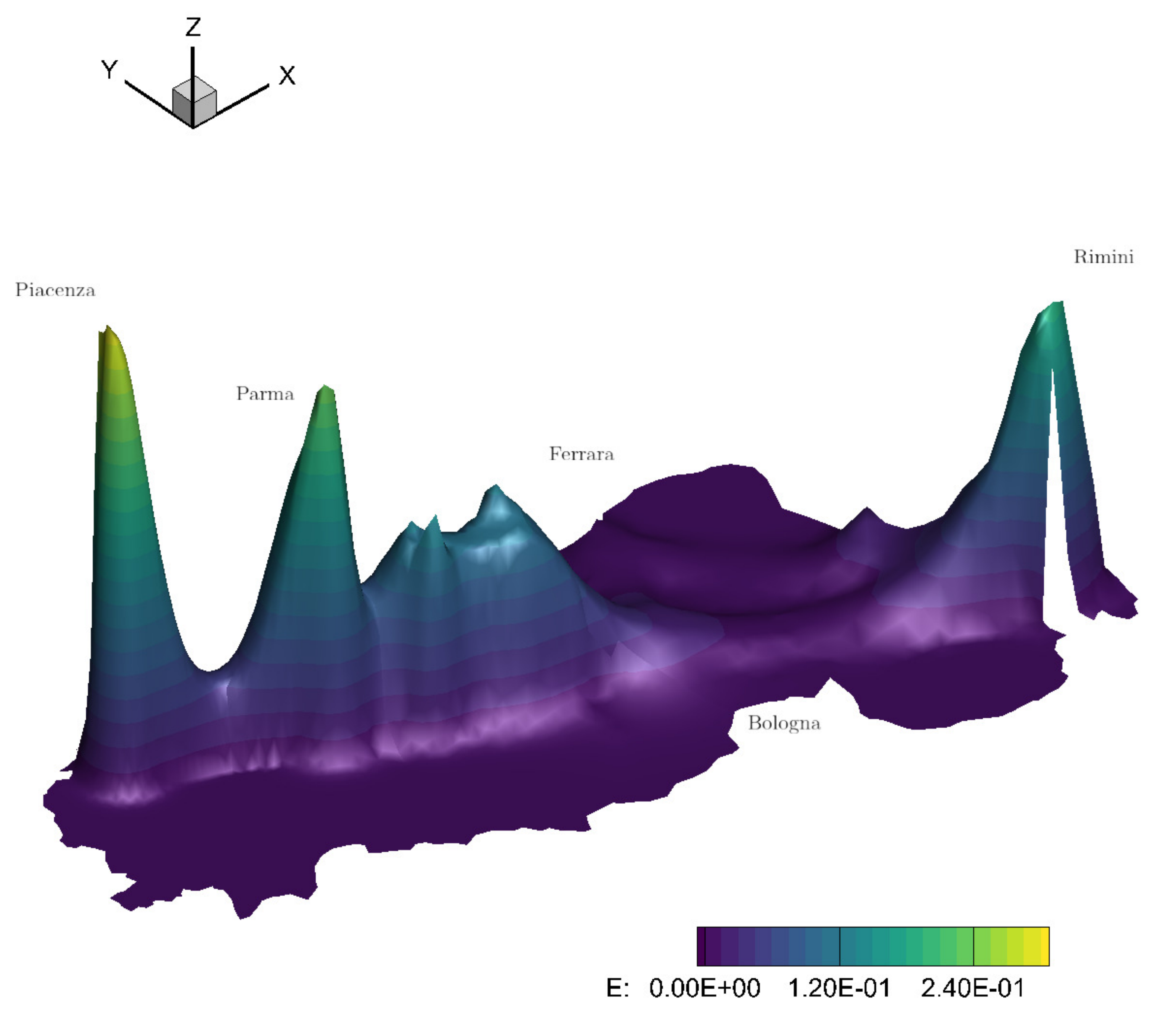} \\	
            \includegraphics[width=0.45\textwidth]{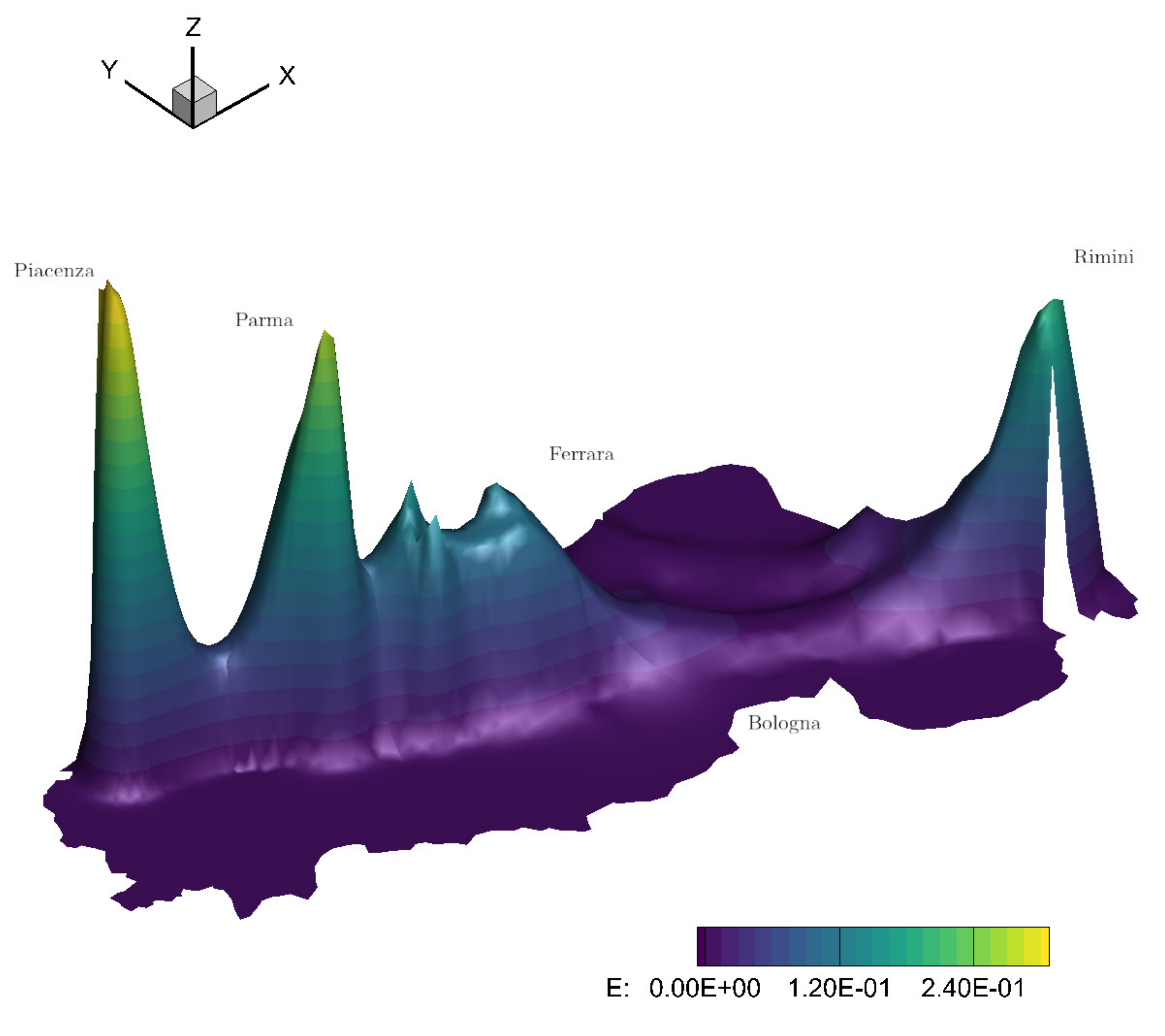} & 
            \includegraphics[width=0.45\textwidth]{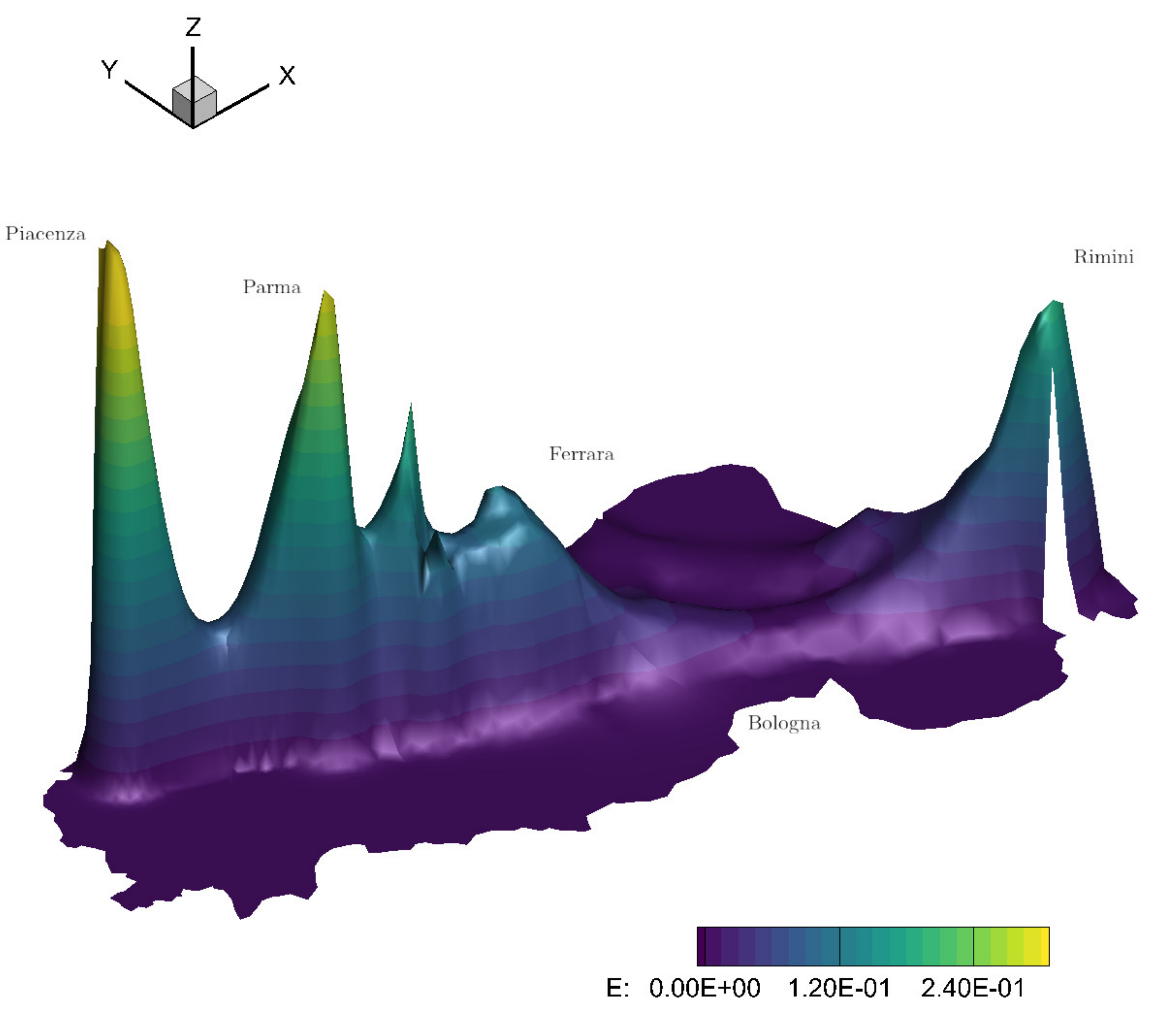} \\		
        \end{tabular} 
        	\end{center}
		\caption{Distribution of exposed population $E$ including asymptomatic at times $t=0$, $t=2$, $t=4$, $t=6$, $t=8$  and $t=10$ (from top left to bottom right).}
	\label{fig.test3_E3D}
\end{figure}

The final time is chosen to be $t=10$, therefore a total number of ten days is simulated. The south-west part of the region is not assigned any population distribution nor any propagation speed since it is mainly covered by the Appennini mountains and almost zero circulation of people is observed.
The units of measure used in our simulations can be conveniently summarized as follows:
\begin{equation}
1 \, \textnormal{km} = 10^{-3} \, \textnormal{L}, \qquad 1 \, \textnormal{person} = 10^{-5} \, \textnormal{P}, \qquad 1 \,  \textnormal{day} = 2 \textnormal{T},
\end{equation}
with [L], [P] and [T] being the length, person and time units used in the simulation, respectively.

\begin{figure}[!htbp]
	\begin{center}
		\begin{tabular}{cc} 
			\includegraphics[width=0.45\textwidth]{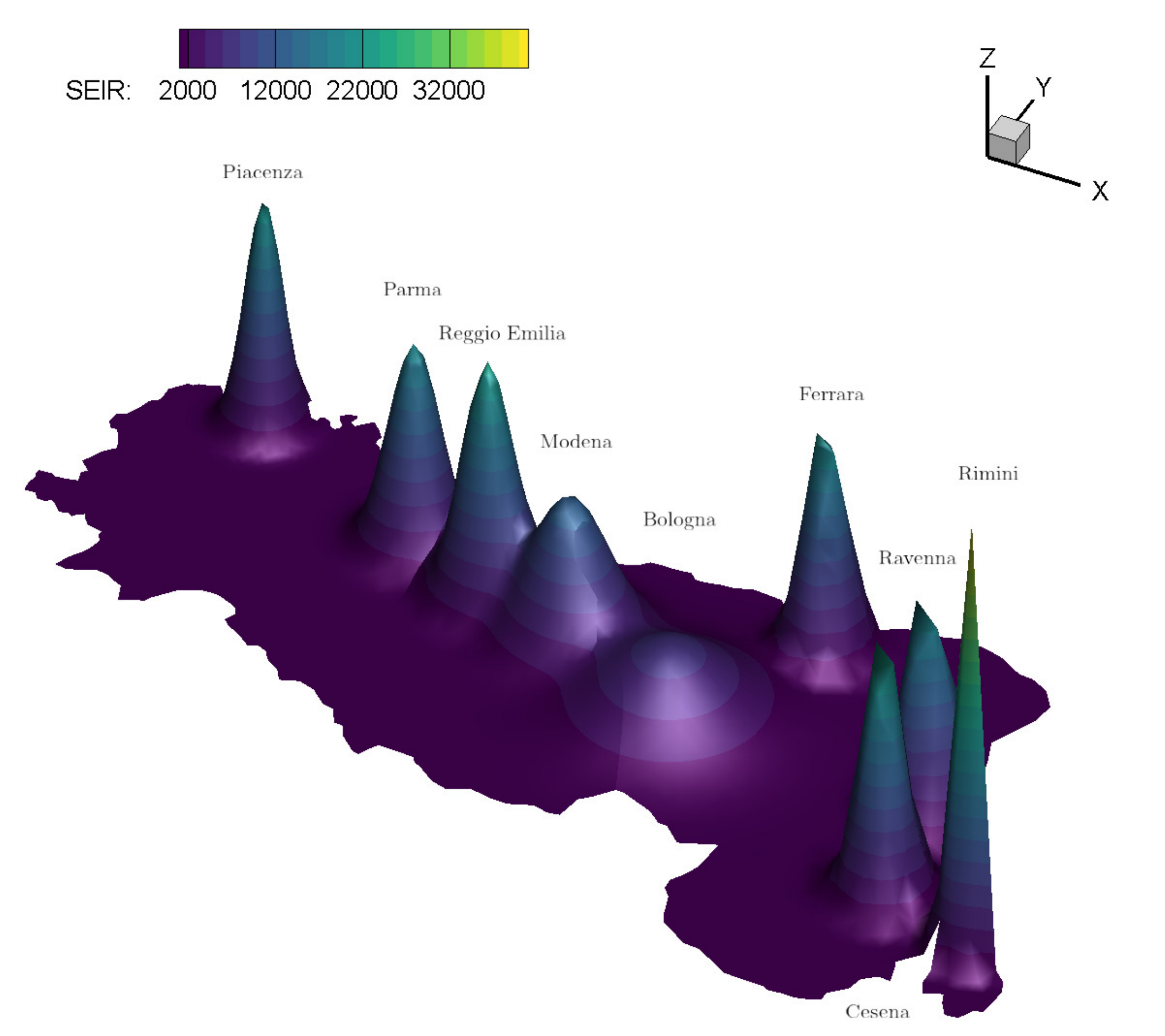} & 
			\includegraphics[width=0.45\textwidth]{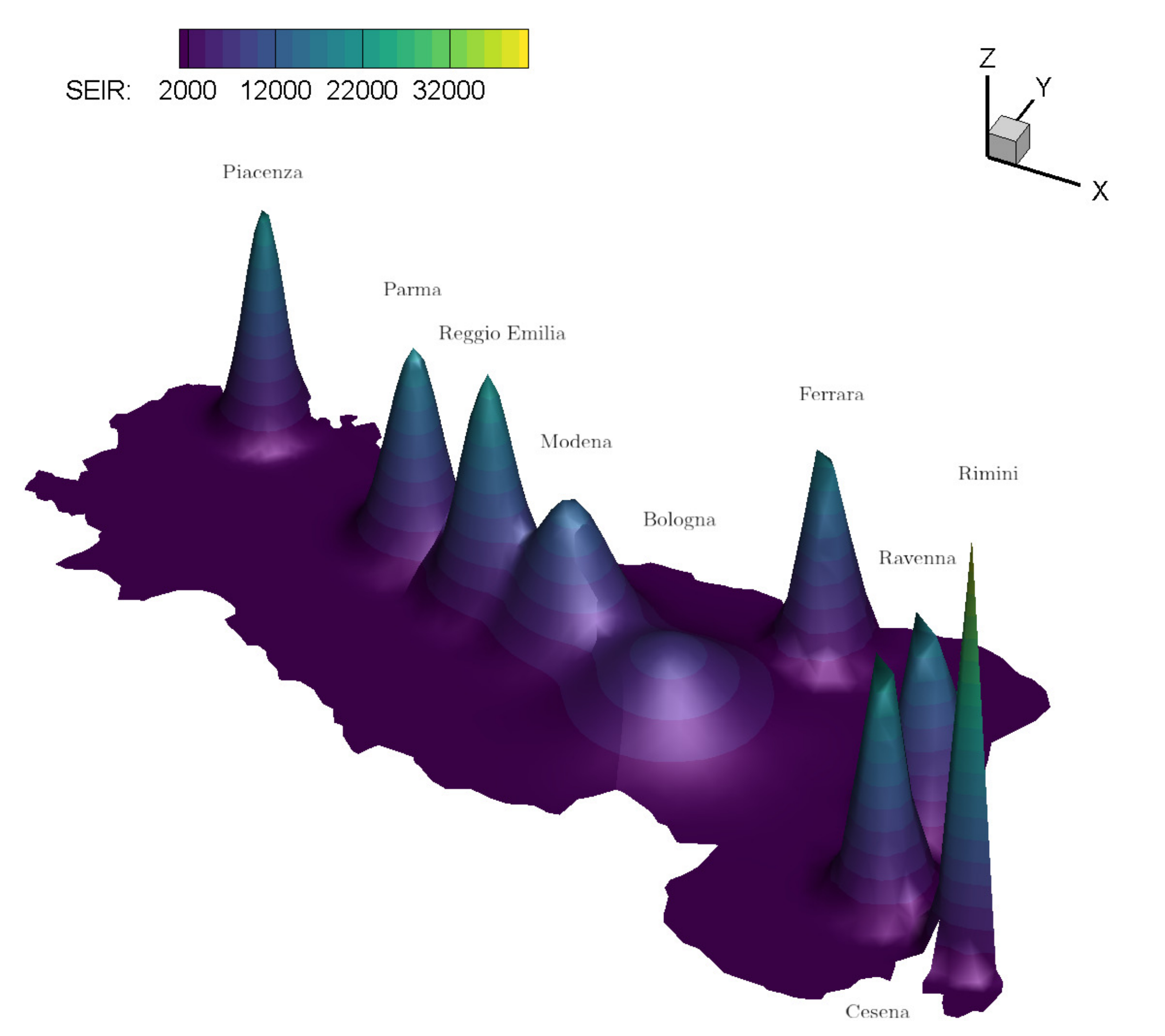} \\
		\end{tabular} 
		\end{center}
		\caption{Distribution of total population $(S+\SO)+ (E+\EO) + (I+\IO) + (R+\RO)$ at initial time $t=0$ (left) and at the final time $t=10$ (right).}
	\label{fig.test3_SEIR3D}
\end{figure}

Figure \ref{fig.test3_E3D} shows the time evolution of the exposed population $E$, including asymptomatic, which is moving from both Piacenza and Rimini towards the center of the region and the city of Bologna, then spreading northern in the direction of Ferrara. The wave of the exposed population is clearly visible, highlighting the hyperbolic regime of the model. Figure \ref{fig.test3_SEIR3D} depicts the total population at the initial and final time, including the non commuters $(\SO,\EO,\IO,\RO)$ that remain at rest in the cities and are only affected by a diffusive process modeled by the relaxation time $\tau_0$. One can notice that the total population of the cities does not change in time, since only $1.45\%-6.8 \%$ of the individuals are moving, which corresponds to the real scenario defined by Table \ref{tab.ER_mobility-abs}. Nevertheless the disease is spreading over the entire region due to the people who daily commute from one city to another one. 

\begin{figure}[!htbp]
	\begin{center}
		\begin{tabular}{cc} 
			\includegraphics[width=0.47\textwidth]{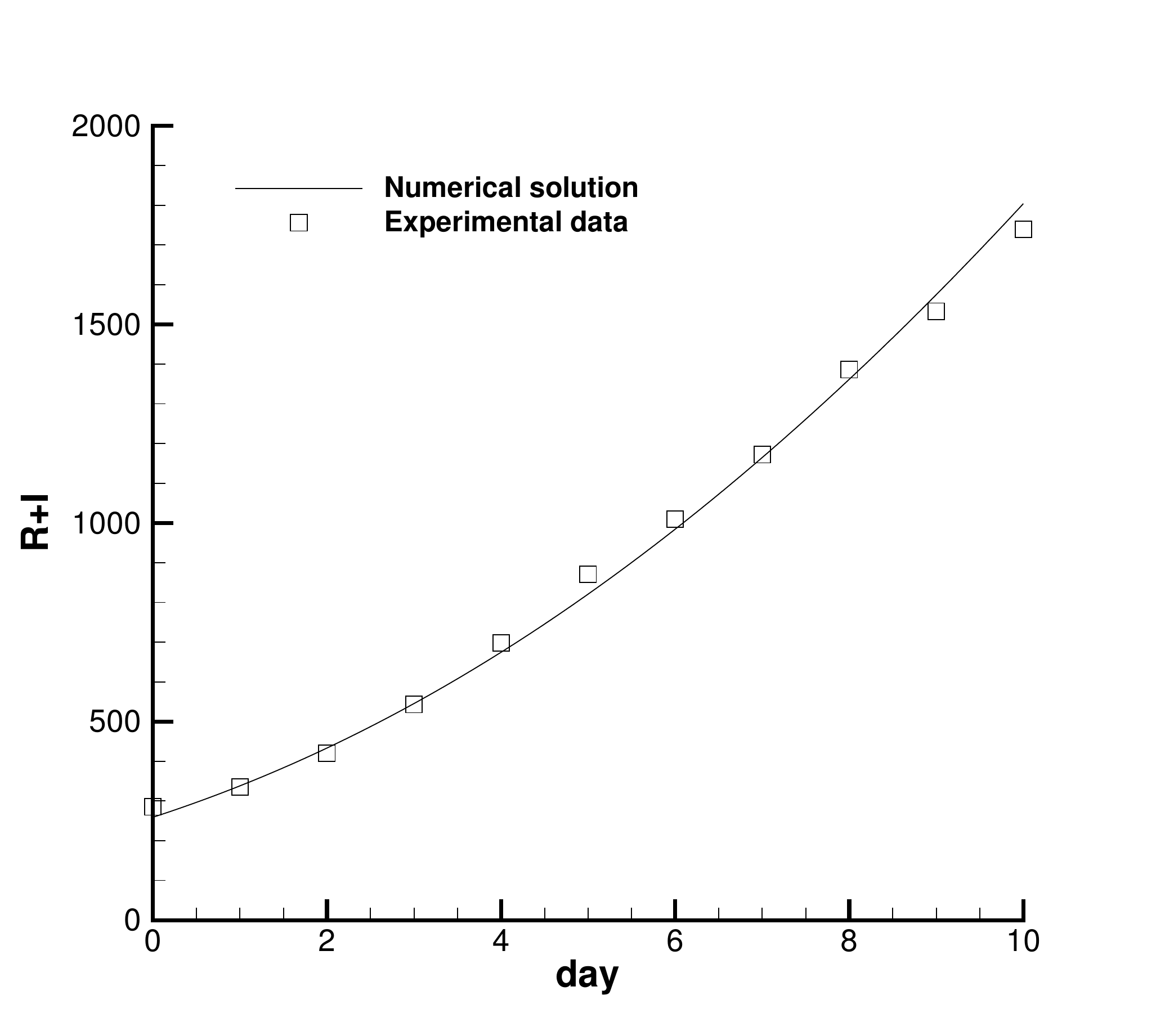} &
			\includegraphics[width=0.47\textwidth]{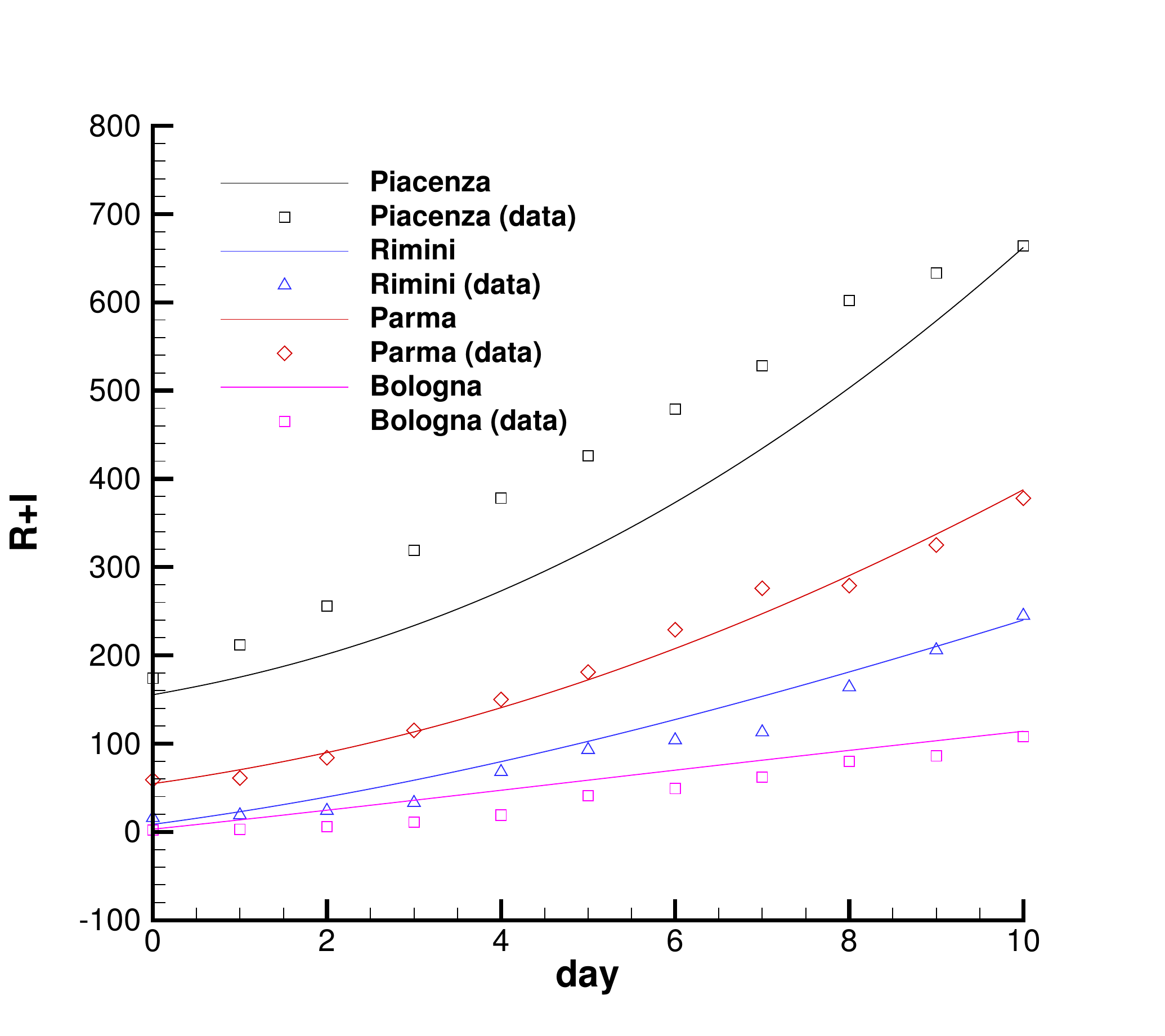} \\    		
		\end{tabular} 
	\end{center}
	\caption{Left: time evolution of total infected and recovered population ($R+I$) compared against experimental data for the region Emilia-Romagna. Right: time evolution of total infected and recovered population ($R+I$) compared against experimental data for the province of Piacenza (black), Parma (red), Bologna (purple) and Rimini (blue).}
	%\caption{Top left: time evolution of the susceptible ($S$), exposed ($E$), infected ($I$) and recovered ($R$) population. Top right: zoom of the time evolution of total exposed ($E$) population.	Bottom left: time evolution of total infected and recovered population ($R+I$) compared against experimental data for the region Emilia-Romagna. Bottom right: time evolution of total infected and recovered population ($R+I$) compared against experimental data for the province of Piacenza (black), Parma (red), Bologna (purple) and Rimini (blue).}
	\label{fig.test3_SEIR}
\end{figure}

\begin{figure}[!htbp]
	\begin{center}
		\begin{tabular}{c}  		
			\includegraphics[width=0.5\textwidth]{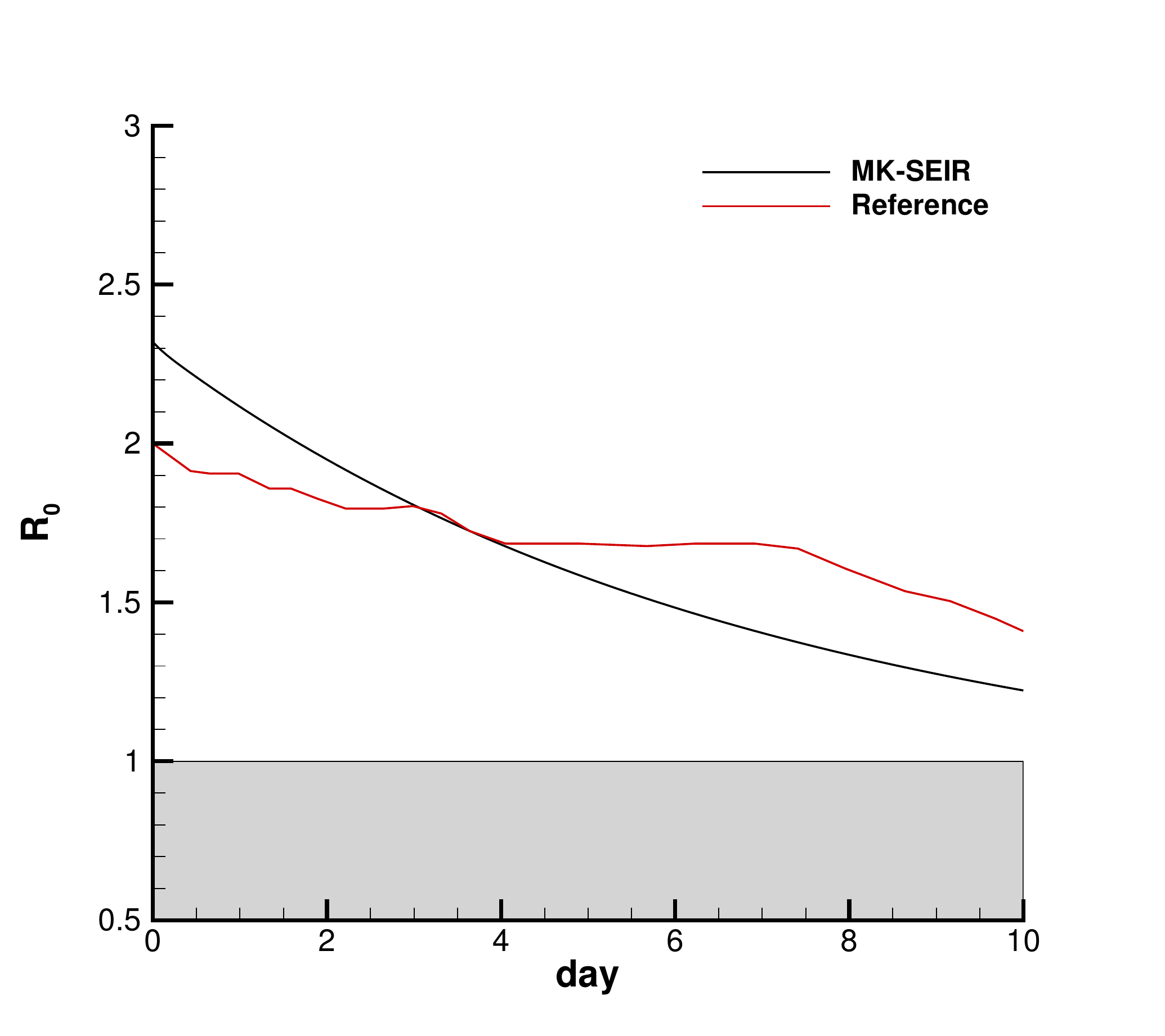}
		\end{tabular} 
	\end{center}
	\caption{Time evolution of the reproduction number at the regional level compared to the study presented in \cite{FBK2020}.}
	\label{fig.test3_R0}
\end{figure}

%Finally, the time evolution of the populations $(S,E,I,R)$ is reported in Figure \ref{fig.test3_SEIR} as well as a detailed view of the growth of the exposed population. 
%\textcolor{red}{The initial decreasing number of the infected people is due to the model initialization which starts the incubation period at the beginning of the simulation, thus a certain amount of time is needed to mimic and recover the real ongoing setting that have started before the initial time of the simulation, which is actually recovered in the second part of the simulation where the number of infected individuals grows again.} 
Finally, Figure \ref{fig.test3_SEIR} plots a comparison against the measured data at the regional level which can be found in \cite{prot_civile} and the same comparison for the province of Piacenza, Parma, Bologna and Rimini.  The time evolution of the reproduction number $R_0$ at the regional level is depicted in Figure \ref{fig.test3_R0}, which qualitatively recovers the results obtained in \cite{FBK2020}. The reproduction number has been monitored with \eqref{eq:R0} and is beyond the threshold of $R_0=1$, thus indicating that the epidemic spreading is growing. An overall very good agreement with experimental measurements can be appreciated both at the regional and at the provincial level, thus demonstrating the capability of the novel model to adapt to real world settings and applications.

%\begin{figure}[!htbp]
%	\begin{center}
%		\begin{tabular}{cc} 
%			\includegraphics[width=0.47\textwidth]{test3_EIR} & 
%			\includegraphics[width=0.47\textwidth]{test3_I} \\
%			\includegraphics[width=0.47\textwidth]{test3_RI_region} &
%			\includegraphics[width=0.47\textwidth]{test3_RI_province} \\    		
%		\end{tabular} 
%	\end{center}
%	\caption{Top left: time evolution of the exposed ($E$), infected ($I$) and recovered ($R$) population. Top right: zoom of the time evolution of total infected ($I$) population.	Bottom left: time evolution of total infected and recovered population ($R+I$) compared against experimental data for the region Emilia-Romagna. Bottom right: time evolution of total infected and recovered population ($R+I$) compared against experimental data for the province of Piacenza (black), Parma (red), Bologna (purple) and Rimini (blue).}
%		\label{fig.test3_SEIR}
%\end{figure}

%------------------------------------------------------------------------------------
%------------------------------------------------------------------------------------
\section{Conclusions}
In this work, we introduced new multiscale kinetic models for the description of the spread of an epidemic disease in a spatially heterogeneous context. The models distinguish between two populations of individuals, a population of commuters moving over long distances (extra-urban) and a population of non commuters acting over short distances (urban). Commuter dynamics are described by a kinetic transport equation, while non commuter dynamics follow a classical diffusion process limited to urban areas. From the epidemiological point of view the model was first presented in a simple SIR compartmental framework and then extended to the case of a SEIR compartmental structure designed to take into account the characteristics of COVID-19. In an appropriate scaling limit, the commuters can be described by traditional diffusive models. This is indeed the scale considered for this portion of the population when approaching an urban center that allows commuters to be consistently described as non-commuters within these highly populated regions.

The proposed models were solved using ad hoc developed finite volume numerical methods acting on unstructured grids capable of effectively solving the kinetic model at small space-time scales. This allowed us to present several numerical results showing the ability of the kinetic model to avoid the unrealistic features of traditional diffusive models based on a single population, such as infinite propagation speed and indistinct movement of the entire population. 
%A numerical application to a realistic geographical scenario has shown the capability of the model to describe correctly the initial phases of the epidemic spread of COVID-19 in a region of northern Italy. 
A last part of the work is dedicated to a careful simulation of the first days of the COVID-19 epidemic in a realistic geographical scenario.

%Emilia-Romagna the second most affected region of Italy during the first epidemic wave of 2020. To our knowledge these are the first space dependent results concerning the epidemic spread of COVID-19 in this region.

In perspective, we would like to study the effectiveness of the model in the case of more complex compartmental models that take into account, for example, hospitalized data and mortality \cite{Gatto, Tang}. Moreover, to make the model more effective for decision-makers, the inclusion of suitable control processes describing lockdown limitations and the incorporation of an age-structured population is crucial to correctly describe the impact of specific kinds of infectious diseases, like the COVID-19 pneumonia \cite{Colombo, HWH00}.
Finally, since data of the spread of epidemics are generally highly heterogeneous and affected by a great deal of uncertainty, future perspectives include the application of uncertainty quantification methods to assess the impact of stochastic inputs in the proposed multiscale kinetic transport model \cite{APZ, JHL}. 
 
%------------------------------------------------------------------------------------
%------------------------------------------------------------------------------------
\section*{Acknowledgments} 
This work has been written within the
activities of GNCS groups of INdAM (National Institute of
High Mathematics). The support of MIUR-PRIN Project 2017, No. 2017KKJP4X “Innovative numerical methods for evolutionary partial differential equations and applications” is acknowledged.

%\begin{figure}
%\begin{center}
%\begin{tabular}{ccc} 
%\includegraphics[scale=0.35]{gauss1.eps}&
%\includegraphics[scale=0.35]{gauss2.eps}&
%\includegraphics[scale=0.35]{gauss4.eps}\\
%\end{tabular}
%\end{center}
%\caption{The four quadrants in the even and odd parities method and the corresponding Gaussian quadrature directions for various numbers of points.}
%\label{fig:gauss}
%\end{figure}

\appendix
\section{An asymptotic-preserving method on unstructured grids}
In this appendix, we report the details of the numerical method used for the discretization of the multiscale kinetic system described in the main part of the paper. To avoid unessential difficulties we illustrate the numerical method in the case of the simpler kinetic model \eqref{eq:kineticc}-\eqref{eq:diffuse}, i.e. the MK-SIR model. The numerical scheme combines a discrete ordinate method in velocity with the even and odd parity formulation \cite{DP,JPT} and achieves asymptotic preservation in time using suitable IMEX Runge-Kutta schemes \cite{Bos1, Bos2}. Namely, to obtain a scheme which consistently captures the diffusion limit and for which the choice of the time discretization step is not related to the smallness of the scaling parameters $\tau_{S,I,R}$. Next, we summarize the key ingredients used to discretize the space variables on a two-dimensional unstructured mesh \cite{BD-BGK,ArepoTN} which permits to deal with realistic geometries.

%------------------------------------------------------------------------------------
\subsection{Even and odd parities formulation}
From a computational viewpoint, it is convenient to rewrite \eqref{eq:kineticc} by 
splitting it into four parts according to the quadrants of the velocity space % (see Figure \ref{fig:gauss}) 
to which velocities belong. Let us denote $v=(\eta,\xi)\in\mathbb{S}^1$, we obtain four equations with non-negative $\xi, \eta \geq 0$. We then define the even and odd parities \cite{JPT}
\[
\begin{split}
r^{(1)}_S(\xi,\eta) &= \frac12(f_S(\xi,-\eta)+f_S(-\xi,\eta)),\quad 
r^{(2)}_S(\xi,\eta) = \frac12(f_S(\xi,\eta)+f_S(-\xi,-\eta))\\
r^{(1)}_I(\xi,\eta) &= \frac12(f_I(\xi,-\eta)+f_I(-\xi,\eta)),\quad 
r^{(2)}_I(\xi,\eta) = \frac12(f_I(\xi,\eta)+f_I(-\xi,-\eta))\\
r^{(1)}_R(\xi,\eta) &= \frac12(f_R(\xi,-\eta)+f_R(-\xi,\eta)),\quad 
r^{(2)}_R(\xi,\eta) = \frac12(f_R(\xi,\eta)+f_R(-\xi,-\eta))
\end{split}
\]
and the scalar fluxes
\[
\begin{split}
j^{(1)}_S(\xi,\eta) &= \frac{\vs}{2} (f_S(\xi,-\eta)+f_S(-\xi,\eta)),\quad 
j^{(2)}_S(\xi,\eta) = \frac{\vs}{2}  (f_S(\xi,\eta)+f_S(-\xi,-\eta))\\
j^{(1)}_I(\xi,\eta) &= \frac{\vi}{2}  (f_I(\xi,-\eta)+f_I(-\xi,\eta)),\quad 
j^{(2)}_I(\xi,\eta) = \frac{\vi}{2}  (f_I(\xi,\eta)+f_I(-\xi,-\eta))\\
j^{(1)}_R(\xi,\eta) &= \frac{\vr}{2}  (f_R(\xi,-\eta)+f_R(-\xi,\eta)),\quad 
j^{(2)}_R(\xi,\eta) = \frac{\vr}{2}  (f_R(\xi,\eta)+f_R(-\xi,-\eta))
\end{split}
\]
An equivalent formulation with respect to \eqref{eq:kineticc} then reads as 
\begin{equation}
\begin{split}
\frac{\partial r^{(1)}_S }{\partial t} +  \xi \frac{\partial j^{(1)}_S}{\partial x}-\eta \frac{\partial j^{(1)}_S}{\partial y} &= -F(r^{(1)}_S, I_T) +\frac1{\tau_S}\left(S-r^{(1)}_S \right)\\
\frac{\partial r^{(2)}_S }{\partial t} +  \xi \frac{\partial j^{(2)}_S}{\partial x}+\eta \frac{\partial j^{(2)}_S}{\partial y} &= -F(r^{(2)}_S, I_T)+\frac1{\tau_S}\left(S-r^{(2)}_S \right)\\
\frac{\partial r^{(1)}_I }{\partial t} +  \xi \frac{\partial j^{(1)}_I}{\partial x}-\eta \frac{\partial j^{(1)}_I}{\partial y} &= F(r^{(1)}_S, I_T)-\gamma r^{(1)}_I+\frac1{\tau_I}\left(I-r^{(1)}_I \right)\\
\frac{\partial r^{(2)}_I }{\partial t} +  \xi \frac{\partial j^{(2)}_S}{\partial x}+\eta \frac{\partial j^{(2)}_I}{\partial y} &= F(r^{(2)}_S, I_T)-\gamma r^{(2)}_I+\frac1{\tau_I}\left(I-r^{(2)}_I \right)\\
\frac{\partial r^{(1)}_R }{\partial t} +  \xi \frac{\partial j^{(1)}_R}{\partial x}-\eta \frac{\partial j^{(1)}_R}{\partial y} &= \gamma r^{(1)}_I+\frac1{\tau_R}\left(R-r^{(1)}_R \right)\\
\frac{\partial r^{(2)}_R }{\partial t} +  \xi \frac{\partial j^{(2)}_R}{\partial x}+\eta \frac{\partial j^{(2)}_R}{\partial y} &= \gamma r^{(2)}_I+\frac1{\tau_R}\left(R-r^{(2)}_R \right)
\end{split}
\label{eq:eop1}
\end{equation}
and 
\begin{equation}
\begin{split}
\frac{\partial j^{(1)}_S}{\partial t} + \vs^2 \xi \frac{\partial r^{(1)}_S}{\partial x} -\vs^2 \eta \frac{\partial r^{(1)}_S}{\partial y} &= - F(j^{(1)}_S, I_T) -\frac1{\tau_S} j^{(1)}_S\\
\frac{\partial j^{(2)}_S}{\partial t} + \vs^2 \xi \frac{\partial r^{(2)}_S}{\partial x} +\vs^2 \eta \frac{\partial r^{(2)}_S}{\partial y} &= - F(j^{(2)}_S, I_T) -\frac1{\tau_S} j^{(2)}_S\\
\frac{\partial j^{(1)}_I}{\partial t} + \vi^2 \xi \frac{\partial r^{(1)}_I}{\partial x} -\vi^2 \eta \frac{\partial r^{(1)}_I}{\partial y} &= F(j^{(1)}_I, I_T) -\gamma j^{(1)}_I -\frac1{\tau_I} j^{(1)}_I\\
\frac{\partial j^{(2)}_I}{\partial t} + \vi^2 \xi \frac{\partial r^{(2)}_I}{\partial x} +\vi^2 \eta \frac{\partial r^{(2)}_I}{\partial y} &= F(j^{(2)}_I, I_T) - \gamma j^{(2)}_I -\frac1{\tau_I} j^{(2)}_I\\
\frac{\partial j^{(1)}_R}{\partial t} + \vr^2 \xi \frac{\partial r^{(1)}_R}{\partial x} -\vr^2 \eta \frac{\partial r^{(1)}_R}{\partial y} &= \gamma j^{(1)}_R  -\frac1{\tau_R} j^{(1)}_R\\
\frac{\partial j^{(2)}_R}{\partial t} + \vr^2 \xi \frac{\partial r^{(2)}_R}{\partial x} +\vr^2 \eta \frac{\partial r^{(2)}_R}{\partial y} &= \gamma j^{(2)}_R  -\frac1{\tau_R} j^{(2)}_R.
\end{split}
\label{eq:eop2}
\end{equation}
Note that, due to symmetry, we need to solve these equations for $\xi$, $\eta$ in the positive quadrant only. Thus the number of unknowns in \eqref{eq:kineticc} and \eqref{eq:eop1}-\eqref{eq:eop2} is effectively the same. Furthermore, setting for $\lambda \in [0,1]$ 
\begin{equation}
\xi = \cos\left(\frac{\lambda\pi}{2}\right),\qquad \eta = \sin\left(\frac{\lambda\pi}{2}\right)
\label{eq:lambda}
\end{equation}
we have
\begin{equation}
S=\frac12\int_{0}^1  (r^{(1)}_S+r^{(2)}_S)\,d\lambda,\quad
I=\frac12\int_{0}^1  (r^{(1)}_I+r^{(2)}_I)\,d\lambda,\quad
R=\frac12\int_{0}^1  (r^{(1)}_R+r^{(2)}_R)\,d\lambda.
\label{eq:dlambda}
\end{equation}
The above densities can be can be approximated by a Gauss-Legendre quadrature rule. This leads to a discrete velocity setting, usually referred to as the discrete ordinate method, where we approximate
\begin{equation}
S \approx S_M = \frac14 \sum_{i=1}^n w_i \left(r^{(1)}_S(\xi_i,\eta_i)+r^{(2)}_S(\xi_i,\eta_i)\right)
\label{eq:gauss} 
\end{equation}
and similarly for the other densities $I$ and $R$.
In \eqref{eq:gauss} we defined 
\[
\xi_i = \cos\left(\frac{(\zeta_i+1)\pi}{4}\right),\quad \eta_i = \sin\left(\frac{(\zeta_i+1)\pi}{4}\right)
\]
so that $w_i$ and $\zeta_i$ are the standard Gauss-Legendre quadrature weights and points in $[-1,1]$. 
%In Figure \ref{fig:gauss} we report the corresponding velocity directions over the circle for $n=1,2,4$ where $M=4n$ is the total number of velocities.

%------------------------------------------------------------------------------------
\subsection{Space discretization on unstructured grids}
We consider a spatial two-dimensional computational domain $\Omega$ which is discretized by a set of non overlapping polygons $P_i, i=1, \dots N_p$. The union of all elements is called the 
\textit{tessellation} $\mathcal{T}_{\Omega}$ of the domain $\Omega$ and can be expressed as  
\begin{equation}
\mathcal{T}_{\Omega} = \bigcup \limits_{i=1}^{N_p}{P_i}, 
\label{trian}
\end{equation}
where $N_p$ is the total number of elements contained in the domain. The mesh is conforming, thus each edge $\lambda$ of an element is always shared by two adjacent control volumes, apart from physical boundaries of the computational domain. Each element $P_i$ is allowed to exhibit an arbitrary number $N_{S_i}$ of edges $\lambda_{j,i}$, thus ranging from triangles to general polygonal shapes. The boundary of the cell is addressed with $\partial P_i$ and is then given by
\begin{equation}
\partial P_i = \bigcup \limits_{j=1}^{N_{S_i}}{\lambda_{ji}}, 
\label{dP}
\end{equation}
where $\lambda_{ji}$ is the edge shared by elements $P_i$ and $P_j$. Further details on the construction of a conforming polygonal tessellation can be found for instance in \cite{BD-BGK,ArepoTN}.

The governing equations are then discretized on the unstructured mesh by means of a finite volume scheme. 
Let the system of PDE be cast in the general form
\begin{equation}
\label{PDE}
\frac{\partial \Q}{\partial t} + \nabla \cdot \F(\Q) = \S(\Q), \qquad (x,y) \in \Omega \subset \mathds{R}^2, \quad t \in \mathds{R}_0^+, \quad \Q \in \Omega_{\Q} \subset \mathds{R}^\nu,     
\end{equation} 
where $\Q=(q_1,q_2,...,q_\nu)$ is the vector of conserved variables defined in the space of the admissible states $\Omega_{\Q} \subset \mathds{R}^\nu$, $\F(\Q)$ is the linear flux tensor and $\S(\Q)$ 
represents the stiff source term. More precisely, the multiscale kinetic SIR model \eqref{eq:kineticc}-\eqref{eq:diffuse} with system \eqref{eq:kineticc} written using the parities in the form \eqref{eq:eop1}-\eqref{eq:eop2} fits the formalism \eqref{PDE} by setting
\begin{equation*}
\Q = \left( r_S^{(1)} , \,  r_S^{(2)} , \,  r_I^{(1)} , \,  r_I^{(2)} , \,  r_R^{(1)} , \,  r_R^{(2)} , \,  j_S^{(1)} , \,  j_S^{(2)} , \,  j_I^{(1)} , \,  j_I^{(2)} , \,  j_R^{(1)} , \,  j_R^{(2)} , \,  \SO , \,  \IO , \,  \RO 
\right)^\top,
\end{equation*}
and
\begin{equation*}
\F = \left( \begin{array}{c} \xi \, j_S^{(1)} \\ \xi \, j_S^{(2)} \\ \xi \,  j_I^{(1)} \\ j_I^{(2)} \\ \xi \, j_R^{(1)} \\ \xi \,  j_R^{(2)} \\ \xi \, \lambda_S^2 \, r_S^{(1)} \\ \xi \,  \lambda_S^2 \, r_S^{(2)}  \\ \xi \,  \lambda_I^2 \, r_I^{(1)} \\ \xi \,  \lambda_I^2 \, r_I^{(2)}  \\ \xi \,  \lambda_R^2 \, r_R^{(1)}  \\ \xi \,  \lambda_R^2 \, r_R^{(2)} \\ -D_S^0 \, S_x \\ -D_I^0 \, I_x \\ -D_R^0 \, R_x \end{array} \right.  
\left. \begin{array}{c} -\eta \, j_S^{(1)}  \\ \eta \, j_S^{(2)}  \\ -\eta \,  j_I^{(1)}  \\ \eta \,  j_I^{(2)} \\ -\eta \,  j_R^{(1)} \\ \eta \,  j_R^{(2)}  \\ -\eta \, \lambda_S^2 \, r_S^{(1)}  \\ \eta \,  \lambda_S^2 \, r_S^{(2)}  \\ -\eta \,  \lambda_I^2 \, r_I^{(1)}  \\ \eta \,  \lambda_I^2 \, r_I^{(2)}  \\ -\eta \,  \lambda_R^2 \, r_R^{(1)}  \\ \eta \, \lambda_R^2 \, r_R^{(2)} \\ -D_S^0 \, S_y \\ -D_I^0 \, I_y \\ -D_R^0 \, R_y \end{array} \right), \qquad 
\S = \left( \begin{array}{c} -F(r^{(1)}_S, I_T) +\frac1{\tau_S}\left(S-r^{(1)}_S \right) \\ -F(r^{(2)}_S, I_T) +\frac1{\tau_S}\left(S-r^{(2)}_S \right) \\ F(r^{(1)}_S, I_T)-\gamma r^{(1)}_I+\frac1{\tau_I}\left(I-r^{(1)}_I \right) \\ F(r^{(2)}_S, I_T)-\gamma r^{(2)}_I+\frac1{\tau_I}\left(I-r^{(2)}_I \right) \\
\gamma r^{(1)}_I+\frac1{\tau_R}\left(R-r^{(1)}_R \right) \\ \gamma r^{(2)}_I+\frac1{\tau_R}\left(R-r^{(2)}_R\right) \\ - F(j^{(1)}_S, I_T) -\frac1{\tau_S} j^{(1)}_S \\
- F(j^{(2)}_S, I_T) -\frac1{\tau_S} j^{(2)}_S \\
F(j^{(1)}_I, I_T) -\gamma j^{(1)}_I -\frac1{\tau_I} j^{(1)}_I \\
F(j^{(2)}_I, I_T) - \gamma j^{(2)}_I -\frac1{\tau_I} j^{(2)}_I \\
\gamma j^{(1)}_R  -\frac1{\tau_R} j^{(1)}_R \\
\gamma j^{(2)}_R  -\frac1{\tau_R} j^{(2)}_R \\ -F(\SO, I_T) \\ F(\SO, I_T) - \gamma \IO \\ \gamma \IO \end{array} \right).
\end{equation*}

As usual for finite volume schemes, data are represented by spatial cell averages, which are defined at time $t^n$ as  
\begin{equation}
\Q_i^n = \frac{1}{|P_i|} \int_{P_i} \Q(\x,t^n) \, d\x,     
\label{eqn.cellaverage}
\end{equation}  
where $|P_i|$ denotes the surface of element $P_i$ at the current time $t^n$. 
Higher order in space is achieved by piecewise high order polynomials. We refer to them as to $\mathbf{w}_i(\x)$ and they have to be reconstructed from the given cell averages \eqref{eqn.cellaverage}. Here, we rely on a second order Central WENO (CWENO) reconstruction procedure along the lines of \cite{ADER-CWENO}. We omit the details for brevity.

\paragraph{Finite volume scheme.} A finite volume method is directly derived by integration of the governing system \eqref{PDE} over a space-time control volume $|P_i|\times [t^n;t^{n+1}]$, thus obtaining
\begin{equation}
\Q_i^{n+1} = \Q_i^n - \frac{\Delta t}{|P_i|}\sum \limits_{P_j \in \mathcal{N}_{S_i}} \,\, {\int \limits_{\lambda_{ij}} \int \limits_{t^n}^{t^{n+1}} 
	\tilde{\mathbf{H}}_{ij}^n \, dt \, d\x}
+ \int\limits_{t^{n}}^{t^{n+1}} \int \limits_{T_i(t)} \S^{n+1} \, d\mathbf{x} \, dt. 
\label{PDEfinal}
\end{equation} 
The term $\tilde{\mathbf{H}}_{ij}^n=\mathbf{F}_{ij} \cdot \mathbf{n}_{ij}$ is a numerical flux function to resolve the discontinuity of the numerical solution at the edges $\lambda_{ij}$ in the normal direction defined by the outward pointing unit normal vector $\mathbf{n}_{ij}$. A simple and robust local Lax-Friedrichs flux is adopted, thus yielding
 \begin{equation}
 \tilde{\mathbf{H}}_{ij}^n =  
 \frac{1}{2} \left( \F(\mathbf{w}_{i,j}^+) + \F(\mathbf{w}_{i,j}^-)  \right) \cdot \mathbf{n}_{ij}  - 
 \frac{1}{2} s_{\max} \left( \mathbf{w}_{i,j}^+ - \mathbf{w}_{i,j}^- \right),  
 \label{eqn.rusanov} 
 \end{equation} 
 where $\mathbf{w}_{i,j}^+,\mathbf{w}_{i,j}^-$ are the high order boundary extrapolated data evaluated through the CWENO reconstruction procedure. The numerical dissipation is given by $s_{\max}$ which is the maximum eigenvalue of the Jacobian matrix in spatial normal direction,  
 \begin{equation}
 \label{eq:An}
 \mathbf{A}_{\mathbf{n}} = \frac{\partial \F}{\partial \Q}.
 \end{equation}
 Notice that for the equations involving the non-commuter populations ($\SO$, $\IO$, $\RO$) in \eqref{eq:diffuse}, the numerical flux must account for a dissipation proportional to the diffusive terms, thus it is supplemented with a numerical viscosity given by the maximum eigenvalue of the viscous operator $s_{max}^V=\max \left( D_S^u, D_I^u, D_R^u \right)$.
 
 Finally, in the diffusion limit the source term $\S(\Q)$ becomes stiff as $(\tau_S,\tau_I,\tau_R)\to 0$, therefore it must be discretized {implicitly} according to \eqref{PDEfinal} in order to overcome too severe time step restrictions. To this aim, a fully second order IMEX method which preserve the asymptotic diffusion limit is proposed and briefly described hereafter.

%------------------------------------------------------------------------------------
%------------------------------------------------------------------------------------
\subsection{Time integration and numerical diffusion limit}
Since the multiscale nature of the dynamics is originated only by the commuters population we restrict our analysis to system \eqref{eq:kineticc} formulated using the parities \eqref{eq:eop1}-\eqref{eq:eop2}. 
For notation simplicity we assume $\tau_{S,I,R}=\tau$ and rewrite \eqref{eq:eop1}-\eqref{eq:eop2} in partitioned form as
\begin{equation}
\begin{split}
	\frac{\partial \u}{\partial t} + \frac{\partial \f(\v)}{\partial x} + \frac{\partial \g(\v)}{\partial y} &= \E(\u) + \frac1{\tau}\left(\U-\u\right)\\
	\frac{\partial \v}{\partial t} + \boldsymbol{\Lambda}^2 \frac{\partial \f(\u)}{\partial x} +\boldsymbol{\Lambda}^2 \frac{\partial \g(\u)}{\partial y} &= \E(\v) - \frac1{\tau}\v,
\end{split}
\label{systcompactform}
\end{equation}
in which 
\begin{equation}
\begin{split}
&\u = \left(r_S^{(1)}, r_S^{(2)}, r_I^{(1)}, r_I^{(2)}, r_R^{(1)}, r_R^{(2)}\right)^T, \quad \v =\left(j_S^{(1)},  j_S^{(2)}, j_I^{(1)}, j_I^{(2)}, j_R^{(1)}, j_R^{(2)}\right)^T,\\ 
&\f(\v) = \xi \v,\quad
\g(\v) = \eta \J \v,\quad \J={\rm diag}\{-1,1,-1,1,-1,1\},\\
&\E(\u)=\left(-F(r_S^{(1)},I_T), -F(r_S^{(2)},I_T),F(r_S^{(1)},I_T)-\gamma r_I^{(1)}, F(r_S^{(2)},I_T)-\gamma r_I^{(2)},\gamma r_I^{(1)},\gamma r_I^{(2)}\right)^T,\\
&\U=\left(S,S,I,I,R,R\right)^T,\quad 
\boldsymbol{\Lambda} ={\rm diag}\{\lambda_S,\lambda_S,\lambda_I,\lambda_I,\lambda_R,\lambda_R\},
\end{split}
\label{eq:variables}
\end{equation}
and $\f(\u)$, $\g(\u)$, $\E(\v)$ are defined similarly.

Following \cite{Bos2}, the Implicit-Explict Runge-Kutta (IMEX-RK) approach that we consider for system \eqref{systcompactform} consists in computing the internal stages
\begin{equation}
\begin{split}
	&\u^{(k)} = \u^n -  \Delta t \sum_{j=1}^{k} a_{kj} \left(\frac{\partial \f(\v^{(j)})}{\partial x} + \frac{\partial \g(\v^{(j)})}{\partial y}-\frac1{\tau}\left(\U^{(j)}-\u^{(j)}\right)\right) + \Delta t \sum_{j=1}^{k-1} \tilde{a}_{kj} \E\left(\u^{(j)}\right)
	\\
	&\v^{(k)} = \v^n -  \Delta t \sum_{j=1}^{k-1} \tilde{a}_{kj} \left(\boldsymbol{\Lambda}^2 \frac{\partial \f(\u^{(j)})}{\partial x} +\boldsymbol{\Lambda}^2 \frac{\partial \g(\u^{(j)})}{\partial y}-\E(\v^{(j)})\right)  + \Delta t \sum_{j=1}^{k} a_{kj} \frac1{\tau}\v^{(j)},
	\end{split}
\label{eq.iterIMEX}
\end{equation}
followed by the numerical solution
\begin{equation}
\begin{split}
	&\u^{n+1} = \U^n -  \Delta t \sum_{k=1}^{s} b_{k} \left(\frac{\partial \f(\v^{(k)})}{\partial x} + \frac{\partial \g(\v^{(k)})}{\partial y}-\frac1{\tau}\left(\U^{(k)}-\u^{(k)}\right)\right) + \Delta t \sum_{k=1}^{s} \tilde{b}_{k} \E\left(\u^{(k)}\right)
	\\
	&\v^{n+1} = \v^n -  \Delta t \sum_{k=1}^{s} \tilde{b}_{k} \left(\boldsymbol{\Lambda}^2 \frac{\partial \f(\U^{(k)})}{\partial x} +\boldsymbol{\Lambda}^2 \frac{\partial \g(\u^{(k)})}{\partial y}-\E(\v^{(k)})\right)  + \Delta t \sum_{k=1}^{s} b_{k} \frac1{\tau}\v^{(k)}.
	\end{split}
\label{eq.finalIMEX}
\end{equation}
Matrices $\tilde A = (\tilde a_{kj})$, with $\tilde a_{kj} = 0 $ for $ j\geq k$, and $A = (a_{kj})$, with $a_{kj} = 0 $ for $ j > k$ are $s \times s$ matrices, with $s$ number of Runge-Kutta stages, defining respectively the explicit and the implicit part of the scheme, and vectors $\tilde b = (\tilde b_1, ...,\tilde b_s)^T$ and $b = (b_1, ...,b_s)^T$ are the quadrature weights. The above IMEX-RK system \eqref{eq.iterIMEX}-\eqref{eq.finalIMEX} is complemented with the same explicit Runge-Kutta scheme defined by $\tilde A$ and $\tilde b$ applied to the non commuters population system \eqref{eq:diffuse}. 

Furthermore, referring to \cite{Bos1, Bos2}, if the following relations hold, 
\begin{equation*}
a_{kj} = b_j, \qquad j = 1,\ldots,s ,\qquad
\tilde a_{kj} = \tilde b_j, \qquad j = 1,\ldots,s-1, 
\end{equation*}
the method is said to be globally stiffly accurate (GSA).
It is worth to notice that this definition states also that the numerical solution of a GSA IMEX-RK scheme coincides exactly with the last internal stage of the scheme. This latter property is fundamental in order to achieve asymptotic-preservation stability in stiff regimes.

\paragraph{Numerical diffusion limit.}
The scheme \eqref{eq.iterIMEX}-\eqref{eq.finalIMEX} permits to treat implicitly the stiff terms and explicitly all the rest, maintaining a consistent discretization of the limit system in the diffusive regime. 
To verify the numerical diffusion limit we assume for simplicity $D_{S,I,R}$ independent from space, the extension to the general case follows straightforwardly.
From the second equation in \eqref{eq.iterIMEX} we have
\[
\tau\v^{(k)} = \tau\v^n -  \Delta t \sum_{j=1}^{k-1} \tilde{a}_{kj} \left(\tau\boldsymbol{\Lambda}^2 \frac{\partial \f(\u^{(j)})}{\partial x} +\tau\boldsymbol{\Lambda}^2 \frac{\partial \g(\u^{(j)})}{\partial y}-\tau \E(\v^{(j)})\right)  + \Delta t \sum_{j=1}^{k} a_{kj} \v^{(j)},
\]
therefore, assuming \eqref{eq:diffcf}, in the limit $\tau\to 0$ yields
\begin{equation}
\sum_{j=1}^{k} a_{kj} \v^{(j)} = \sum_{j=1}^{k-1} \tilde{a}_{kj} \left(2\D \frac{\partial \f(\U^{(j)})}{\partial x} +2\D \frac{\partial \g(\U^{(j)})}{\partial y}\right),
\label{eq:vasym}
\end{equation}
where $\D={\rm diag}\left\{D_S,D_S,D_I,D_I,D_R,D_R\right\}$ and we used the fact that from the first equation in \eqref{eq.iterIMEX} as $\tau\to 0$ we have $\u^{(j)}=\U^{(j)}$. Note that \eqref{eq:vasym} implies that $j_{S,I,R}^{(1)}=j_{S,I,R}^{(2)}$ in $\v^{(j)}$.

Using the identity $\u^{(j)}=\U^{(j)}$ into the first equation in \eqref{eq.iterIMEX} we get 
\begin{equation}
	\U^{(k)} = \U^n -  \Delta t \sum_{j=1}^{k} a_{kj} \left(\frac{\partial \f(\v^{(j)})}{\partial x} + \frac{\partial \g(\v^{(j)})}{\partial y}\right) + \Delta t \sum_{j=1}^{k-1} \tilde{a}_{kj} \E\left(\U^{(j)}\right).
\label{eq:uasym}	
\end{equation}
Applying \eqref{eq:vasym} into \eqref{eq:uasym} thanks to the definitions of $\f$ and $\g$ gives
\begin{equation}
\begin{split}
	\U^{(k)} =& \U^n -  2\Delta t\D \sum_{j=1}^{k-1} \tilde{a}_{kj} \left(\xi^2 \frac{\partial^2 \U^{(j)}}{\partial x^2} + 2\xi\eta\J \frac{\partial^2 \U^{(j)}}{\partial x \partial y} + \eta^2 \frac{\partial^2 \U^{(j)}}{\partial y^2}\right)\\
	&+ \Delta t \sum_{j=1}^{k-1} \tilde{a}_{kj} \E\left(\U^{(j)}\right).
	\end{split}
\label{eq:uasym2}	
\end{equation}

Finally, integrating over $\lambda$ defined by \eqref{eq:lambda} and summing up the components of $\U^{(k)}$ by pairs yields 
\begin{equation}
\begin{split}
	S^{(k)} =& S^n -  \Delta t D_S \sum_{j=1}^{k-1} \tilde{a}_{kj} \left(\frac{\partial^2 S^{(j)}}{\partial x^2} + \frac{\partial^2 S^{(j)}}{\partial y^2}\right)- \Delta t \sum_{j=1}^{k-1} \tilde{a}_{kj} F(S^{(j)},I_T^{(j)}),\\
	I^{(k)} =& I^n -  \Delta t D_I \sum_{j=1}^{k-1} \tilde{a}_{kj} \left(\frac{\partial^2 I^{(j)}}{\partial x^2} + \frac{\partial^2 I^{(j)}}{\partial y^2}\right)+ \Delta t \sum_{j=1}^{k-1} \tilde{a}_{kj} \left(F(S^{(j)},I_T^{(j)})-\gamma I^{(j)}\right),\\
	R^{(k)} =& R^n -  \Delta t D_R \sum_{j=1}^{k-1} \tilde{a}_{kj} \left(\frac{\partial^2 R^{(j)}}{\partial x^2} + \frac{\partial^2 R^{(j)}}{\partial y^2}\right)+ \Delta t \sum_{j=1}^{k-1} \tilde{a}_{kj} \gamma I^{(j)}
	\end{split}
\label{eq:RKdiffuse}	
\end{equation}
and thus, the internal stages correspond to the stages of the explicit scheme applied to the reaction-diffusion system \eqref{eq:diff}. Thanks to the GSA property this is enough to guarantee that the scheme is asymptotic-preserving.

\subsection{Numerical convergence analysis}\label{app4}
\begin{figure}[!htbp]
	\begin{center}
		\begin{tabular}{ccc} 
			\includegraphics[width=0.3\textwidth]{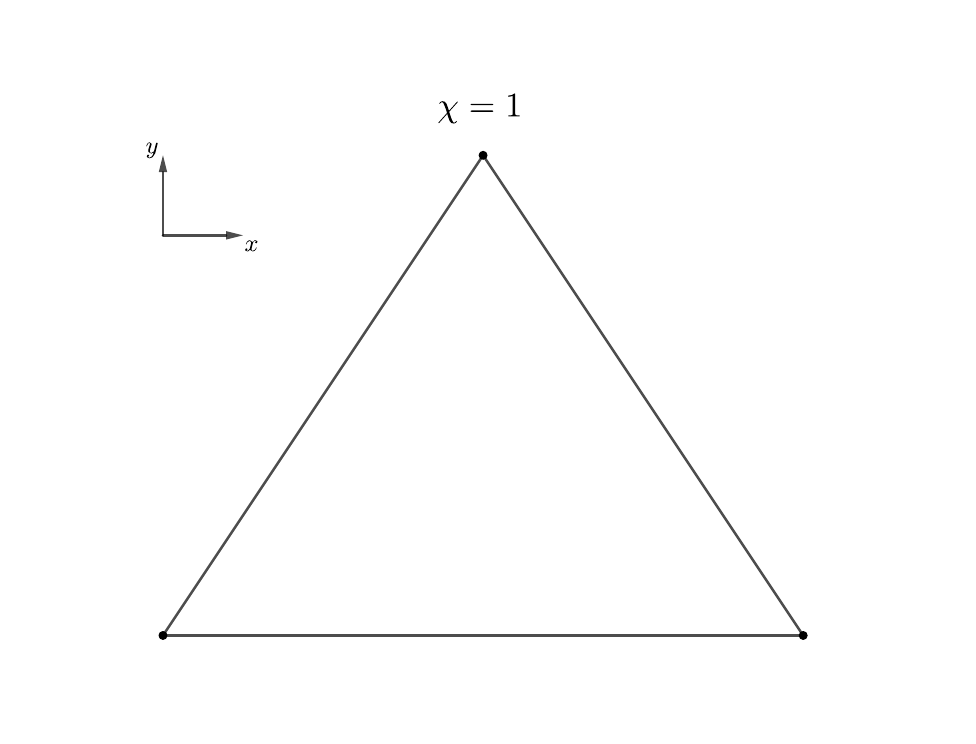} & 
			\includegraphics[width=0.3\textwidth]{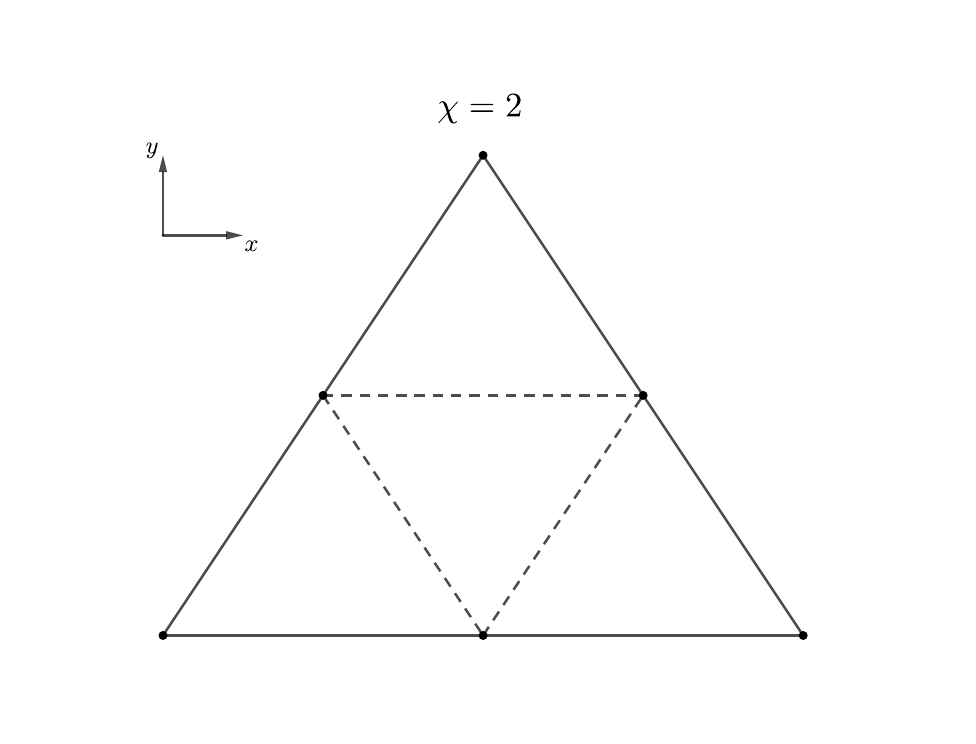} &
			\includegraphics[width=0.3\textwidth]{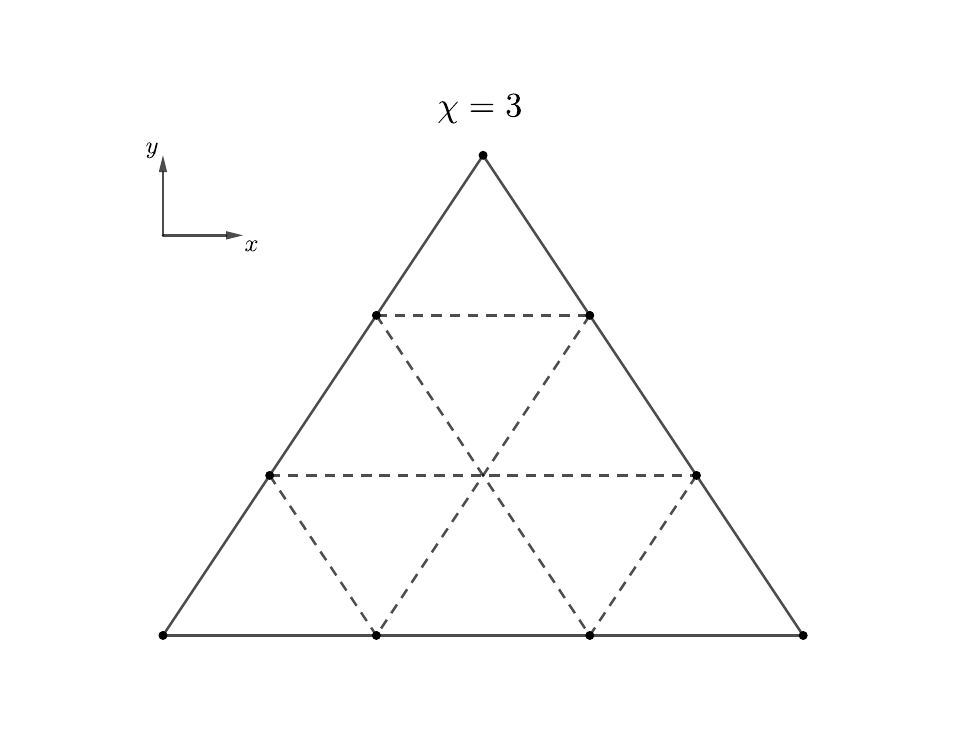} \\   
		\end{tabular} 
	\end{center}
	\caption{Convergence analysis: Isotropic mesh refinement used for convergence analysis with triangular meshes and refinement factor $\chi=1$ ($N_{R}=1$), $\chi=2$ ($N_{R}=4$) and $\chi=3$ ($N_{R}=9$). }
	\label{fig.MeshRef}
\end{figure}

The convergence rate of the novel numerical scheme is studied by considering the MK-SIR model with only commuters and a test problem with no discontinuities neither in the distribution of the populations nor in the associated flux functions. The computational domain is the square $\Omega=[-1;1]^2$ with zero-flux boundary conditions and the initial condition reads
\begin{equation}
	S = \sin(2\pi x) \sin(2\pi y), \qquad I = 1 - S, \qquad R=0,
\end{equation}    
with all other quantities initially set to zero. The contact and recovery rate are set to $\beta=10$ and $\gamma=4$, respectively, while the final time of the simulation is $t_f=0.1$. Three different relaxation times are considered, thus yielding a fully hyperbolic system with $\tau=1$ and $\lambda^2=1$, a mildly diffusive system with $\tau=10^{-2}$ and $\lambda^2=10^2$, and a purely diffusive system  with $\tau=10^{-4}$ and $\lambda^2=10^4$. The computational domain is initially discretized with $N_E=304$ triangles. Mesh refinement is then carried out for triangular meshes relying on conforming finite element discretizations, hence each element is split into sub-elements with an isotropic refinement factor $\chi$. Specifically, a total number of sub-elements $N_{R}=\chi^2$ is generated, as depicted in Figure \ref{fig.MeshRef}. 

\begin{table}[htbp] 
	\caption{Numerical convergence results for the kinetic transport model with second order of accuracy in space and time and discrete number of velocities $M=4$. The errors are measured in $L_1$ norm and refer to the variables $S$, $I$, $j^{(1)}_S$ and $j^{(1)}_I$ for relaxation times $\tau=1.0$, $\tau=10^{-2}$ and $\tau=10^{-4}$.} 
	\begin{center} 	
		\begin{small}
			\renewcommand{\arraystretch}{1.05}
			\begin{tabular}{c|cc|cc|cc|cc} 				
				% tau=1
				\multicolumn{9}{c}{$\tau=1.0$} \\
				& \multicolumn{2}{c|}{$S$} & \multicolumn{2}{c|}{$I$} & \multicolumn{2}{c|}{$j^{(1)}_S$} & \multicolumn{2}{c}{$j^{(1)}_I$} \\
				\hline
				$h(\Omega)$ & $L_1$ & $\mathcal{O}(L_1)$ &  $L_1$ & $\mathcal{O}(L_1)$ & $L_1$ & $\mathcal{O}(L_1)$ & $L_1$ & $\mathcal{O}(L_1)$ \\ 
				\hline
				7.24E-02 & 4.3160E-03 & -    & 2.2020E-02 & -    & 2.4944E-03 & -    & 1.3150E-02 & - \\ 
				4.82E-02 & 1.9773E-03 & 1.93 & 1.0274E-02 & 1.88 & 1.3645E-03 & 1.49 & 7.5176E-03 & 1.38 \\ 
				3.62E-02 & 1.1979E-03 & 1.74 & 6.4905E-03 & 1.60 & 8.0601E-04 & 1.83 & 4.5374E-03 & 1.76 \\ 
				2.89E-02 & 8.0575E-04 & 1.78 & 4.5021E-03 & 1.64 & 4.7722E-04 & 2.35 & 2.6960E-03 & 2.33  \\ 
				\multicolumn{9}{c}{} \\
				% tau=1e-2
				\multicolumn{9}{c}{$\tau=10^{-2}$} \\
				\hline
				$h(\Omega)$ & $L_1$ & $\mathcal{O}(L_1)$ &  $L_1$ & $\mathcal{O}(L_1)$ & $L_1$ & $\mathcal{O}(L_1)$ & $L_1$ & $\mathcal{O}(L_1)$ \\ 
				\hline
				7.24E-02 & 1.1339E-03 & -    & 6.7944E-03 & -    & 2.8801E-03 & -    & 1.7520E-02 & - \\ 
				4.82E-02 & 4.8583E-04 & 2.09 & 2.9178E-03 & 2.08 & 1.3359E-03 & 1.89 & 8.1327E-03 & 1.89 \\ 
				3.62E-02 & 2.4811E-04 & 2.34 & 1.4901E-03 & 2.34 & 7.1628E-04 & 2.17 & 4.3649E-03 & 2.16 \\ 
				2.89E-02 & 1.3812E-04 & 2.62 & 8.2940E-04 & 2.63 & 4.1486E-04 & 2.45 & 2.5317E-03 & 2.44  \\ 
				\multicolumn{9}{c}{} \\
				% tau=1e-4
				\multicolumn{9}{c}{$\tau=10^{-4}$} \\
				\hline
				$h(\Omega)$ & $L_1$ & $\mathcal{O}(L_1)$ &  $L_1$ & $\mathcal{O}(L_1)$ & $L_1$ & $\mathcal{O}(L_1)$ & $L_1$ & $\mathcal{O}(L_1)$ \\ 
				\hline
				7.24E-02 & 2.9065E-03 & -    & 1.6949E-02 & -    & 6.4945E-03 & -    & 3.9721E-02 & - \\ 
				4.82E-02 & 1.4582E-03 & 1.70 & 8.7873E-03 & 1.62 & 3.3835E-03 & 1.61 & 2.0671E-02 & 1.61 \\ 
				3.62E-02 & 8.0997E-04 & 2.04 & 4.9511E-03 & 1.99 & 1.9730E-03 & 1.87 & 1.2093E-02 & 1.86 \\ 
				2.89E-02 & 4.7890E-04 & 2.36 & 2.9495E-03 & 2.32 & 1.2228E-03 & 2.14 & 7.5199E-03 & 2.13  \\ 
			\end{tabular}
		\end{small}
	\end{center}
	\label{tab.conv_n4}
\end{table}

Errors are measured in $L_1$ norm as
\begin{equation}
	L_1 = \int_{\Omega}\left\|s_{e}(\x)-s_{h}(\x)\right\| dA,
\end{equation}
with $s_{e}$ and $s_{h}$ denoting the reference and the numerical solution of a generic variable $s$ of the system. Results are shown for variables $S$, $I$, $j^{(1)}_S$ and $j^{(1)}_I$, demonstrating that the formal order of accuracy is achieved in all regimes thanks to the asymptotic preserving property exhibited by the second order IMEX scheme.  

\begin{table}[htbp] 
	\caption{Numerical convergence results for the kinetic transport model with second order of accuracy in space and time and discrete number of velocities $M=8$. The errors are measured in $L_1$ norm and refer to the variables $S$, $I$, $j^{(1)}_S$ and $j^{(1)}_I$ for relaxation times $\tau=1.0$, $\tau=10^{-2}$ and $\tau=10^{-4}$.} 
	\begin{center} 	
		\begin{small}
			\renewcommand{\arraystretch}{1.05}
			\begin{tabular}{c|cc|cc|cc|cc} 				
				% tau=1
				\multicolumn{9}{c}{$\tau=1.0$} \\
				& \multicolumn{2}{c|}{$S$} & \multicolumn{2}{c|}{$I$} & \multicolumn{2}{c|}{$j^{(1)}_S$} & \multicolumn{2}{c}{$j^{(1)}_I$} \\
				\hline
				$h(\Omega)$ & $L_1$ & $\mathcal{O}(L_1)$ &  $L_1$ & $\mathcal{O}(L_1)$ & $L_1$ & $\mathcal{O}(L_1)$ & $L_1$ & $\mathcal{O}(L_1)$ \\ 
				\hline
				7.24E-02 & 3.9041E-03 & -    & 1.9272E-02 & -    & 1.1499E-03 & -    & 5.8528E-03 & - \\ 
				4.82E-02 & 1.7662E-03 & 1.96 & 8.8050E-03 & 1.93 & 5.8676E-04 & 1.66 & 3.1030E-03 & 1.48 \\ 
				3.62E-02 & 1.0044E-03 & 1.96 & 5.1038E-03 & 1.90 & 3.3150E-04 & 1.98 & 1.7731E-03 & 1.95 \\ 
				2.89E-02 & 6.0960E-04 & 2.24 & 3.1314E-03 & 2.19 & 2.1407E-04 & 1.96 & 1.1824E-03 & 1.82  \\ 
				\multicolumn{9}{c}{} \\
				% tau=1e-2
				\multicolumn{9}{c}{$\tau=10^{-2}$} \\
				\hline
				$h(\Omega)$ & $L_1$ & $\mathcal{O}(L_1)$ &  $L_1$ & $\mathcal{O}(L_1)$ & $L_1$ & $\mathcal{O}(L_1)$ & $L_1$ & $\mathcal{O}(L_1)$ \\ 
				\hline
				7.24E-02 & 8.9391E-04 & -    & 5.2662E-03 & -    & 1.3332E-03 & -    & 8.2209E-03 & - \\ 
				4.82E-02 & 4.6649E-04 & 1.60 & 2.8084E-03 & 1.55 & 7.6802E-04 & 1.36 & 4.7691E-03 & 1.34 \\ 
				3.62E-02 & 2.7284E-04 & 1.86 & 1.6548E-03 & 1.84 & 4.7508E-04 & 1.67 & 2.9561E-03 & 1.66 \\ 
				2.89E-02 & 1.6479E-04 & 2.26 & 1.0015E-03 & 2.25 & 2.9910E-04 & 2.07 & 1.8616E-03 & 2.07  \\ 
				\multicolumn{9}{c}{} \\
				% tau=1e-4
				\multicolumn{9}{c}{$\tau=10^{-4}$} \\
				\hline
				$h(\Omega)$ & $L_1$ & $\mathcal{O}(L_1)$ &  $L_1$ & $\mathcal{O}(L_1)$ & $L_1$ & $\mathcal{O}(L_1)$ & $L_1$ & $\mathcal{O}(L_1)$ \\ 
				\hline
				7.24E-02 & 2.8333E-03 & -    & 1.6481E-02 & -    & 3.1802E-03 & -    & 1.9511E-02 & - \\ 
				4.82E-02 & 1.4078E-03 & 1.72 & 8.4724E-03 & 1.64 & 1.6563E-03 & 1.61 & 1.0163E-02 & 1.61 \\ 
				3.62E-02 & 7.7612E-04 & 2.07 & 4.7389E-03 & 2.02 & 9.5751E-04 & 1.90 & 5.8910E-03 & 1.90 \\ 
				2.89E-02 & 4.5623E-04 & 2.38 & 2.8063E-03 & 2.35 & 5.8739E-04 & 2.19 & 3.6247E-03 & 2.18  \\  
			\end{tabular}
		\end{small}
	\end{center}
	\label{tab.conv_n8}
\end{table}


\begin{thebibliography}{99.}

\bibitem{Albietal} G.~Albi, N.~Bellomo, L.~Fermo, S.-Y.~Ha, J.~Kim, L.~Pareschi, D.~Poyato, J.~Soler.  Vehicular traffic, crowds, and swarms: from kinetic theory and multiscale methods to applications and research perspectives.
 Math. Models Methods Appl. Sci.  29  (2019),  no. 10, 1901--2005.

\bibitem{APZ} G.~Albi, L.~Pareschi, M.~Zanella.
\newblock Control with uncertain data of socially structured compartmental epidemic models.
\newblock preprint arXiv:2004.13067, 2020 

\bibitem{ABLN} L.J.S.~Allen, B.M.~Bolker, Y.~Lou, A.L.~Nevai.
\newblock Asymptotic profiles of the steady states for an SIS epidemic reaction–diffusion model.
\newblock \emph{Discrete Contin. Dyn. Syst.} 21, 1--20, (2008).

\bibitem{BGRCV} D.~Balcan, B.~Gonçalves, H.~Hu, J.J.~Ramasco, V.~Colizza, A.~Vespignani. Modeling the spatial spread of infectious diseases: the GLobal Epidemic and Mobility computational model. 
\newblock \emph{J. Comput. Sci.} 1(3):132--145, (2010).

\bibitem{BCV13} E.~Barbera, G.~Consolo, G.~Valenti. 
Spread of infectious diseases in a hyperbolic reaction-diffusion susceptible-infected-recovered model. {\em Physical Review E}, 88, 052719 (2013).

%\bibitem{barthlsq} T.J.~Barth, P.O.~Frederickson. Higher order solution of the {Euler} equations on unstructured grids using quadratic reconstruction. \emph{AIAA paper no. 90-0013}, (1990).

\bibitem{bellomo2020multiscale}
N.~Bellomo, R.~Bingham, M.~A.~J. Chaplain, G.~Dosi, G.~Forni, D.~A. Knopoff,
  J.~Lowengrub, R.~Twarock, and M.~E. Virgillito.
\newblock A multi-scale model of virus pandemic: Heterogeneous interactive
  entities in a globally connected world.
\newblock {\em Mathematical Models and Methods in Applied Sciences}, 30(8):1591--1651, (2020)
  2020.
  
\bibitem{Bellomo2000a}  
N.~Bellomo, S.~St\"ocker, Sabine. Development of Boltzmann models in mathematical biology.
\emph{Modeling in applied sciences}, Model. Simul. Sci. Eng. Technol., Birkh\"auser Boston, 225--262, (2000). 


\bibitem{Bert} G.~Bertaglia, L.~Pareschi. Hyperbolic models for the spread of epidemics on networks: kinetic description and numerical methods, \emph{ESAIM Math. Model. Numer. Anal.} to appear, 2020.

\bibitem{Bos1} S.~Boscarino, L.~Pareschi, G.~Russo. Implicit-Explicit Runge-Kutta schemes for hyperbolic systems and kinetic equations in the diffusion limit. \emph{SIAM J. Sci. Comp.} 35:22--51, (2013).

\bibitem{Bos2} S.~Boscarino, L.~Pareschi, G.~Russo. A unified IMEX Runge-Kutta approach for hyperbolic systems with   multiscale relaxation. \emph{SIAM J. Numer. Anal.} 55(4):2085--2109, (2017).

\bibitem{BD-BGK} W.~Boscheri, G.~Dimarco. High order central WENO-Implicit-Explicit Runge Kutta
schemes for the BGK model on general polygonal meshes. \emph{Journal of Computational Physics}, 422:109766 (2020).


\bibitem{cs} K.M.~Case, P.F.~Zweifel. Existence and Uniqueness Theorems for the
Neutron Transport Equation. \emph{J. Math. Physics} 4(11):1376--1385, (1963).
	
\bibitem{Cap} V.~Capasso.
\newblock Global solution for a diffusive nonlinear deterministic epidemic model.
\newblock \emph{SIAM J. Appl. Math.} 35:274--284, (1978).  

\bibitem{CS78} V.~Capasso, G.~Serio, A generalization of the Kermack-McKendrick deterministic epidemic model. {\em Math. Biosci.} 42, 43 (1978).

\bibitem{Cer} C.~Cercignani, R.~Illner, M.~Pulvirenti. 
\newblock \emph{The Mathematical Theory of Diluted Gases}
\newblock Springer, New York, (1994).

\bibitem{CMPS} F.A.C.C.~Chalub, P.A.~Markovich, B.~Perthame, C.~Schmeiser.
\newblock Kinetic models for chemotaxis and their drift-diffusion limits. 
\newblock \emph{Monatsh. Math.} 142, 123--141 (2004).



\bibitem{CV} V.~Colizza, A.~Vespignani.
\newblock A. Epidemic modeling in metapopulation systems with heterogeneous coupling
pattern: Theory and simulations. 
\newblock \emph{J. Theor. Biol.} 251:450--467, (2008).

\bibitem{Colombo} R.M.~Colombo, M.~Garavello, F.~Marcellini, E.~Rossi. An age and space structured SIR model describing the COVID-19 pandemic. \emph{J. Math. Ind.}, 10(1): 22, 2020.

%\bibitem{CLL} R.~Cui, K-Y.~Lam, Y.~Lou.
%\newblock Dynamics and asymptotic profiles of steady states of an epidemic model in advective environments.
%\newblock \emph{Journal of Differential Equations}, 263(4):2343--2373, (2017).


\bibitem{Deli} 
M.~Delitala.
Generalized kinetic theory approach to modeling spread- and evolution of epidemics.
{\em Mathematical and Computer Modelling}, 39(1):1--12, (2004).

\bibitem{Diekmann} 
O.~Diekmann, J.~Heesterbeek, M.~Roberts, The construction of next-generation matrices for compartmental epidemic models. \emph{J. Roy. Soc. Interface}, 7:873--885, (2010).

\bibitem{DP} G.~Dimarco, L.~Pareschi.  Numerical methods for kinetic equations. \emph{Acta Numer.}  23:369--520, (2014).

\bibitem{DPTZ} G.~Dimarco, L.~Pareschi, G.~Toscani, M.~Zanella. Wealth distribution under the spread of infectious diseases, Phys. Rev. E 102, 022303, (2020).

\bibitem{ADER-CWENO} M.~Dumbser, W.~Boscheri, M.~Semplice, G.~Russo. Central weighted ENO schemes for hyperbolic conservation laws on fixed and moving unstructured meshes. \emph{SIAM Journal of Scientific Computing}, 39: A2564-A2591 (2017).

\bibitem{Dumbser2007693} M.~Dumbser, M.~Kaeser. Arbitrary high order non-oscillatory finite volume schemes on unstructured meshes for linear hyperbolic systems. \emph{Journal of Computational Physics}, 221:693--723 (2007).


\bibitem{FMW} W.~E.~Fitzgibbon, J.~J.~Morgan, G.~F.~Webb.
\newblock An outbreak vector-host epidemic model with spatial structure: the 2015–2016 Zika outbreak in Rio De Janeiro.
\newblock \emph{Theor. Biol. Med. Model.} 14: 7, (2017).

  \bibitem{Franco2020}
E.~Franco.
\newblock A feedback {SIR (fSIR)} model highlights advantages and limitations
  of infection-based social distancing.
\newblock {\em arXiv:2004.13216}, 2020.

\bibitem{FWF} E.~Frias-Martinez, G.~Williamson, V.~Frias-Martinez. 
\newblock An Agent-Based Model of Epidemic Spread using
Human Mobility and Social Network Information. 
\newblock \emph{In Proceedings of the 3rd International Conference on
Social Computing} (SocialCom11), Boston, MA, USA, 49--56, (2011).

\bibitem{ArepoTN} E.~Gaburro, W.~Boscheri, S.~Chiocchetti, C.~Klingenberg, V. Springel, M.~Dumbser. High order direct Arbitrary-Lagrangian-Eulerian schemes on moving Voronoi meshes with topology changes. \emph{Journal of Computational Physics}, 407:109167 (2020).

\bibitem{Gatto} M.~Gatto, E.~Bertuzzo, L.~Mari, S.~Miccoli, L.~Carraro, R.~Casagrandi, A.~Rinaldo. Spread and dynamics of the COVID-19 epidemic in Italy: Effects of emergency containment measures.
\emph{Proceedings of the National Academy of Sciences}, 117(19):10484--10491, (2020).

%\bibitem{FL71} K.O.~Friedrichs, P. D. Lax. Systems of Conservation Equations with a Convex Extension, {\em Proceedings of the National Academy of Sciences} 68(8):1686--1688, (1971).

\bibitem{GJL} F.~Golse, S.~Jin, C.~Levermore. The convergence of numerical transfer schemes in diffusive regimes I: Discrete-ordinate method. \emph{SIAM Journal on Numerical analysis}, 36(5):1333--1369, (1999).

\bibitem{HS} T.~Hillen, A.~Swan. 
\newblock The diffusion limit of transport equations in biology. 
\newblock In: Preziosi~L., Chaplain~M., Pugliese~A. (eds) Mathematical Models and Methods for Living Systems. 
\newblock \emph{Lecture Notes in Mathematics} 2167, Springer, Cham, (2016). 

\bibitem{HWH00} H.W.~Hethcote, The Mathematics of Infectious Diseases. {\em SIAM Review} 42(4):599--653, (2000).

\bibitem{JHL} S.~Jin, H.~Lu, L.~Pareschi. Efficient stochastic asymptotic-preserving implicit-explicit methods for transport equations with diffusive scalings and random inputs.
 \emph{SIAM J. Sci. Comput.}  40(2):A671--A696, (2018).

\bibitem{JPT} S.~Jin, L.~Pareschi, G.~Toscani. Uniformly accurate diffusive relaxation schemes for multiscale transport equations. \emph{SIAM Journal on Numerical Analysis}, 38(3):913--936, (2000).

\bibitem{KGV} J.P.~Keller, L.~Gerardo-Giorda, A.~Veneziani. Numerical simulation of a susceptible-exposed-infectious space-continuous model for the spread of rabies in raccoons across a realistic landscape. \emph{J. Biol. Dyn.}, 7(1):31--46, (2014).

\bibitem{KM05} A. Korobeinikov,  P. K. Maini. Non-linear incidence and stability of infectious disease models. {\em Mathematical Medicine and Biology: A Journal of the IMA},  22, 113--128, (2005). 

\bibitem{LK} E.W.~Larsen, J.B.~Keller. Asymptotic solution of neutron transport problems for small free mean paths. 
\emph{J. Math. Phys.}, 15:75--81, (1974).



%\bibitem{LPW} H. Li, R. Peng, F-B. Wang.
%\newblock Varying total population enhances disease persistence: Qualitative analysis on %a diffusive SIS epidemic model.
%\newblock \emph{Journal of Differential Equations}, 262(2):885--913, (2017).

\bibitem{Liu} Q-X.~Liu, Z.~Jin. Formation of spatial patterns in an epidemic model with constant removal rate of the infectives.
\emph{Journal of Statistical Mechanics: Theory and Experiment}, 2007(05):P05002--P05002,(2007).

\bibitem{MWW} P.~Magal, G.F.~Webb, X.~Wu. 
\newblock Spatial spread of epidemic diseases in geographical settings: Seasonal influenza epidemics in Puerto Rico.
\newblock \emph{Discrete \& Continuous Dynamical Systems - B},25(6):2185--2202, (2019).


%\bibitem{MR98} I.~Muller, T.~Ruggeri. {\em Rational Extended Thermodynamics}, Springer, New York, (1998).

\bibitem{Per} B.~Perthame.
\newblock \emph{Transport Equations in Biology}.
\newblock Birkh\"auser, Boston, (2007)

\bibitem{PuSa} M.~Pulvirenti, S.~Simonella.
A kinetic model for epidemic spread. \emph{Math. Mech. Complex Systems}, 8(3):249--260, 2020.

\bibitem{FBK2020} F.~Riccardo, M.~Ajelli et al. Epidemiological characteristics of COVID-19 cases in Italy and estimates of the reproductive numbers one month into the epidemic. \emph{MedRxiv, Cold Spring Harbor Laboratory Press}, doi:10.1101/2020.04.08.20056861 (2020).

\bibitem{REIM} S.~Riley, K.~Eames, V.~Isham, D.~Mollison, P.~Trapman.
\newblock Five challenges for spatial epidemic models.
\newblock \emph{Epidemics}, 10:68--71, (2015).

\bibitem{SaSiLu} Md.~Samsuzzoha, M.~Singh, D.~Lucy.
Numerical study of a diffusive epidemic model of influenza with variable transmission coefficient. \emph{
Applied Mathematical Modelling},
35(2):5507--5523, (2011).

\bibitem{Sun} G.~Sun. Pattern formation of an epidemic model with diffusion. \emph{Nonlinear Dyn.} 69:1097--1104, (2012). 


\bibitem{Tang} B.~Tang, X.~Wang, A.~Li, N.L.~Bragazzi, S.~Tang, Y.~Xiao, J.~Wu. Estimation of the transmission risk of the 2019-nCoV and its implication for public health interventions. \emph{J. Clinical Med.} 9, 462 (2020).

\bibitem{Veneziani2021}  A.~Viguerie, G.~Lorenzo, F.~Auricchio, D.~Baroli, T.J.R.~Hughes, A.~Patton, A.~Reali, T.E.~Yankeelov, A.~Veneziani. Simulating the spread of COVID-19 via a spatially-resolved susceptible–exposed–infected–recovered–deceased (SEIRD) model
with heterogeneous diffusion. \emph{Applied Mathematics Letters}, 101:106617, (2021).

\bibitem{Veneziani2020}  A.~Viguerie, A.~Veneziani, G.~Lorenzo,  D.~Baroli, N. Aretz-Nellesen, A.~Patton,  T.E.~Yankeelov, A.~Reali,T.J.R.~Hughes, F.~Auricchio. Diffusion–reaction compartmental models formulated in a continuum mechanics framework: application to COVID-19, mathematical analysis, and numerical study. \emph{Comput Mech}, 66, 1131--1152, (2020).



\bibitem{Webb} G.F.~Webb. A reaction-diffusion model for a deterministic diffusion epidemic. \emph{J. Math. Anal. Appl.}, 84: 150--161, (1981). 

\bibitem{Wang2020}  J.~Wang, F.~Xie, T.~Kuniya. Analysis of a reaction-diffusion cholera epidemic model in a spatially heterogeneous environment. \emph{Communications in Nonlinear Science and Numerical Simulation}, 80:104951, (2020).

\bibitem{RS}
R.~Yano. Kinetic modeling of local epidemic spread and its simulation.
\emph{J. Sci. Comput.}, 73:122--156, (2017). 


\bibitem{istat} \url{https://www4.istat.it/it/archivio/209722}

\bibitem{prot_civile} \url{https://github.com/pcm-dpc/COVID-19}

\bibitem{mobilityER} \url{https://sasweb.regione.emilia-romagna.it/statistica/SceltaAnno.do?analisi=matPend2011_2015}






  




\end{thebibliography}
\end{document}